\documentclass[11pt]{article}

\usepackage{amsmath,amsthm,epsfig, amscd}

\def\ftoday{le \space\number\day \space\ifcase\month\or

  janvier\or f\'evrier\or mars\or avril\or mai\or juin\or

  juillet\or ao\^ut\or septembre\or octobre\or novembre\or d\'ecembre\fi

  \space\number\year}

\def\real{I\kern-0.20em R}

\def\integer{I\kern-0.20em N}

\def\relative{{\rm \rlap Z\kern 2.2pt Z}}

\def\cc{\kern-.25em{\c c}}

\def\bc{\begin{center}}

\def\ec{\end{center}}

\def\=def{\stackrel{{\rm def}}{=}}

\newcommand\vfo[2]{{\cal X}_{#1}^{#2}}

\newcommand\fvfo[2]{\widehat{\cal X}_{#1}^{#2}}

\newcommand\pvf[3]{{\cal P}_{#1}^{#2,#3}}

\newcommand\hvf[2]{{\cal P}_{#1}^{#2}}

\newcommand\lie[1]{{\frak #1}}

\renewcommand\hom[1]{\text{Hom}_{\Bbb #1}}

\newcounter{indconst}

\newcounter{auxconst}

\def\bit{\begin{itemize}}

\def\eit{\end{itemize}}

\def\ben{\begin{enumerate}}

\def\een{\end{enumerate}}

\def\bde{\begin{description}}

\def\ede{\end{description}}

\def\beq{\begin{equation}}

\def\eeq{\end{equation}}

\def\bfi{\begin{figure}[hbt] \begin{center}}

\def\efi{\end{center} \end{figure}}

\def\bce{\begin{center}}

\def\ece{\end{center}}

\newtheorem {theos} {Theorem}[subsection]

\newtheorem {lemms} [theos]{Lemma}
\newtheorem {props}[theos] {Proposition}
\newtheorem {defis}[theos] {Definition}
\newtheorem {rems}[theos]{Remark}
\newtheorem {ex}[theos] {Exemple}

\input amssym.def
\input amssym

\numberwithin{equation}{subsection}

\begin{document}
\title{A KAM phenomenon for singular holomorphic vector fields}
\author{Laurent Stolovitch \thanks{CNRS UMR 5580, Laboratoire Emile Picard,
Universite Paul Sabatier, 118 route de Narbonne,
31062 Toulouse cedex 4, France. E-mail : {\tt stolo@picard.ups-tlse.fr}}}

\maketitle



\begin{abstract}
Let $X$ be a germ of holomorphic vector field at the origin of ${\bf C}^n$ and vanishing there. We assume that $X$ is a  ``nondegenerate" good perturbation of a singular completely
integrable system. The latter is associated to a family of linear diagonal vector fields which is assumed to have nontrivial polynomial first integrals. We show that $X$ admits many invariant analytic subsets in a neighborhood of the origin. These are biholomorphic to the intersection of a
polydisc with an analytic set of the form ``resonant monomials =
constants". Such a biholomorphism conjugates the restriction of $X$ to one of its invariant varieties to the restriction of a linear diagonal vector field to a toric variety. Moreover, we show that the set of ``frequencies" defining the invariant sets is of positive measure.
\end{abstract}
\tableofcontents
\section{Introduction}

In this article, we continue our earlier work on germs of singular holomorphic vector fields in ${\bf C}^n$. Our aim is to give a better understanding of the behavior of complex flows in a neighborhood of an isolated singular point (which will be $0$) of such a  vector field.

As it is well known, the behavior of trajectories at the vicinity of the singular point is very difficult to describe.
These difficulties are closely related, in the one hand, to the problem of {\it small divisors} and, on the other hand, to the problem of {\it symmetries} and {\it first integrals}.

The vector fields for which the situation is well understood are the completely integrable ones, in a sense which will be recalled later on.
One of their main features is that they are holomorphically normalizable in a
neighborhood of the origin. The analysis of the behavior of the flow can
therefore be carried out on the normal form and then, pulled back by the biholomorphism: all the fibers of an associated algebraic map (the moment map or the resonant map), when intersected with a fixed polydisc around the origin, are invariant by the flow of the normal form. Moreover, its restriction is nothing but the restriction of a linear diagonal vector field whose eigenvalues depend only on the fiber.

The situation we shall deal with concerns the perturbed case. 
By this we mean the following: we choose a holomorphic singular completely integrable system. Let us perturb it in some way. What can be said about the behavior of the flow of the perturbed system~?
Roughly speaking, we shall show that, generically, a large set of the deformed fibers is still invariant under the perturbed flow.

\subsection{Classical hamiltonian framework: Complete integrability and KAM theory}

First of all, let us recall Liouville's theorem \cite{Arn1} which concerns hamiltonian systems. 
Let $H_1,\ldots, H_n$ be smooth functions on a smooth symplectic manifold $M^{2n}$; 
let $\pi:M^{2n}\rightarrow {\bf R}^n$ be the {\it moment map} defined to be $\pi(x)=(H_1(x),\ldots,H_n(x))$. We assume that, for all $1\leq i,j\leq n$, 
the Poisson brackets $\{H_i,H_j\}$ vanish. Let $c\in {\bf R}^n$ be a regular value of $\pi$; we assume that 
$\pi^{-1}(c)$ is compact and connected. 
Then there exists a neighborhood $U$ of $\pi^{-1}(c)$ and a symplectomorphism $\Phi$ from $U$ to 
$\pi(U)\times\Bbb T^n$ such that, in this new coordinate system, the symplectic
vector field $X_{H_i}$ associated to each $H_i$ 
is tangent to the fiber $\{d\}\times \Bbb T^n$. It is constant on it and the constant depends only on the fiber. They define quasi-periodic motions on each torus.
The family of hamiltonian vector fields $X_{H_1},\ldots,X_{H_n}$ is said to be {\bf completely integrable}.

Nevertheless, completely integrable systems are pretty rare when one looks at problems arising from physics and in particular, celestial mechanics. One often encounters small perturbations of integrable systems. So, the natural question to be asked is: what can be said about these nonintegrable systems~? Do these systems still have invariant tori on which the motion is quasi-periodic~? 
Of course, the perturbation is assumed to be hamiltonian.

The answer was given almost fifty years ago by the celebrated KAM theorem. It
is named after its authors Kolmogorov-Arnold-Moser \cite{Kolmogorov,Kolmogorov-congres,arnold-kol,arnold-celest,moser-anneau}. Roughly
speaking, this theorem states that, if the integrable vector field which is to
be perturbed is nondegenerate in some sense and if the perturbation is small
enough and still hamiltonian then 
there is a set of ``positive measure" of invariant tori for the perturbed
hamiltonian and it gives rise to quasi-periodic motions of these tori. The
constants defining the quasi-periodic motions of the tori satisfy some
diophantine condition. 

Let $(\theta, I)$ be symplectic coordinates (angles-actions) of $\Bbb T^n\times
{\bf R}^n$ ($\Bbb T^n$ denotes the $n$-dimensional torus). Assume that the flow
of the unperturbed hamiltonian $H_0$
\begin{displaymath}
\begin{cases} \dot\theta = \omega(I)\\ \dot I = 0\end{cases}
\end{displaymath}
with $\omega(I)$ belongs to ${\bf R}^n$ and is such that $\det(\frac{\partial
  \omega_i}{\partial I_j})(0)\neq $ (this is the classical nondegeneracy
  condition). Let us consider a small analytic perturbation of $H_0$:
\begin{displaymath}
\begin{cases} \dot\theta = \omega(I)+\epsilon f(\theta,I, \epsilon)\\ \dot I = \epsilon g(\theta,I, \epsilon)\end{cases}.
\end{displaymath}
According to the nondegeneracy condition, for any $k\geq 1$, there is an
analytic change of coordinates $(\phi,J)$ such that
\begin{displaymath}
\begin{cases} \dot\phi = \omega_k(J)+\epsilon^k f_k(\phi,J, \epsilon)\\ \dot J = a_k(J,\epsilon)+\epsilon^k g_k(\phi,J, \epsilon)\end{cases}.
\end{displaymath}
It is defined on some open set in the $J$ coordinates.
This is known as the Lindstedt procedure. We shall call this a {\bf Lindstedt
normal form} up to order $k$. One can get rid of the fast
variables (angles) up to any order of the perturbation. Moreover, if we assume
that the perturbation is still hamiltonian then we have 
\begin{displaymath}
\begin{cases} \dot\phi = \omega_k(J)+\epsilon^k f_k(\phi,J, \epsilon)\\ \dot J = \epsilon^k g_k(\phi,J, \epsilon)\end{cases}.
\end{displaymath}
The KAM procedure says that Lindstedt normalization process can be carried out
``until the end" if the slow variable $J$ belongs to some well chosen set: there is a symplectic change of coordinates such
that, if $J_0\in {\bf R}^n$ belongs to this set, we have
\begin{displaymath}
\begin{cases} \dot\phi = \omega_{\infty}(J_0)\\ \dot J = 0\end{cases}.
\end{displaymath}
This shows that the torus $\Bbb T^n\times \{J_0\}$ is an invariant manifold in
the new coordinates. We refer to \cite{arnold-ds3}[chapter 5].
Moreover, it is required that $\omega_{\infty}(J_0)$ be diophantine, that is 
$$
\exists C, \nu >0,\; \forall Q\in {\bf Z}^n\setminus \{0\}, \;|(Q,\omega_{\infty}(J_0))|> \frac{C}{|Q|^{\nu}},
$$ 
where $(.,.)$ denotes the usual scalar product of ${\bf R}^n$ and $|Q|$ denotes
the sum of the absolute values of the coordinates of $Q$. 

Both the nondegeneracy condition and the diophantine condition have been improved by H. R\"{u}ssmann \cite{russmann-weak}.

A very nice introduction to these results can be found in the exposition at
S\'eminaire Bourbaki of J.-B. Bost \cite{bost-bourb}; it contains a proof of
the KAM theorem based on the Nash-Moser theorem (see also \cite{zehnder1,zehnder2,eliasson-kam}). Other surveys on that topic are \cite{arnold-ds3},
\cite{Arn1}[appendix 8] and in particular, the book \cite{bhs-book} by Broer,
Huitema and Sevryuk, which contains an extensive bibliography. About Lindsedt
expansion, one can consult the article \cite{eliasson-lindstedt}.

Since then, a lot of work has been done on that subject. A closely related theme is the existence of invariant circles of twist mappings of the annulus \cite{Russmann-anneau1,Russmann-anneau2,Herman-anneau1,Herman-anneau2,Yoccoz-bourb-herman} as well as the bifurcation of elliptic fixed point (smooth case) \cite{Chenciner1,Yoccoz-Bourb-Chenciner}.
These topics together with their links with celestial mechanics are explained in the books by C.L. Siegel and J.~Moser \cite{siegel-moser}, by S. Sternberg \cite{Sternb1,Sternb2} and by J.~Moser \cite{moser-book}.

All this literature is concerned with nonsingular hamiltonian dynamical systems. Few results have been obtained in the singular case \cite{arnold-equil}. 

\subsection{Singular complete integrability}

From now on, we shall be concerned with singular holomorphic vector fields in
a neighborhood of the origin of ${\bf C}^n$, $n\geq 2$. Let us recall one of the statements of a previous article \cite{stolo-ihes}.

Let $\lie g$ be a $l$-dimensional commutative Lie algebra over ${\bf C}$. Let $\lambda_1,\ldots,\lambda_n$ 
be complex linear forms over $\lie g$ such that the Lie morphism $S$ from $\lie g$ to the Lie algebra 
of linear vector fields of ${\bf C}^n$ defined by $S(g)=\sum_{i=1}^n\lambda_i(g)x_i\partial/\partial x_i$ 
is injective. 
For any $Q=(q_1,\ldots,q_n)\in {\bf N}^n$ and $1\leq i\leq n$, we define the weight $\alpha_{Q,i}(S)$ of $S$ to be the linear form 
$\sum_{j=1}^n{q_j\lambda_j(g)}-\lambda_i(g)$. Let us set
$|Q|=q_1+\cdots+q_n$. Let $\|.\|$ be a norm on the ${\bf C}$-vector space 
of linear forms on $\lie g$.
Let us define a sequence of positive real numbers 
$$
\omega_k=\inf\left\{\|\alpha_{Q,i}\|\neq 0, 1\leq i\leq n, 2\leq |Q|\leq 2^k\right\}. 
$$
We define a {\bf diophantine condition} relative to $S$ to be 
$$
(\omega(S))\quad\quad -\sum_{k\geq 0}\frac{\ln \omega_k}{2^k}<+\infty.
$$
Let $\vfo n k$ (resp. $\fvfo n k$) be the Lie algebra of germs of holomorphic (resp. formal) vector fields vanishing at order greater than or equal to $k$ at 
$0\in {\bf C}^n$. Let $\left(\fvfo n 1\right)^S$ (resp. $\widehat{\cal O}_n^S$) 
be the formal centralizer of $S$ (resp. the ring of formal first integrals), that is the set of formal vector fields $X$ 
(resp. formal power series $f$) such that $[S(g),X]=0$ (resp. ${\cal L}_{S(g)}(f)=0$) for all $g\in \lie g$. 

A nonlinear deformation $S+\epsilon$ of $S$ is a Lie morphism from $S$ to $\vfo n 1$ such that $\epsilon\in \hom C(\lie g, \vfo n 2)$. 
Let $\hat\Phi$ be a formal diffeomorphism of $({\bf C}^n,0)$ which is assumed to be tangent to $Id$ at $0$. We define 
$\hat \Phi^*(S+\epsilon)(g):=\hat\Phi^*(S(g)+\epsilon(g))$ to be the conjugate of $S+\epsilon$ by $\hat \Phi$.
After having defined the notion of formal normal form of $S+\epsilon$ relative to $S$, we can state the following
\begin{theos}\cite{Stolo-intg-cras,stolo-ihes}\label{theo-scintg}
Let $S$ be an injective diagonal morphism such that the condition $(\omega(S))$ holds. Let $S+\epsilon$ be a nonlinear holomorphic 
deformation of $S$. Let us assume it admits an element of $\hom C\left(\lie g, \widehat{\cal O}_n^S\otimes_{{\bf C}}S(\lie g)\right)$ as 
a formal normal form. 
Then there is a formal normalizing diffeomorphism $\hat\Phi$ which is holomorphic in a neighborhood of $0$ in ${\bf C}^n$.
\end{theos}

Such a nonlinear deformation is called a {\bf holomorphic singular completely
  integrable system}. Let us make a few remarks about this. Assume that the
  ring $(\widehat{\cal O}_n)^S$ of formal first integrals of $S$ doesn't reduce to
  the constants. Then is is generated, as an algebra, by some monomials 
$x^{R_1},\ldots, x^{R_p}$ of ${\bf C}^n$ (\cite{stolo-ihes}[proposition 5.3.2,
  p.163]), $R_i\in {\bf N}^n$. These are the {\bf resonant monomials}. We define the {\bf resonant map} $\pi$ to be the map which associates to a point $x$ of ${\bf C}^n$, the values of the monomials at this point; that is 
$$
\pi : x\in {\bf C}^n  \mapsto (x^{R_1},\ldots, x^{R_p})\in{\cal C}_S\subset {\bf C}^p,
$$
where ${\cal C}_S$ is the algebraic subvariety of ${\bf C}^p$ defined by the algebraic relations among the $x^{R_i}$'s. The fibers of this mapping will be called the {\bf resonant varieties} (they may have singularities).

\begin{figure}[hbtp]\label{figure1}
  \begin{center}
    \leavevmode
    \input{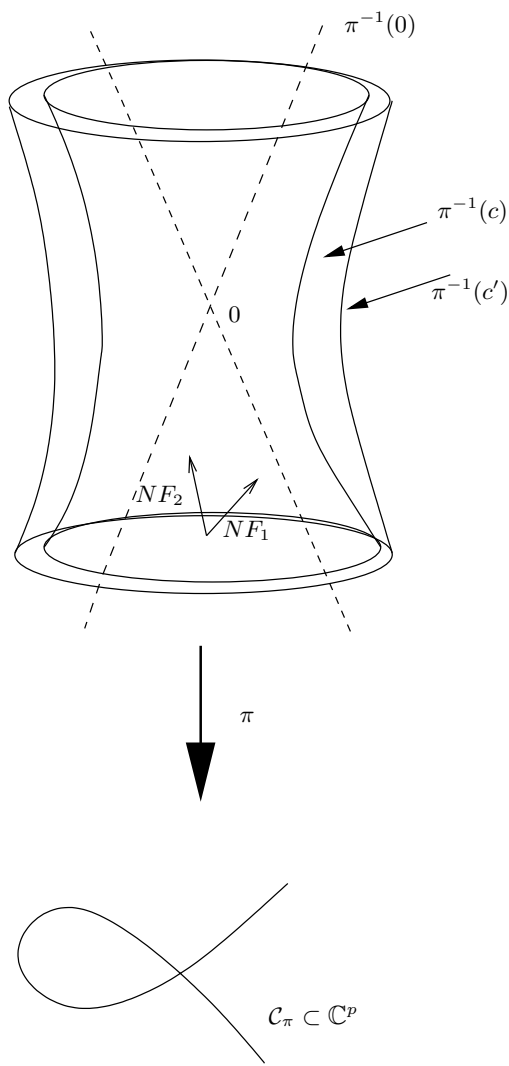} 
    \caption{Singular complete integrability: in the new holomorphic coordinate system, all the fibers (intersected with a fixed polydisc) are left invariant by the vector fields and their motion on it is a linear one}
    \label{fig:compl}
  \end{center}
\end{figure}

The conclusion of the previous theorem has the following geometric interpretation: let $D$ be a polydisc, centered at the origin and included in
the range of the holomorphic normalizing diffeomorphism. In the sequel, when
we say fiber of $\pi$, we mean its intersection with $D$. Our previous result
implies that, in the new holomorphic coordinates, the holomorphic vector
fields are tangent to each fiber over $\pi(D)$, they are pairwise commuting
and, when restricted to a fiber, they are just the restriction to the fiber of
a linear diagonal vector field whose eigenvalues depend only on the fiber (see
figure \ref{figure1}).

This reminds us of the classical complete integrability theorem of hamiltonian systems. The fibers, which can be regarded as the toric varieties, play the r\^ole of the {\it classical tori}. The flows associated to the restrictions to the fibers of the linear vector fields to which the original vector fields are conjugate to, play the r\^ole of the {\it quasi-periodic motions} on the tori.

\subsection{A KAM phenomenon for singular holomorphic vector fields}

With respect to what has already been said, the natural question one may ask is the following: starting from a holomorphic singular completely integrable system in a neighborhood of the origin of ${\bf C}^n$ (a common fixed point), we consider a holomorphic perturbation (in some sense) of one its vector fields. Does this perturbation still have invariant varieties in some neighborhood of the origin~?
Are these varieties biholomorphic to resonant varieties~? To which vector field on a resonant variety does the biholomorphism conjugate the restriction of the perturbation to an invariant variety~? Is there a ``big set" of surviving invariant varieties~?

The aim of this article is to answer these questions. Before fixing notation and giving precise statements, let us give a taste of what it is all about.

Let $S:\lie g \rightarrow \hvf n 1$ be as above. This defines a collection of
 linear diagonal vector fields on ${\bf C}^n$ we shall work with. Let $X$ be a
 holomorphic vector field in a neighborhood of the origin in ${\bf C}^n$. Let
 $X_0$ be a {\bf nondegenerate singular integrable} vector field (in the
 sense of R\"ussmann). 
We mean that $X_0$ is of the form
$$
X_0=\sum_{j=1}^l a_j(x^{R_1},\ldots, x^{R_p})S_j,\quad\quad a_j\in {\cal O}_n^S
$$
where the range of the map $(a_1,\ldots, a_l)$ from $({\bf C}^n,0)$ to 
$({\bf C}^l,0)$ is not included in any complex hyperplane.

Then, we consider a small holomorphic perturbation $X$ of $X_0$. Let us set
$X=X_0+R_{m_0}$ where $R_{m_0}$ is a germ of holomorphic vector field at the origin and of order greater than or equal to $m_0$ at that point. One of the difficulties is that
there are no natural actions-angles coordinates to play with. Nevertheless, we shall construct
something similar: we add new
variables $u_1,\ldots, u_p$ which correspond to the resonant monomials (which are assumed to be algebraically independent). These are the ``slow variables". 
To the holomorphic vector field $X$ in $({\bf C}^n,0)$, we shall associate a
holomorphic vector field $\tilde X$ in $({\bf C}^{n+p},0)$ where the
coordinate along $\frac{\partial }{\partial u_j}$ is the Lie derivative of
the resonant monomial $x^{R_j}$ along $X$. This vector field is tangent to the
variety 
$$
\Sigma=\{(x,u)\in {\bf C}^n\times{\bf C}^p,|\;\;u_j=x^{R_j},\quad j=1,\ldots,p\}
$$
and its restriction to it is nothing but $X$. We shall say that $\tilde X$ is
{\bf fibered over $X$}. We shall conjugate $\tilde X$ by germs of diffeomorphisms
which preserve the variety $\Sigma$. Such a germ will be built in the
following way: let $\Phi(x,u):=y=x+U(x,u)$ be a family of germs of biholomorphisms of
$({\bf C}^n,0)$, tangent to the identity at the origin and  parametrized over an open set ${\cal U}$. Let us set 
$v:= u + \pi(y)-\pi(x)$ and $\tilde \Phi(x,u):=(y,v)$. The latter is a germ of
diffeomorphism at $(0,b)$ and tangent to the identity at this point, for any
$b$. It leaves $\Sigma$ invariant. We shall say the $\tilde \Phi$ is {\bf fibered
over $\Phi$}. We shall define the notion of {\bf Lindstedt-Poincar\'e normal
form} of $\tilde X$ of order $k$ as follows: there exists a fibered
diffeomophism $\tilde \Phi_k$ such that 
\begin{displaymath}
(\tilde \Phi_k)_*\tilde X = \begin{cases} \dot y = NF^k(y,v)+R_{k+1}(y,v)+r_{\Sigma,k}(y,v)\\ 
\dot v = \pi_*(NF^k(y,v)+R_{k+1}(y,v)+r_{\Sigma,k}(y,v))\end{cases}.
\end{displaymath}
Here, $[X_0, \widetilde{ NF^k(y,v)}]$ vanishes on $\Sigma$, $R_{k+1}$ is of order greater than or equal to $k+1$ and $r_{\Sigma}$ vanishes on $\Sigma$.

If we were dealing with an integrable symplectic vector field $X_0$, we would require the
perturbation to be also symplectic. The analogue, in our general setting, is
an assumption on the {\bf Lindstedt-Poincar\'e normal form} of $X$. Namely, we
require that 
$$
NF^k(y,v)=\sum_{j=1}^l a_j^k(v)S_j(y).
$$
Hence, the Lindstedt-Poincar\'e normal form reads
\begin{displaymath}
(\tilde \Phi_k)_*\tilde X = \begin{cases} \dot y = \sum_{j=1}^l a_j^k(v)S_j(y)+R_{k+1}(y,v)+r_{\Sigma,k}(y,v)\\ 
\dot v = \pi_*(R_{k+1}(y,v)+r_{\Sigma,k}(y,v))\end{cases}.
\end{displaymath}
A perturbation $X$ of $X_0$ which has a Lindstedt-Poincar\'e normal form of this
type for any $k$ will be called {\bf good perturbation} of $X_0$.

For such a small perturbation, we shall prove that there is a neighborhood $U$ of the origin of ${\bf C}^n$ and there are compact sets of positive $2p$-measure belonging to the range of the resonant map $\pi$ with the following properties: let $K$ be such a compact set,
\begin{itemize}
\item for each $b\in K$, for each connected component of $\pi^{-1}(b)\cap U$, 
$X$ has an invariant holomorphic subset of some open set of ${\bf C}^n$ biholomorphic to the connected component of the fiber;
\item when $X$ is restricted to this invariant subset, the biholomorphism conjugates $X$ to the restriction of a linear vector field to the connected component of the fiber.
\end{itemize}

\begin{figure}[hbtp]
  \begin{center}
    \leavevmode
    \input{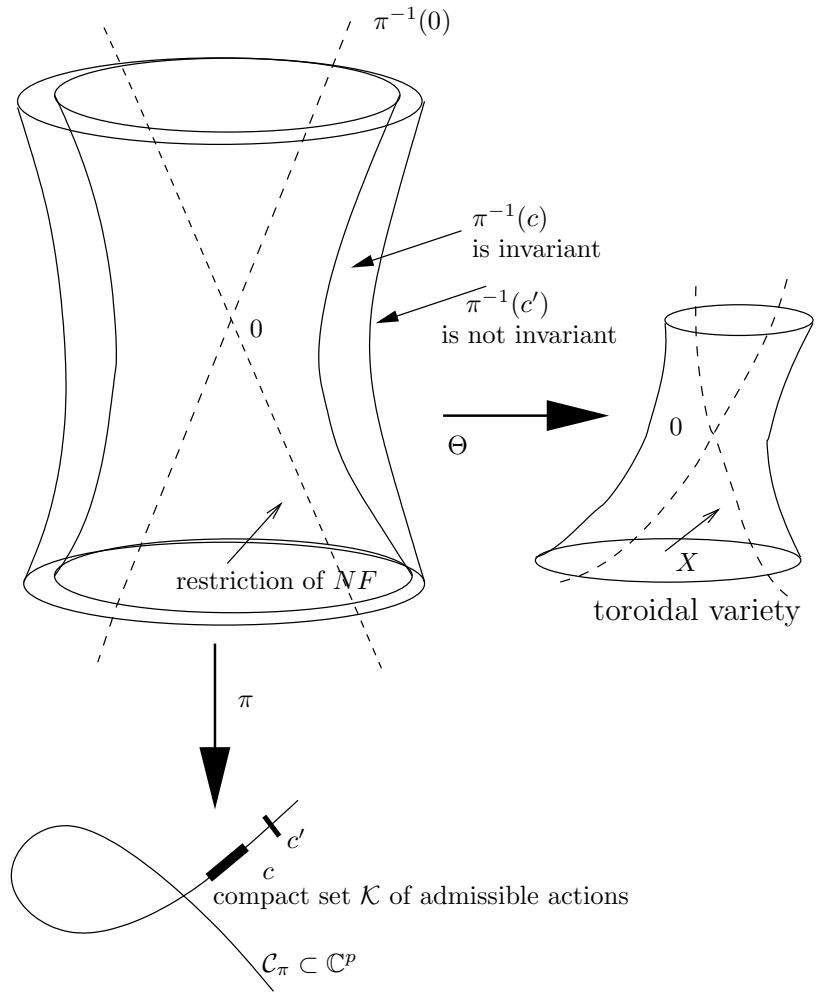} 
    \caption{KAM phenomenon: any fiber (when restricted to a fixed polydisc) over a compact set of positive measure is biholomorphic to an analytic subset left invariant by $X$ and the biholomorphism conjugates a linear motion on the fiber to the motion of $X$ on the corresponding invariant set }
    \label{fig:kam}
  \end{center}
\end{figure}
The germs of the ideas of the proof can be found in the work of Bibikov and Pliss
\cite{bibikov-pliss} (see also \cite{bibikov-book}[chap. 3]) although the
authors work with the Poincar\'e-Dulac normal form and this leads to an incomplete argument. The authors are concerned with peculiar systems of differential equations of the form
\begin{eqnarray*}
\frac{ dx_j}{dt}& = &i\lambda_jx_j+if_j(x_1,\ldots,x_n,y_1,\ldots,y_n)\\ 
\frac{ dy_j}{dt}& = &-i\lambda_jy_j-if_j(y_1,\ldots,y_n,x_1,\ldots,x_n) \quad(j=1,\ldots n),
\end{eqnarray*}
where $i^2=-1$ and the $\lambda_j$'s are real and uncommensurable numbers. The
authors are looking for invariant manifolds of the form $x_jy_j=constant$,
with $j=1,\ldots, n$. 

I would like to thank  M. Chaperon, Y. Colin de Verdi\`ere, B. Malgrange, J.-P. Ramis and J.-C. Yoccoz for their encouragements. I also thank F. Fauvet and P. Thomas for correcting some of my English language mistakes.


\section{Notation}

Let $R=(r_1,\ldots,r_n)\in \left({\bf R}^*_+\right)^n$; the open polydisc
centered at $0\in {\bf C}^n$ with polyradius $R$ will be denoted by
$D_n(0,R)=\{z\in {\bf C}^n\;|\;|z_i|<r_i\}$. If $r>0$ then $D_n(0,r)$ will
denote the polydisc $D_n(0,{(r,\ldots,r)})$. If $Q=(q_1,\ldots,q_n)\in \Bbb
N^n$, $|Q|=q_1+\cdots+q_n$ will denote the norm of $Q$.

Let ${\cal U}$ be an nonvoid connected open set of ${\bf C}^p$; then ${\cal
  O}_p({\cal U})$ (resp. ${\cal O}_p(\overline{\cal U})$) will denote the ring of holomorphic functions on ${\cal U}$ (resp. in a neighborhood of the closed set $\overline{\cal U}$). Let $f\in {\cal O}_p({\cal U})$ be such a function, we shall set $\|f\|_{\cal U}:=\sup_{x\in {\cal U}}|f(x)|$.

\subsection{Norms}

Let ${\cal U}$ be a nonvoid connected open set of ${\bf C}^p$. 
Let us set
$$
B_{\cal U}:={\cal O}_p({\cal U})\otimes_{{\bf C}}{\bf C} [[x_1,\ldots,x_n]].
$$
It is an algebra over ${\bf C}$. Let $f=\sum_{Q\in {\bf N}^n}{f_Q(u)x^Q}$ be an element of $B_{\cal U}$ where the $f_Q$'s belong to ${\cal O}_p({\cal U})$.  
We shall set 
$$
\bar f:= \sum_{Q\in {\bf N}^n}{\|f_Q\|_{\cal U}x^Q}\in {\bf C}[[x_1,\ldots,x_n]].
$$
The {\bf order} of such an element is the smallest integer $k\in {\bf N}$ such that there exists $Q\in {\bf N}^n$ of norm equal to $k$ and 
$f_Q\not\equiv 0$; this integer will be denoted by $ord_{\cal U}(f)$. Let $k$ be an integer; the {\bf $k$-jet} of $f$ is 
$$
J^{k}_{{\cal U}}(f):=\sum_{Q\in {\bf N}^n, |Q|\leq k}{f_Q(u)x^Q}.
$$

We shall say that an element $g$ dominates $f$ if for any multiindex $Q\in {\bf N}^n$, $\|f_Q\|_{\cal U}$ is less than or equal to $\|g_Q\|_{\cal U}$;
in this case, we shall write $f\prec_{\cal U} g$. Let $R=(r_1,\ldots,r_n)\in \left({\bf R}^*_+\right)^n$; let us set 
$$
|f|_{{\cal U},R}:=\sum_{Q\in {\bf N}^n}{\|f_Q\|_{\cal U}R^{Q}}=\bar f(r_1,\ldots,r_n).
$$
We have the following properties:
\begin{eqnarray*}
fg & \prec_{\cal U} & \bar f\bar g,\\
\text{if}\;f\prec_{\cal U} g & \text{then} & |f|_{{\cal U},R}\leq |g|_{{\cal U},R}.\\
\end{eqnarray*}
Let $f\in B_{{\cal U}}$ such that $ord_{{\cal U}}(f)\geq m$. We assume that 
$|f|_{{\cal U},r}$ is finite for some positive $r$. Then, for any positive $R<r$, we have
\begin{eqnarray}
|f|_{{\cal U},R} & \leq & \left(\frac{R}{r}\right)^m |f|_{{\cal U},r}\\
sup_{z\in {\cal U}\times D_{r}}|f(z)| & \leq & |f|_{{\cal U},r}.
\end{eqnarray}

We shall define the set 
$$
{\cal H}_n({\cal U},R)=\left\{f\in {\cal O}_n(\overline{\cal U})[[x_1,\ldots, x_n]]\;|\;|f|_{{\cal U},R}<+\infty\right\}.
$$
\begin{rems}[Important remark]\label{importante_remarque}
Let $x^{R_1},\ldots, x^{R_p}$ be a finite number of (nontrivial) monomials of ${\bf C}^n$. Let $\Sigma$ be the submanifold of ${\bf C}^n\times {\bf C}^p$ defined to be
$$
\Sigma = \left\{ (x,u)\in {\bf C}^n\times {\bf C}^p\;|\; u_i=x^{R_i}\;\; i=1,\ldots, p\right\}.
$$
Let $f(x,u)=\sum_{Q\in {\bf N}^n}f_Q(u)x^Q$ be a holomorphic function on
$D_n(0,r)\times D_p(b,r')$ which is assumed to intersect $\Sigma$. We shall
denote by {\bf $f_{|\Sigma}$ the function obtained by replacing each monomial $x^{R_i}$ by $u_i$}. If $|b|+r'\leq r^{|R_i|}$, we have
$$
|f_{|\Sigma}|_{D_p(b,r'),r}\leq |f|_{D_p(b,r'),r}.
$$
This is due to the fact that 
$$
\|u_i\|_{D_p(b,r')}\leq |b|+r'\leq r^{|R_i|}=|x^{R_i}|_{D_p(b,r'),r}.
$$
\end{rems}
\subsection{Spaces of vector fields and spaces of functions}

Let us set some notation which will be used throughout this article. 
Let $k$ be a positive integer:
\begin{itemize}
\item $\hvf n k$ denotes the ${\bf C}$-space of homogeneous polynomial 
vector fields on ${\bf C}^n$ and of degree $k$,
\item $\pvf n m k$ denotes the ${\bf C}$-space of polynomial 
vector fields on ${\bf C}^n$, of order $\geq m$ and of degree $\leq k$ ($m\leq k$),
\index{$\pvf n m k$ : ${\bf C}$-space of polynomial 
vector fields on ${\bf C}^n$, of order $\geq m$ and of degree $\leq k$}
\item $\fvfo n k$ denotes the ${\bf C}$-space of formal 
vector fields on ${\bf C}^n$ and of order $\geq k$ at $0$,
\item $\vfo n k$ denotes the ${\bf C}$-space of germs of holomorphic 
vector fields on $({\bf C}^n,0)$ and of order $\geq k$ at $0$,
\item $p_n^k$ denotes the ${\bf C}$-space of homogeneous polynomials 
on ${\bf C}^n$ and of degree $k$,
\item $\widehat{\cal O}_n$ denotes the ring of formal power series in ${\bf C}^n$,
\item ${\cal O}_n$ denotes the ring of germs at $0$ of holomorphic functions in ${\bf C}^n$.
\end{itemize}


\section{Normal forms of vector fields, invariants and nondegeneracy}
\subsection{Normal forms}

Let $X\in \fvfo n 1$ be a formal vector field of ${\bf C}^n$ and $s$ its linear
part (assumed not to be zero) at the origin.
\begin{defis}
$X$ is said to be a normal form if  $[s,X]=0$.
\end{defis}
\begin{props}[Poincar\'e-Dulac formal normal form]\cite{Arn2}
Let $X$ belong $ \fvfo n 1$ as above. Then there exists a formal diffeomorphism $\hat \Phi$, vanishing at $0$ and tangent to identity at this point such that 
$$
\hat\Phi^*X=s+\hat N,\quad [s,\hat N]=0;
$$
where $\hat N$ is a nonlinear formal vector field.
\end{props}
This means that we can find a formal change of coordinates in which $X$ is
transformed into a normal form. In general, the normalizing diffeomorphism is
not unique and diverges.

\subsection{Invariants and toric varieties}
Let $\lie g$ be a complex $l$-dimensional commutative Lie algebra. Let $S:\lie g\rightarrow \hvf n 1$ be a Lie morphism from $\lie g$ to the Lie algebra of linear vector fields of ${\bf C}^n$. It is assumed to be injective and  semi-simple. This means that we can find a basis $\{g_1,\ldots, g_l\}$ of $\lie g$ such that, for each index $i$, $S(g_i)$ is a linear diagonal vector field of ${\bf C}^n$ and the family $\{S(g_i)\}_{i=1,\ldots,l}$ is linearly independent over ${\bf C}$. 
Let us define 
\begin{eqnarray*}
\widehat{\cal O}_n^S & = & \left\{f\in \widehat{\cal O}_n\,|\, \forall g\in\lie g,\;{\cal L}_{S(g)}(f)=0\right\},\\
\left(\fvfo n 1\right)^S & = & \left\{Y\in \fvfo n 1\,|\, \forall g\in\lie g,\;[S(g),Y]=0\right\},
\end{eqnarray*}
to be the ring of formal first integrals of $S$ and the formal centralizer of
$S$ respectively. Here, ${\cal L}_{S(g)}(f)$ denotes the Lie derivative of $f$
along $S(g)$.

\subsubsection{Weights, weight spaces of the morphism $S$}

The morphism $S$ induces a representation $\rho_k$ of $\lie g$ in $\hvf n k$ defined to be $\rho_k(g)p=[S(g),p]$.
The linear forms $\alpha_{Q,i}:= (Q,\lambda)-\lambda_i$, $|Q|=k$ are the {\bf weights} of order $k$ of this representation. There is a decomposition of $\hvf n k$ into direct sums of {\bf weight spaces} of this representation:
$$
\hvf n k=\bigoplus_{\alpha}\left(\hvf n k\right)_{\alpha}
$$
where 
$$
\left(\hvf n k\right)_{\alpha}=\{p\in \hvf n k\,|\,\forall g\in \lie g,\,[S(g),p]=\alpha(g)p\}\neq \{0\}
$$ 
denotes the $\alpha$-weight space associated to the weight $\alpha$. The set
of nonzero weights of $S$ into $\pvf n {k} {m}$ will be denoted by ${\cal W}_{n,*}^{k,m}$. We refer to our previous article \cite{stolo-ihes}[Chapter 5-The fundamental structures- p.158-168]
for more details about this topic.

Let $\|.\|$ be a norm in ${\lie g}^*$. We define the sequence of positive numbers
$$
\omega_k(S)=\inf\left\{\|\alpha\|, \;\text{for all {\bf nonzero} weights of } S \text{ into } \pvf n {2^k+1} {2^{k+1}}\right\}.
$$

\subsubsection{Toric varieties}

If $\widehat{\cal O}_n^S\neq {\bf C}$, we know from \cite{stolo-ihes}[prop. 5.3.2] that there exists a finite number of monomials $x^{R_1},\ldots, x^{R_p}$ such that $$\widehat{\cal O}_n^S={\bf C}[[x^{R_1},\ldots, x^{R_p}]].$$ Moreover, $\left(\fvfo n 1\right)^S$ is a $\widehat{\cal O}_n^S$-module of finite type. Let us define the map 
$$
\pi :x\in {\bf C}^n \mapsto  (x^{R_1},\ldots, x^{R_p})\in  {\bf C}^p.
$$
Let ${\cal C}_{\pi}$ denote the algebraic subvariety of ${\bf C}^p$ defined by the algebraic relations among the monomials $x^{R_1},\ldots, x^{R_p}$: 
$$
{\cal C}_{\pi}=\{u\in {\bf C}^p\,|\, P(u)=0,\; \forall P\in{\bf C}[u_1,\ldots,u_p]\text{ such that } P(x^{R_1},\ldots, x^{R_p})=0\}.
$$
It is well known that such a variety is a {\bf toric variety}. In fact, there is an action of the algebraic torus 
$\left({\bf C}^*\right)^n$ on ${\bf C}^p$ defined by $x.u=(u_1x^{R_1},\ldots, u_px^{R_p})$ whenever $u=(u_1,\ldots, u_p)\in {{\bf C}}^p$ and $x=(x_1,\ldots x_n)\in \left({\bf C}^* \right)^n$.
The affine variety ${\cal C}_{\pi}$ is the Zariski closure  of the orbit of 
the point ${\bf 1}=(1,\ldots, 1)$. If $b\in {\bf C}^p$ belongs to this orbit, then the set 
$$\pi^{-1}(b)=\{x\in \left({\bf C}^*\right)^n\,|\,x^{R_1}=b_1,\ldots, x^{R_p}=b_p\}$$ is isomorphic the isotropy subgroup of ${\bf 1}$; that is 
the subgroup $G$ of $\left({\bf C}^*\right)^n$ defined by $G=\{x\in \left({\bf C}^*\right)^n\,|\,x.{\bf 1}={\bf 1}\}$.

Let ${\cal I}_{\pi}$ (resp. $\widehat{\cal I}_{\pi}$) denote the ideal of definition of ${\cal C}_{\pi}$ (resp. ${\cal I}_{\pi}\otimes\widehat{\cal O}_n$).
Therefore, we have an isomorphism $\phi_{\pi}:{\bf C}[[u_1,\ldots,
u_p]]/\widehat{\cal I}_{\pi}\stackrel{\sim}{\rightarrow}\widehat{\cal O}_n^S$
defined to be $\phi_{\pi}([F])=F\circ \pi$ for any representative $F\in{\bf C}[[u_1,\ldots, u_p]]$ of $[F]$. We shall also denote by $\phi_{\pi}$ the induced map on $\left({\bf C}[[u_1,\ldots, u_p]]/\widehat{\cal I}_{\pi}\right)^k$, for any positive integer $k$.

Let $R$ be the $p\times n$-matrix which rows are the exponents of the monomials defining $\pi$:
$$
R=\begin{pmatrix}R_1\\ \vdots\\ R_p\end{pmatrix}=\begin{pmatrix}R_{1,1}&\cdots & R_{1,n}\\ \vdots& & \vdots\\ R_{p,1}& \cdots &R_{p,n}\end{pmatrix}.
$$ 
Let $r$ be the rank of $R$. Let $M=(m_{i,j})_{\substack{1\leq i\leq p\\ 1\leq j\leq n}}$ be a $p\times n$-matrix and let $k$ be smaller than or equal to $\min (n,p)$. Let $I=\{i_1<\cdots<i_k\}$ (resp. $J=\{j_1<\cdots<j_k\}$) be a set of increasing integers of $[1,p]$ (resp. $[1,n]$). We shall say that both $I$ and $J$ are $k$-sets.
We shall write $M_{I,J}$ for the matrix $M_{I,J}=(m_{i_s,j_t})_{1\leq s,t\leq k}$.
\begin{props}
Let $r$ be a positive number and let $D_r$ be the polydisc centered at the origin of ${\bf C}^n$ with radius $r$.
Let $H=\{x\in {\bf C}^n\,|\,x_1\cdots x_n=0\}$ be the union of the coordinate hyperplanes. Let $D_r\setminus H$ be the set of points of $D_r$ not belonging to $H$. Then, 
\begin{enumerate}
\item $D_r\setminus H$ is an open connected set of ${\bf C}^n$,
\item $\pi_{|D_r\setminus H}$ is of constant rank equal to the rank of $R$.
\end{enumerate}
\end{props}
\begin{proof}
The first point is classical since $H$ is of complex codimension $1$.
Let $k$ be a positive integer smaller than or equal to $\min (p,n)$ and let $I=\{i_1<\cdots<i_k\}$ (resp. $J=\{j_1<\cdots<j_k\}$) be a set of increasing integers of $[1,p]$ (resp. $[1,n]$). Let $x\in D_r\setminus H$. We have 
$$
(D\pi(x))_{I,J} = \left( r_{j,i}\frac{x^{R_i}}{x_j}\right)_{\substack{i\in I\\j\in J}}.
$$
According to the multilinearity property of the determinant, we obtain
$$
x_{j_1}\cdots x_{j_k}\det (D\pi(x))_{I,J} = \det\left( r_{j,i}x^{R_i}\right)_{\substack{i\in I\\j\in J}} = x^{R_{i_1}}\cdots x^{R_{i_k}}\det\left( r_{i,j}\right)_{\substack{i\in I\\j\in J}}.
$$
Since the zero set of $x_{j_1}\cdots x_{j_k}$ (resp. $x^{R_{i_1}}\cdots x^{R_{i_k}}$) is included in $H$ then, for all $x\in D_r\setminus H$, the determinant of $(D\pi(x))_{I,J}$ vanishes if and only if the determinant of $R_{I,J}$ does. Let $s$ be the rank of $R$. Then, for all positive integers $s<t\leq \min (n,p)$, and all $t$-sets $I,J$, the determinant of $R_{I,J}$ vanishes; thus, so does $\det (D\pi(x))_{I,J}$. Moreover, there are $s$-sets $I,J$ such that the determinant of $R_{I,J}$ is not zero; thus, $\det (D\pi(x))_{I,J}\neq 0$. As a consequence, 
the rank of $\pi_{D_r\setminus H}$ is equal to the rank of $R$.
\qed\end{proof}
Let us recall the rank theorem:
\begin{theos}\cite{chirka}[p.307]
Let $D$ be an open connected set in ${\bf C}^n$ and let $f:D\rightarrow {\bf C}^p$ be a holomorphic map such that $\text{rank}_z f=r$ for all $z\in D$. Then, for any point $a\in D$, there exists a neighborhood $U$ of $a$ in $D$ and a neighborhood $V$ of $f(a)$ in ${\bf C}^p$ such that:
\begin{enumerate}
\item $f(U)$ is an $r$-dimensional complex submanifold in $V$,
\item there is an $r$-dimensional complex plane $L$ passing through $a$ such that 
$f$ maps biholomorphically $L\cap U$ onto $f(U)$.
\end{enumerate}
\end{theos}
We may apply this result to our situation: $\pi$ is of constant rank $s$ on $D_r\setminus H$; therefore for any point $a\in D_r\setminus H$, there is a neighborhood $U_a$ of $a$ in $D_r\setminus H$ such that $\pi(U_a)$ is an $r$-dimensional submanifold of ${\bf C}^p$.

\subsection{Nondegeneracy of vector fields}
{\bf From now on, we shall assume that the resonant monomials $x^{R_1},\ldots, x^{R_p}$ are algebraically independent.}
\begin{defis}
\begin{itemize}
\item A formal map $\hat f : {\bf C}^p\rightarrow {\bf C}^l$ is a $l$-tuple of formal power series of ${\bf C}^p$.
\item A formal map $\hat f : {\bf C}^p\rightarrow {\bf C}^l$ is said to be {\bf nondegenerate} (in the sense of R\"ussmann) if for any nonzero $c=(c_1,\ldots, c_l)\in {\bf C}^l$, $<c,\hat f>:= \sum_{j=1}^l c_j\hat f_j$ doesn't vanish identically.
\end{itemize}
\end{defis}
\begin{rems}
If the $k$-jet, $J^k(\hat f)$, of a formal map $\hat f$ is nondegenerate, so is $J^{k+m}(\hat f)$.
\end{rems}
\begin{lemms}
Let $\hat F=(\hat F_1,\ldots, \hat F_l): {\bf C}^n\hat\rightarrow {\bf C}^l$ be a formal map such that for all $1\leq i\leq l$, $\hat F_i\in \widehat{\cal O}_n^S$.
The following statements are equivalent:
\begin{enumerate}
\item $\hat F: {\bf C}^n\rightarrow {\bf C}^l$ is a nondegenerate formal map,
\item $\phi_{\pi}^{-1}(\hat F): {\bf C}^p\rightarrow {\bf C}^l$ is a nondegenerate formal map.
\end{enumerate}
\end{lemms}
\begin{proof}
There is a formal power series $F$ of ${\bf C}^p$ such that $F\circ \pi=\hat F$. Since the resonant monomials are algebraically independent, $<c,\hat F>=0$ if and only if $<c,F>=0$.
\qed\end{proof}

\begin{defis}
A holomorphic vector field $\sum_{j=1}^la_jS_j$ will be said to be nondegenerate if the map 
$(a_1,\ldots, a_l)$ is nondegenerate.
\end{defis}


\subsection{Fibered vector fields and diffeomorphisms along $\pi$}\label{fibrered}

Let $\Sigma$ be the graph of $\pi$; that is the algebraic subvariety of ${\bf C}^n\times {\bf C}^p$ defined by
$$
\Sigma:=\left\{(x,u)\in {\bf C}^n\times {\bf C}^p\;|\;u_i=x^{R_i},\; i=1,\ldots,p\right\}.
$$

Let ${\cal U}$ be an open set in ${\bf C}^p$ and let $f$ be a holomorphic function over $D_n(0,r)\times {\cal U}$ which vanishes as well as its derivative at $x=0$. It can be written as 
$$
f(x,u)=\sum_{Q\in {\bf N}^n_2}f_Q(u)x^Q.
$$
Moreover, it can be written also as
$$
f=\sideset{}{'}\sum_{Q\in {\bf N}^n}f'_Q(x^{R},u)x^Q
$$
where the sum ranges over $\{0\}$ together with the set of $Q\in {\bf N}^n$ for which $x^Q$ is not divisible by any $x^{R_i}$.
Two such series $f$ and $g$ will be said {\bf equivalent modulo $\Sigma$} if 
$$
f'_Q(x^{R},u)=g'_Q(x^{R},u)\;\; \text{on } \Sigma.
$$

We shall denote $[f]$ (or ${\bf f}$) the equivalence class and write it as the series 
$$
{\bf f}:=[f]:=\sideset{}{'}\sum_{Q\in {\bf N}^n}[f'_Q(x^{R},u)]x^Q
$$
where $[f'_Q(x^{R},u)]$ denotes the equivalence class modulo $\Sigma$ of $f'_Q(x^{R},u)$.
Let $$X=\sum_{i=1}^{n}X_i\partial /\partial x_i+\sum_{k=1}^{p}X_k\partial /\partial u_k$$ be a holomorphic vector field on $D_n(0,r)\times {\cal U}$. We shall define $[X]$ (or ${\bf X}$) to be
$$
{\bf X}:=[X]:=\sum_{i=1}^{n}[X_i]\partial /\partial x_i+\sum_{k=1}^{p}[X_k]\partial /\partial u_k.
$$
\begin{defis}
The restriction $f_{|\Sigma}$ of $f$ to $\Sigma$ is the formal power series
defined to be 
$$
f_{|\Sigma}:=\sideset{}{'}\sum_{Q\in {\bf N}^n}[f'_Q(u,u)]x^Q,
$$
where the sum ranges over $\{0\}$ together with the set of $Q\in {\bf N}^n$ for which $x^Q$ is not divisible by any $x^{R_i}$.
\end{defis}
Let $X$ be a holomorphic family over ${\cal U}$ of holomorphic of vector fields in a neighborhood $V$ of the origin in ${\bf C}^n$.
We mean that 
$$
X(x,u)=\sum_{i=1}^nX_i(x,u)\frac{\partial }{\partial x_i}
$$
where the $X_i(x,u)$'s belong to ${\cal O}_n(V)\otimes_{{\bf C}}{\cal O}_p({\cal U})$.

Let us set 
$$
\tilde X(x,u):= (X,\pi_*X)(x,u)=\sum_{i=1}^nX_i(x,u)\frac{\partial }{\partial x_i}+\sum_{k=1}^pX(x^{R_k})(x,u)\frac{\partial}{\partial u_k}.
$$
As usual, $X(x^{R_k})$ denotes the Lie derivative of $x^{R_k}$ along $X$, namely:
$$
X(x^{R_k})(x,u)= \sum_{i=1}^nX_i(x,u)\frac{\partial x^{R_k}}{\partial x_i}.
$$
The vector field $\tilde X$ will be said to be {\bf fibered over $\pi$ (with respect to $X$)}. One of the features of $\tilde X$ is that it is tangent to $\Sigma$ and its restriction to it is precisely the vector field $X(\pi(x),x)$.

Let $\Phi$ be a holomorphic family over ${\cal U}$ of germs of holomorphic diffeomorphisms of $({\bf C}^n,0)$ tangent to identity at $0$ in ${\bf C}^n$. 
We shall write 
$$
\Phi(x,u)= x+\phi(x,u),
$$
where $\phi$ belongs to ${\cal O}_p({\cal U})\otimes{\cal M}_n^2$.
Let us set
$$
\tilde \Phi(x,u):=\left(\Phi(x,u),u+(\pi(\Phi(x,u))-\pi(x))\right).
$$
If $b$ belongs to ${\cal U}$, then $\tilde \Phi(0,b)=(0,b)$. Moreover, $\tilde
\Phi$ is tangent to identity at that point. Thus, it is a local diffeomorphism
at $(0,b)$. The diffeomorphism $\tilde X$ will be said to be {\bf fibered over $\pi$ (with respect to $\Phi$)}. It leaves the variety $\Sigma$ (globally) invariant. In fact, if $u=\pi(x)$ then 
$$
\tilde \Phi(\pi(x),x)=(\Phi(\pi(x),x),\pi(\Phi(\pi(x),x)).
$$
In the sequel, {\bf $pr_1$ (resp. $pr_2$) will denote the projection onto ${\bf C}^n$ (resp. ${\bf C}^p$)}.
\begin{lemms}
If both $\tilde X$ and $\tilde \Phi$ are fibered over $\pi$ (with respect to $X$ and $\Phi$ respectively) then so is $\tilde\Phi_*\tilde X$. Moreover, the restriction of $\tilde \Phi_*\tilde X$ on $\Sigma$ does depend only on the restrictions of $\tilde X$ and $\tilde \Phi$ to $\Sigma$.
\end{lemms}
\begin{proof}
In fact, we have 
\begin{eqnarray*}
\tilde\Phi_*\tilde X(y,v)& = & D\tilde\Phi(x,u)\tilde X(x,u)\\
&=& \left(D_x\Phi(x,u)X(x,u)+D_u\Phi(x,u)\pi_*X(x,u),\right.\\
& & \left. D_x(\pi(\Phi(x,u))-\pi(x))X(x,u)+D_u(u+\pi(\Phi(x,u)))\pi_*X(x,u)\right).
\end{eqnarray*}
But precisely, the quantity
$$
D_x(\pi(\Phi(x,u))-\pi(x))X(x,u)+D_u(u+\pi(\Phi(x,u)))\pi_*X(x,u) 
$$
is equal to 
$$
D\pi(\Phi(x,u))\left(D_x\Phi(x,u)X(x,u)+D_u\Phi(x,u)\pi_*X(x,u)\right).
$$

It is sufficient to prove this with
$X=(f(x^R,u)+f_{\Sigma}(x^R,u))x^Q\partial/\partial x_i$ and $(\tilde y,\tilde
v):=\tilde \Phi(x,u)$ with
\begin{eqnarray*}
\tilde y & = & x+x^P(g(x^R,u)+g_{\Sigma}(x^R,u))\partial/\partial x_i\\ 
\tilde v & = & u +\pi(x+x^P(g(x^R,u)+g_{\Sigma}(x^R,u))\partial/\partial x_i)-\pi(x)).
\end{eqnarray*}
Here, we assume that neither $x^Q$ nor $x^P$ are divisible by any $x^{R_k}$. Moreover, both
$f_{\Sigma}$ and $g_{\Sigma}$ belong to the ideal generated by the $x^{R_k}-u_k$'s. 
Let us set
\begin{eqnarray*}
(y,v) & := & \tilde \Psi(x,u)\\
y & := & x+x^Pg(x^R,u)\partial/\partial x_i\\
v & := & u +\pi(x+x^Pg(x^R,u)\partial/\partial x_i)-\pi(x)\\
Y(x,u) & := & f(x^R,u)x^Q\partial/\partial x_i.
\end{eqnarray*}
The previous computation shows that
\begin{eqnarray*}
pr_1(\tilde\Phi_*\tilde X(y,v)) & = &
D_x\Phi(x,u)(f(x^R,u)+f_{\Sigma}(x^R,u))x^Q(\frac{\partial}{\partial
  x_i})\\
& & +D_u\Phi(x,u)(\sum_{k=1}^p(f(x^R,u)+f_{\Sigma}(x^R,u))x^Q\frac{\partial x^{R_k}}{\partial
  x_i}\frac{\partial}{\partial u_k}).
\end{eqnarray*}
Let us write $g_{\Sigma(x^R,u)}=\sum_{k=1}^p(x^{R_k}-u_k)g_{\Sigma,k}(x^R,u)$. Let us define
the differential operator 
$$
\frac{\partial g(x^R,u)}{\partial v_k} := \frac{\partial g}{\partial t_{k}}(x^R,u)+\frac{\partial g}{\partial u_k}(x^R,u)
$$
where $g(x^R,u)=\pi_*(g(t,u))$.
Then, we obtain
\begin{eqnarray*}
pr_1(\tilde\Phi_*\tilde X(\tilde y,\tilde v)) & = & \left[x^P\sum_{k=1}^p\left(\frac{\partial
  g(x^R,u)}{\partial v_k}+\sum_{m=1}^p(x^{R_m}-u_m)\frac{\partial
  g_{\Sigma,m(x^R,u)}}{\partial v_k}\right)\frac{\partial x^{R_k}}{\partial
  x_i}\right.\\
& & \left. +(g(x^R,u)+g_{\Sigma}(x^R,u))\frac{\partial x^P}{\partial x_i}+1\right](f(x^R,u)+f_{\Sigma}(x^R,u))x^Q\frac{\partial}{\partial y_i}.
\end{eqnarray*}
This is due to the fact that $\frac{\partial (x^{R_m}-u_m)}{\partial v_k}=0$
for all integers $1\leq k,m\leq p$.
As a consequence, if $(x,u)$ belongs to $\Sigma$, then $(\tilde y,\tilde v)=(y,v)$. Moreover, we have
$$
pr_1(\tilde\Phi_*\tilde X(y,v)) = pr_1(\tilde\Psi_*\tilde Y(y,v)) .
$$

\qed\end{proof}

%

\subsection{Lindstedt-Poincar\'e normal forms of fibered vector fields along $\pi$}

Let $X=s_1(x)+R(x)$ be a germ of vector field in a neighborhood of the origin
in ${\bf C}^n$ which is a nonlinear perturbation of the polynomial vector
field $s_1=\sum_{j=1}^la_j(x^R)S_j$. Let us assume that the order $m_0$ of $R$ is
greater than the degree of $s_1(x)$.

Let ${\cal U}$ be an open set in ${\bf C}^p$. Let $Y(x,u)$ be a family, over
${\cal U}$, of germs of vector fields of ${\bf C}^n$, whose restriction to
$\Sigma$  is equal to $X$. It is required that $Y(x,u)=s_1(x,u)+ R(x,u)$ where
$s_1(x,u)$ is defined to be $\sum_{j=1}^la_j(u)S_j$ and where $R(x,u)$ is
nonlinear (in $x$). We define the holomorphic vector field $\tilde X$ on an
open set of ${\bf C}^n\times{\bf C}^p$ to be
$$
\tilde X (x,u)= (Y(x,u), \pi_*Y(x,u))=(s_1(x,u)+ R(x,u), \pi_*R(x,u)).
$$
\begin{lemms}
Let $U(x,u)$ (resp. $s(x,u)$) be a representative, over ${\cal U}$, of some
nonlinear vector field (resp. $s_1(u,x)$). Then, the class of the Lie
bracket\\
$[s(x,u),\tilde U(x,u)]$ depends only on the classes $[s]$ and $[\tilde
U]$. Moreover, if $U(x,u)$ belongs to $\alpha$-weight space of $S$ then so do
$\tilde U$ and $[s_1(x,u), \tilde U(x,u)]$.
\end{lemms}
\begin{proof}
It is sufficient to prove this for $U(x,u)=\tilde b(x,u)x^Q\frac{\partial }{\partial x_i}$
and $s(x,u)=\sum_{j=1}^l\tilde a_j(x,u)S_j$ where $\tilde a_j=a_j(x^R,u)+a_{j,\Sigma}(x^R,u)$
and $\tilde b=b(x^R,u)+b_{\Sigma}(x^R,u)$. We have
\begin{eqnarray*}
[s(x,u), \tilde U(x,u)] & = & \sum_{j=1}^l \tilde a_j\tilde b \left[S_j, x^Q\frac{\partial }{\partial x_i}\right]-\tilde bx^Q\frac{\partial \tilde a_j}{\partial x_i}S_j\\
& & +\sum_{j=1}^l\left[\tilde a_jS_j, \sum_{k=1}^p \tilde bx^Q\frac{\partial x^{R_k}}{\partial x_i}\frac{\partial}{\partial u_k}\right].
\end{eqnarray*}
Since the $x^{R_k}$'s are common first integrals of the $S_j$'s, we have $S_j(\tilde b)=0$. This leads to
\begin{eqnarray*}
\left[\tilde a_jS_j, \sum_{k=1}^p \tilde bx^Q\frac{\partial x^{R_k}}{\partial x_i}\frac{\partial}{\partial u_k}\right] & = & \sum_{k=1}^p\left[ \tilde a_j\tilde bS_j\left(x^Q\frac{\partial x^{R_k}}{\partial x_i}\right)\frac{\partial}{\partial u_k}\right.\\
& &\left. -\tilde bx^Q\frac{\partial x^{R_k}}{\partial x_i}\frac{\partial \tilde a_j}{\partial u_k}S_j\right].
\end{eqnarray*}
It is then sufficient to show that the class of
$$
B:=\frac{\partial \tilde a_j}{\partial x_i}+\sum_{k=1}^p\frac{\partial x^{R_k}}{\partial x_i}\frac{\partial \tilde a_j}{\partial u_k}=\sum_{k=1}^p\frac{\partial x^{R_k}}{\partial x_i}\left(\frac{\partial \tilde a_j}{\partial v_k}(x^{R},u)\right)
$$
depends only on the class of $\tilde a_j$. Let us write $a_{j,\Sigma}=\sum_{m=1}^p(x^{R_m}-u_m)a_{j,\Sigma,m}$. Since $\partial (x^{R_m}-u_m)/\partial v_k=0$, we have 
$$
\frac{\partial \tilde a_j(x^{R},u)}{\partial v_k}= \frac{\partial a_j(x^{R},u)}{\partial v_k}+\sum_{m=1}^p (x^{R_m}-u_m)\frac{\partial a_{j,\Sigma, m}(x^{R},u)}{\partial v_k},
$$
which ends the proof.

Let us assume that $U$ belongs to the $\alpha$-weight space of $S$. According to the previous computations, we have 
$$
[S_k,\tilde U]  = [S_k, U]+\sum_{m=1}^p S_k(U(x^{R_m}))\frac{\partial }{\partial u_m}.
$$
The result follows from the equality $S_k(U(x^{R_m}))=[S_k,U](x^{R_m})= \alpha(g_k)U(x^{R_m})$.

\qed\end{proof}
\begin{defis}[Lindstedt-Poincar\'e normal form]
The fibered vector field $\tilde X$ over $\pi$ with respect $X$ is said to be {\bf normalized up to order $m$ along $\Sigma$}, if there exists an open set ${\cal U}_m$ and a fibered diffeomorphism $\tilde \Phi$ over ${\cal U}_m$ such that 
$$
\tilde \Phi_*\tilde X = \widetilde{NF^m}+\widetilde{R_{m+1}}+\widetilde{r_{\Sigma}}
$$
where $\widetilde{R_{m+1}}$ is a fibered vector field with respect to a
vector field of order greater than or equal to $m+1$. The vector field $\widetilde{NF^m}$ commutes with $s_1(x,u)$ and does not depend on the choice of a representant of $X$.

The vector field $X$ is said to be {\bf Lindstedt-Poincar\'e normalized up to order $m$} if one of its fibered vector field $\tilde X$ is normalized up to order $m$.
We shall also say that the restriction $\widetilde{NF^m}_{|\Sigma}$ is a {\bf Lindstedt-Poincar\'e normal form} up order $m$ of $X$.
\end{defis}
\begin{lemms}
If $s_0$ is nondegenerate then $\tilde X$ admits a Lindstedt-Poincar\'e normal form up to any order.
\end{lemms}
\begin{proof}
Let us assume that $\tilde X$ is normalized up to order $m$ and let us write
$$
pr_1\tilde X = NF^m(x,u)+R_{m+1}(x,u)+r_{\Sigma}(x,u)
$$
where $u$ belongs to an open set ${\cal U}_m$.
Let us make a change of variables of the form 
$$
\tilde \Phi := \begin{pmatrix}y\\ v\end{pmatrix}=\begin{pmatrix}x+U(x,u)\\ u+\pi(x+U(x,u))-\pi(x)\end{pmatrix}=:\begin{pmatrix}x\\ u\end{pmatrix}+ W.
$$
On the one hand, we have
$$
\tilde \Phi_*\tilde X = \tilde X + [W,\tilde X]+ \frac{1}{2}[W, [W,\tilde X]]+\cdots.
$$
On the other hand, $\tilde r'_{\Sigma}=\tilde \Phi_*\tilde r_{\Sigma}$ vanishes on $\Sigma$ as well as
\begin{eqnarray*}
pr_1\left([\widetilde{NF^m}, W]\right) & = & [s_1(x,u), W_x] -\sum_{j=1}^l W_u(a_j(u))S_j+ [NF^m-s_1, W_x]\\
& &  + \sum_{i=1}^n \pi_*NF^{m}(W_{x,i})\frac{\partial}{\partial x_i}-\sum_{i=1}^n W_u(NF^{m}_i)\frac{\partial}{\partial x_i}.
\end{eqnarray*}
Here, $W_u$ (resp. $W_x$) stands for $\pi(x+U(x,u))-\pi(x)$ (resp. $U$). Let
$B_{m+1}$ denote the sum of the projection of $R_{m+1}$ onto the weight spaces associated to
weights of degree $m+1$. {\bf Its restriction to $\Sigma$ doesn't depend on the choice of the representative $R_{m+1}$}. Therefore, if $U$ is of order greater than or equal to $m+1$, we may write
$$
pr_1\left(\tilde \Phi_*\tilde X (y,v)\right)= NF^m(y,v)+B^{m+1}(y,v)-[s_1(y,v), U]+R_{m+2}(y,v)+r'_{\Sigma}(y,v),
$$
where the vector field $R_{m+2}(y,v)$ is of order greater than or equal to $m+2$.

For any weight $\alpha$ of $S$ in $\hvf n {m+1}$, let us set 
$$
A_{\alpha}:=\sum_{j=1}^la_j(u)\alpha(g_j).
$$
Since $s_1$ is nondegenerate, if $\alpha$ is nonzero then $A_{\alpha}$
doesn't vanish identically. For if it did then the image of map $(a_1,\ldots,a_l)$ would be contained in a complex hyperplane. Let ${\cal U}_{m+1}$ be the complement, in  ${\cal U}_{m}$, of the zero set of the $A_{\alpha}$'s where $\alpha$ ranges over the set of nonzero weights of $S$ in $\hvf n {m+1}$. For such an $\alpha$, let us set, for $u\in  {\cal U}_{m+1}$,
$$
U_{\alpha}(x,u):=\frac{B_{\alpha}^{m+1}(x,u)}{A_{\alpha}(u)}\quad\text{and} \quad U=\sum_{\alpha\neq 0}U_{\alpha}.
$$
We have 
$$
[s_1(x,u),U_{\alpha}]=\sum_{j=1}^la_j(u)[S_j,U_{\alpha}]=\left(\sum_{j=1}^la_j(u)\alpha(g_j)\right)U_{\alpha}=B_{\alpha}^{m+1}.
$$
Hence, we have
$$
pr_1\left(\tilde \Phi_*\tilde X (y,v)\right)= NF^m(y,v)+B^{m+1}_0(y,v)+R_{m+2}(y,v)+r'_{\Sigma}(y,v).
$$
Let us set 
$$
NF^{m+1}(y,v):=NF^m(y,v)+B^{m+1}_{0,|\Sigma}(y,v).
$$
Then, we have
$$
\left[s_1,NF^{m+1}\right]=0.
$$

Let $\tilde X'$ be another fibered vector field with respect to $X$ and which
is normalized up to order $m$. Let us assume that it differs from $pr_1(\tilde X)$ by a
vector fields vanishing on $\Sigma$. Hence, it can be written as 
$$
pr_1\tilde X' = NF^m(x,u)+R_{m+1}(x,u)+r'_{\Sigma}(x,u).
$$
Then, $\tilde \Phi_*\tilde X'$ will differ from $\tilde \Phi_*\tilde X$ by a
fibered vector field vanishing on $\Sigma$.

\qed\end{proof}
\begin{defis}
The vector field $\tilde X$ is said to be a {\bf good deformation relative to
  $S$} if  $\tilde X$ admits a Lindstedt-Poincar\'e normal form of order $m$ of the form $$\sum_{j=1}^la_{j}^m(u)S_j,$$  for any order $m$ greater than or equal to $2$.
\end{defis}
\begin{rems}
Let us expand the meaning of this definition. In coordinates, the vector field $(\tilde \Phi_m)_*X$ can be written as
\begin{eqnarray*}
\frac{d x_i}{dt} & = & \left(\sum_{j=1}^l a_j^m(u)\lambda_{j,i}\right)x_i + R_{m+1,i}(x,u)+r_{\Sigma,i}\\
\frac{d u_k}{dt} & = & R_{m+1,i}(x^{R_k})+r_{\Sigma,i}(x^{R_k}),
\end{eqnarray*}
where $i$ ranges from $1$ to $n$ and $k$ ranges from $1$ to $p$.
It is a perturbation of order greater than or equal to $m+1$ (in $x$) of the ``integrable"-one along $\Sigma$
\begin{eqnarray*}
\frac{d x_i}{dt} & = & \left(\sum_{j=1}^l a_j^m(u)\lambda_{j,i}\right)x_i +r_{\Sigma,i}\\
\frac{d u_k}{dt} & = & r_{\Sigma,i}(x^{R_k})
\end{eqnarray*}
where $i$ ranges from $1$ to $n$ and $k$ ranges from $1$ to $p$.
The $u_k$'s are to be thought as the actions (the ``slow variables") whereas the nonresonant monomials $x^Q$ are the functions of the ``fast variables".
\end{rems}

%


\section{Main results}

Let $\lie g$ be a complex $l$-dimensional commutative Lie algebra. Let $S:\lie g\rightarrow \hvf n 1$ be a Lie morphism from $\lie g$ to the Lie algebra of linear vector fields of ${\bf C}^n$. It is assumed to be injective and  semi-simple. This means that, up to a linear change of coordinates, there are linear forms $\lambda_1,\ldots, \lambda_n\in \lie g^*$ such that for all $g\in \lie g$, 
$$
S(g)= \sum_{i=1}^n{\lambda_i(g)x_i\frac{\partial}{\partial x_i}}.
$$ 
In the sequel, $\{g_1,\ldots, g_l\}$ will denote a fixed basis of $\lie g$ and we shall set $S_i=S(g_i)$. The family $\{S_i\}_{i=1,\ldots,l}$ is linearly independent over ${\bf C}$. We shall denote by ${\cal W}_{n,*}^{k,m}$ the set 
of nonzero weights of $S$ in $\pvf n k m$; that is the set set of nonzero
linear forms $\sum_{i=1}^nq_j\lambda_j(g)-\lambda_i(g)$, for which $(q_1,\ldots,
q_n)\in {\bf N}^n$, $1\leq i\leq n$ and $k\leq q_1+\cdots+q_n\leq m$. We {\bf
  assume} that the ring of formal first integrals $\widehat{{\cal O}_n^S}$ is
not reduced to ${\bf C}$. We recall that $\widehat{{\cal O}_n^S}=\Bbb
C[[x^{R_1},\ldots, x^{R_p}]]$ where $x^{R_1},\ldots, x^{R_p}$ are monomials of ${\bf C}^n$. We
shall assume that they are {\bf algebraically independent}.

Let $X\in \vfo n 1$ be a germ of vector field of $({\bf C}^n,0)$ vanishing at
the origin and which linear part belongs to $S(\lie g)$. Let us assume that
$X$ is a {\bf good perturbation} of order $m_0+1\geq 2$ of a {\bf nondegenerate} vector field 
$$X_0=\sum_{j=1}^la_jS_j$$ where the $a_j$'s belong to ${\cal O}_n^S$. We shall write $a=(a_1,\ldots, a_l)$. Hence, its Lindstedt-Poincar\'e normal form, at each order $m\geq m_0$, is of the form
$$NF^m(x,u)=\sum_{j=1}^la_j^m(u)S_j(x),$$ where the $a_j^m$'s are holomorphic
functions on some open set $U_m$ of ${\bf C}^p$. Moreover, $X-X_0$ is flat up to order $m_0$
at the origin. 

Let $\omega=\{\omega_k\}_{k\in {\bf N}^*}$ be a sequence a positive numbers such that 
\begin{itemize}
\item $\omega_k\leq 1$,
\item $\omega_{k+1}\leq \omega_k$,
\item the series $\sum_{k>0}\frac{-\ln \omega_k}{2^k}$ converges.
\end{itemize}
Such a sequence will be called a {\bf diophantine sequence}.

Let $\rho$ be a sufficiently small positive number less than $1/2$. Let ${\cal K}$ be a nonvoid compact set of $\pi(D_n(0,\rho))$.
Let $\gamma$ be a positive real number and less than some $\gamma'$. We define the decreasing sequence $\{{\cal K}_k(X,{\cal K},\omega,\gamma)\}_{k\in {\bf N}}$ of compact sets of $\pi(D_n(0,\rho))$ as follows:
\begin{eqnarray*}
{\cal K}_0 & = & {\cal K},\\
{\cal K}_k & = & \left\{b\in  {\cal K}_{k-1}\,|\,\forall \alpha\in {\cal W}_{n,*}^{2^k+1,2^{k+1}},\;\;\left|\alpha\left(\sum_{j=1}^l{a_j^{2^k}(b)g_j}\right)\right|\geq \gamma \omega_{k+1}\right\}.
\end{eqnarray*}

\begin{theos}\label{theo1}
Under the assumptions above, let $\omega$ be a diophantine sequence; let
${\cal K}$ and $\gamma$ be defined as above. If $m_0$ is large enough and if the set ${\cal K}_{\infty}(X,{\cal K},\omega,\gamma):=\bigcap_{k\in {\bf N}^n} {\cal K}_k$ is nonvoid then, for any $b\in {\cal K}_{\infty}$, 
\begin{enumerate}
\item the sequence $\{NF^m(x,b)\}$ converges to a linear diagonal vector field
  \\$NF(x,b)=\sum_{j=1}^l\tilde a_j(b)S_j(x)$, 
\item there is a biholomorphism of analytic subsets of open sets in ${\bf C}^n$, $$\Theta_b:\pi^{-1}(b)\cap D_0(\rho)\rightarrow V_{b}\subset {\bf C}^n$$ which conjugates the restriction of $NF(x,b)$ to $\pi^{-1}(b)\cap D_0(\rho)$ to the restriction of $X$ to 
$V_{b}$.
\end{enumerate}
As a consequence, $X$ is tangent to the {\bf toroidal} analytic subset
$V_{b}$; its restriction to it is conjugated to the restriction to the toric
analytic subset $\pi^{-1}(b)\cap D_n(0,\rho)$ of the {\bf linear diagonal vector field} $\sum_{j=1}^l\tilde a_j(b)S_j$.
\end{theos}

\begin{defis}
Let $\omega=\{\omega_k\}_{k\geq 1}$ be a diophantine sequence and $\mu_0$ be a positive integer. We shall say that 
$S$ is {\bf strictly diophantine relatively to $(\omega,\mu_0)$} if 
$$
\lim_{k\rightarrow+\infty}\left(2^k+n+1\right)^{n+1}\left(\frac{\omega_k}{\omega_k(S}\right)^{2/\mu_0}=0.
$$
\end{defis}
\begin{theos}\label{theo3}
Let ${\cal K}$ be a compact set of $\pi(D_n(0,\rho))$ of positive $2p$-measure. Assume that $S$ is strictly diophantine relatively to the sequence $(\omega=\{\omega_i\}_{i\geq
  1},\mu_0)$ where $\mu_0$ is the index of nondegeneracy of $X_0$ with respect to ${\cal K}$ (see section 9.1). Then, under the assumptions of theorem \ref{theo1}, ${\cal K}_{\infty}$ is nonvoid and has a positive $2p$-measure.
\end{theos}

The definition of strict diophantineness is derived from R\"ussmann assumptions
in his work on KAM theory in the symplectic case \cite{russmann-weak}. It should be noticed that, in the hamiltonian situation, the numbers $\omega_k(S)$ don't appear in the condition since they are bounded from below (see for instance \cite{stolo-ihes} [theorem 10.2.1-10.2.2]).

\begin{rems}
Similar statements were announced in \cite{Stolo-kam-cras}. At that time,
we made the assumption that the formal Poincar\'e-Dulac normal form of the
deformation $X$ was of the form $\sum_j\hat a^jS_j$, with $\hat a_j\in
\widehat{\cal O}_n^S$. But, we failed to prove that this assumption implies
the assumption we actually do on the Lindsedt-Poincar\'e normal form. Hence,
the result about the holomorphy of the normal form is not proved yet.
\end{rems}

\section{Sketch of the proof}

Let $\tilde X$ be a fibered deformation of $X$ over a neighborhood of the origin in ${\bf C}^p$.
The global mechanism of our proof is to normalize, in the sense of Lindstedt-Poincar\'e, the vector field $\tilde X$. To do so, we shall use a slightly modified Newton process. Assuming that $\tilde X$ is normalized up to order $m$, we shall normalize it up to order $2m$ by a holomorphic diffeomorphism of some polydisc $D_n(0,R_m)\times D_p(b, t_m)$, tangent to identity at $(0,b)$. It is also fibered over $\pi$. 
Let us set 
$$
\tilde X(x,u)=NF^m+\tilde R_{m+1}+\tilde r_{\Sigma}
$$
where $NF^m(x,u)=\sum_{i=1}^la_i^m(u)S_i(x)$, $\tilde R_{m+1}$ is a fibered over a holomorphic vector field of order greater than or equal to $m+1$ and $\tilde r_{\Sigma}$ is fibered over a vector field vanishing on $\Sigma$. Let us write 
$$
\tilde R_{m+1}=\tilde B^{m+1,2m}+\tilde R_{2m+1}
$$
where $\tilde B^{m+1,2m}$ denotes the projection of $\tilde R_{m+1}$ on the weight spaces (of $S$) of degree less than or equal to $2m$:
$$
\tilde B^{m+1,2m}=\tilde B_0^{m+1,2m}\oplus\sum_{\alpha\neq 0}\tilde B_{\alpha}^{m+1,2m}
$$
where the sum is taken over the nonzero weights of the representation of $S$ into the space ${\cal O}_p(D_p(b, t_m))\otimes\pvf n {m+1} {2m}$, $B_{\alpha}^{m+1,2m}$ being the projection of $B^{m+1,2m}$ onto the $\alpha$-weight space. 

In order to normalize $X$ up to order $2m$, one is led to solve, for each nonzero weight $\alpha$ of order less than or equal to $2m$, the following {\bf cohomological equation}:
$$
[NF^m, U_{m,\alpha}](x,u)+D_uNF_m(x,u)D\pi(x)U_{m,\alpha}(x,u)= B_{\alpha}^{m+1,2m}(x,u),
$$
where the unknown, $U_{m,\alpha}$, is to be a holomorphic vector field of order greater than or equal to $m+1$, belonging to the $\alpha$-weight space. Let us set 
$$
A_{m,\alpha}(u)=\sum_{j=1}^la_j^{m}(u)\alpha(g_j).
$$
Let $V_m$ be an open set where $A_{m,\alpha}$ doesn't vanish and assume that $b$ belongs to $V_m$. On this set, we find that 
$$
U_{m,\alpha}(x,u) = \left(Id - \frac{D_m}{A_{m,\alpha}(u)}\right)\frac{B_{\alpha}^{m+1,2m}(x,u)}{A_{m,\alpha}(u)}
$$
is the solution of the cohomological equation. The ${\cal O}_n\times{\cal O}_p(\{b\})$-linear operator $D_m$ is nilpotent of order $2$ on ${\cal O}_p(\{b\})\otimes\vfo n {m+1}$.

Let us set 
$$
U_m=\sum_{\alpha\neq 0 \text{ and of degree }\leq 2m}U_{m,\alpha},
$$
and
$$
\Phi_m=Id +U_m, \quad\tilde \Phi_m(x,u)=\left(x+U_m(x,u), u+\pi(x+U_m(x,u))-\pi(x)\right).
$$
This map is a holomorphic diffeomorphism in a neighborhood of $(0,b)$ and it normalizes $\tilde X$ up to order $2m$. In fact, let us set $NF^{2m}=NF^{m}+\tilde B^{m+1,2m}_{0|\Sigma}$. Then, 
$$
(\tilde\Phi_m)^*\tilde X=NF^{2m}+\tilde R_{2m+1}+r_{\Sigma}
$$
is normalized up to order $2m$ ($\tilde R_{2m+1}$ is of order greater than or equal to $2m+1$ in $x$).

We have to control the behavior of the estimates of $\tilde \Phi_m$ and $(\tilde\Phi_m)_*\tilde X$ when $m$ ranges from some integer on. This is the goal of the section entitled ``The induction process".
In order to obtain good estimates we shall consider special sets.
Let us set $m=2^k$ and le $b$ belong to 
$$
{\cal K}_k:=\left\{c\in  {\cal K}_{k-1}\,|\,\forall \alpha\in {\cal W}_{n,*}^{2^k+1,2^{k+1}}\left|\alpha\left(\sum_{j=1}^l{a_j^{2^k-1}(c)g_j}\right)\right|\geq \gamma \omega_{k+1}\right\}.
$$
Let use set
$$
t_m := \gamma\frac{\omega_{k+1}}{2l\Lambda(2m+1)}.
$$

Let $U$ be an open connected set of ${\bf C}^p$ and let $r$ be a positive number. Let $f$ belong to ${\cal O}_p(U)\otimes {\bf C}[[x_1,\ldots, x_n]]$:
$$
f=\sum_{Q\in {\bf N}^n}f_Q(u) x^Q
$$
where the $f_Q$'s are holomorphic functions on $U$. We set 
$$
|f|_{U,r}=\sum_{Q\in {\bf N}^n}\sup_{u\in U}|f_Q(u)| r^{|Q|}.
$$
If $f$ is just a formal power series, then $|f|_{U,r}$ doesn't depend on $U$ and is nothing but the polydisc norm
$$
|f|_{r}=\sum_{Q\in {\bf N}^n}|f_Q| r^{|Q|}.
$$

We show that the function $A_{m,\alpha}(u)$ doesn't vanish on $D_n(0,R_m)\times D_p(b, t_m)$. Therefore, $U_{m,\alpha}$ can be seen as a holomorphic vector field which coefficients are holomorphic functions on $D_n(0,R_m)\times D_p(b, t_m)$. Moreover, it is of order greater than or equal to $m+1$ in $x$. We show the following estimate (proposition \ref{cohom-estimate}): 

{\it let $r>1/2$. If $\|D(a^{m}(u))\|_{D_p(b, t_m)}\leq 1$, there exists $c_1>0$ such that the solution of the cohomological equation satisfies
$$
|U_{m,\alpha}|_{D_p(b, t_m),r} \leq \frac{c_1}{\gamma^2\omega_{k+1}^2}|B_{\alpha}^{m+1,2m}|_{D_p(b, t_m),r}.
$$}

We should emphasize the r\^ole of the ring of invariants ${\cal O}_n^S$, the element of which are the natural ``constants". The small divisors are no longer complex numbers but rather elements of this ring: these are the functions $A_{m,\alpha}$ defined on an appropriate set.

We assume that $1/2< r\leq 1$ and that $m=2^k$ for some positive integer $k$. We define the positive numbers
$$
\gamma_k=\left(\frac{c_1}{\gamma^2\omega_{k+1}^2}\right)^{-1/m},\quad\theta_k := \gamma_km^{-2/m},\quad r_i:=\theta_k^ir,\;\;i\geq 1.
$$
If $m$ is sufficiently large, we set
\begin{eqnarray*}
{\cal NF}_{m,b}(r)& = & \left\{X\in {\cal O}_p(D_p(b,t_{m/2}))\otimes\hvf n {1} \;|\; |X|_{D_p(b,t_m),r}<1-\frac{1}{m^3},\right.\\
& & \left. |D_u(X)|_{D_p(b,t_m),r}<1-\frac{1}{m^2}\right\},\\
{\cal B}_{m+1,b}(r) & = & \left\{X\in {\cal O}_p(D_p(b,t_m))\otimes\vfo n {m+1}\;|\; |X|_{D_p(b,t_m),r}<\frac{2^5n}{m^4}\right\}.\\
\end{eqnarray*}
We assume that ${2^5n}/{m^4}<1$. We shall show the following (proposition \ref{reccurence}):

{\it let $b$ belong to ${\cal K}_{k-1}$. We assume $\tilde X=NF^m+\tilde R_{m+1}+\tilde r_{\Sigma}$ is normalized up to order $m$ and that $(NF^m,R_{m+1})\in {\cal NF}_{m,b}(r)\times {\cal B}_{m+1,b}(r)$. If $b$ belongs to ${\cal K}_k$ and if $m$ is sufficiently large (say $\geq 2^{k_0}$ independent of $r$), then 
\begin{itemize}
\item $D_n(0,r_5)\times D_p(b, t_{2m})\subset \tilde \Phi_m(D_n(0,r)\times D_p(b, t_m))$,
\item $(NF^{2m}, R_{2m+1})\in {\cal NF}_{2m,b}(r_5)\times {\cal B}_{2m+1,b}(r_5)$.
\end{itemize}}

This proposition is fundamental for the induction process since it enables us to control the norms. Moreover, it states that, at each step, there is a ``good" thickening neighborhood of the toric variety on which the analysis can be done. Its proof is rather long and technical.

Now, we are able to give a sketch of the proof the existence of the invariants varieties. Since the sequence $\{\omega_k\}$ is diophantine, we show that the sequence of radii defined by 
$$
R_{k+1}=\theta_k^5R_k,
$$
converges to some positive $R$ which can be assume to be greater than $1/2$. Let $\{\tilde\Psi_k\}$ be the sequence of holomorphic diffeomorphisms defined by
$$
\tilde \Psi_k=\tilde \Phi_{2^k}\circ \tilde \Psi_{k-1}.
$$
The diffeomorphism $\tilde \Psi_k$ normalizes $\tilde X$ up to order $2^{k+1}$. Thanks to our estimates, we show that the sequence of inverse diffeomorphisms $\{\Psi_k^{-1}\}$, when restricted to $D_n(0,\rho)\times\pi^{-1}(b)$, converges (for the compact topology of the analytic set) to a holomorphic
diffeomophism on $\Theta_b$ on $D_n(0,\rho)\times\pi^{-1}(b)$ (for some well
chosen $\rho<1/2$).

What about the conjugacy problem? By construction, both $\tilde \Psi_m$ and $\tilde X$ are tangent to $\Sigma$. Hence, according to the induction process, we show that the sequence $\{NF^m\}$, restricted to $(D_n(0,\rho)\cap \pi^{-1}(b))\times\pi^{-1}(b)$, converges to a linear vector field $$NF_b(x)=\sum_{j=1}^la_j(b)S_j(x).$$ By definition, it is tangent to the {\bf toric variety} $\pi^{-1}(b)\cap D_n(0,\rho)$. Moreover, $\Theta_b$ conjugates its restriction to $\pi^{-1}(b)\cap D_n(0,\rho)$ to the restriction of $X$ to the analytic set $\Theta_b(\pi^{-1}(b)\cap D_n(0,\rho))$. That is, $X$ has an invariant analytic subset, namely $\Theta_b(\pi^{-1}(b)\cap D_n(0,\rho))$, which is biholomorphic to $\pi^{-1}(b)\cap D_n(0,\rho)$. This ends the sketch of the proof of the first theorem.

The proof of the last theorem is an adaptation of R\"ussmann work on KAM theory \cite{russmann-weak}[part 4]. Its goal is to give a sufficient condition which will ensure that the compact set ${\cal K}_{\infty}$ is not empty. To do so, we shall show that our assumptions are sufficient to ensure that each compact set ${\cal K}_k$ has a positive measure. The main tool is R\"ussmann theorem \ref{mes-epsilon} which can rephrased as follows: given a compact set ${\cal K}$ in ${\bf R}^n$, a neighborhood ${\cal B}$ of ${\cal K}$ and a real valued $C^{\mu_0+1}$-function $g$ on ${\cal B}$. Assume roughly that $g$ as well as all its $\mu_0$ derivatives do not vanish simultaneously on ${\cal K}$. Then, there is an upper bound for the measure of the inverse image of $]-\epsilon,\epsilon[$ by any small $C^{\mu_0}$-perturbation of $g$ defined on ${\cal B}$; $\epsilon$ has to be small and the upper bound depends only on $g$, n, and the size of ${\cal B}$ relatively to ${\cal K}$. Its proof is based on Bakhtin theorem \cite{bakhtin}.


\section{Solution of the cohomological equation}

Let us set 
$$
\tilde X_m:=\tilde \Psi^m_*\tilde X \in {\cal O}_p({\cal U}_m)\otimes \vfo n 1 (D_n(0, R_m))
$$ 
which is assumed to be normalized up to order $m$. Let us assume that its Lindstedt-Poincar\'e normal form of order $m$ can be written as 
$$
NF^m(x,u)=\sum_{i=1}^la_i^{m}(u)S_i
$$
with the $a_i^{m}$'s belong to ${\cal O}_p({\cal U}_m)$. Hence, we have 
$$
pr_1(X_m)= NF^m(x,u)+R_{m+1}(x,u)+r_{\Sigma}(x,u),
$$ 
where $R_{m+1}$ belongs to ${\cal O}_p({\cal U}_m)\otimes\vfo n {m+1}$ and where $r_{\Sigma}$ denotes a vector field vanishing on the subvariety $\Sigma$. 

Let us decompose $R_{m+1}$ along the weight spaces of $S$ and let us write
$R_{m+1}=B^{m+1,2m}+ R'_{2m+1}$ where $B^{m+1,2m}(x,u)$ denotes the sum of the
projections of $R_{m+1}$ along the weight spaces associated to a weight of degree less
than or equal to $2m$. We should emphasize that $B^{m+1,2m}(x,u)$ may not be a
polynomial (in $x$).
The vector field $R'_{2m+1}$ is of order greater than or equal to $2m+1$, with coefficients in ${\cal O}_p({\cal U}_m)$. 

\subsection{Cohomological equations }

Let $x=\Phi_m(y,v)=y+U_m(y,v)$ be a family of holomorphic diffeomorphisms of
$({\bf C}^n,0)$ where $U_m\in {\cal O}_p({\cal U'}_m)\otimes\vfo n {m+1}$. Let $\tilde
\Phi_m$ be the fibered diffeomorphism over $\pi$ associated to
$\Phi_m$. Here, ${\cal U'}_m$ denotes an open set such that
$$\tilde\Phi_m(D_n(0, R_m)\times {\cal U'}_m)\subset D_n(0, R_m)\times {\cal
  U}_m.$$ We shall write
$(x,u)=\tilde \Phi_m(y,v)$. 
We have 
$$
D\left(\tilde \Phi_m^{-1}\right)(y,v)(\tilde\Phi_m)_*\tilde X(y,v) = \tilde X\left(\tilde \Phi_m^{-1}(y,v)\right)
$$
Let us write 
\begin{eqnarray*}
pr_1\left((\tilde\Phi_m)_*\tilde X(y,v)\right) &=&
NF^m(y,v)+B'(y,v)+C'(y,v)+r'_{\Sigma}\\
&=:& Z(y,v)+r'_{\Sigma}
\end{eqnarray*}
where $B'$ belongs to ${\cal O}_p({\cal U}'_m)\otimes\left(\vfo n {m+1}\right)^S$ and $C'$ belongs to $
{\cal O}_p({\cal U}'_m)\otimes\vfo n {2m+1}$. By assumption, we have
$D\pi(x)NF_m(x,u)=0$, hence we have the following relations:
\begin{eqnarray*}
(Id +D_yU_m(y,v))Z(y,v)& & \\
+D_vU_m(y,v)D\pi(y)(B'(y,v)+C'(y,v))& = & NF^m(x,u)+B(x,u)+C(x,u)\\
& = & (NF^m+B)(y,v)\\ 
& & +D(NF^m)(y,v)(\tilde \Phi_m^{-1}(y,v)-(y,v))\\
& & +(B(\tilde\Phi_m^{-1}(y,v))-B(y,v))\\
& & +\left(NF^m(\tilde \Phi_m^{-1}(y,v))-NF^m(y,v)\right.\\
& & \left.-D(NF^m)(y,v)(\tilde \Phi_m^{-1}(y,v)-(y,v))\right)\\
& & -pr_1\left(D\left(\tilde \Phi_m^{-1}\right)(y,v)\tilde r'_{\Sigma}(y,v)\right).
\end{eqnarray*}

Therefore, we have
\begin{eqnarray}
C'(y,v) + B'(y,v)-B(y,v)& & \label{action-diff}\\
-[NF^m,U_m](y,v)&& \nonumber\\
-D_uNF_m(y,v)D\pi(y)U_m(y,v) &= &  -D_y(U_m)(y,v)(B'(y,v)+C'(y,v))\nonumber\\
&& - C'(y,v)\nonumber\\
&& +(B(\tilde\Phi_m^{-1}(y,v))-B(y,v))\nonumber\\
&& \left(NF^m(\tilde\Phi_m^{-1}(y,v))-NF^m(y,v)\right.\nonumber\\
&& \left.-D(NF^m)(y,v)(\tilde\Phi_m^{-1}(y,v)-(y,v))\right)\nonumber\\
&& +D(NF^m)(y,v)(0, \pi(y+U_m)-\pi(y)\nonumber\\
&& -D\pi(y)U_m)\nonumber\\
&& -D_vU_m(y,v)D\pi(y)(B'(y,v)+C'(y,v))\nonumber\\
&& -pr_1\left(D\left(\tilde \Phi_m^{-1}\right)(y,v)\tilde r'_{\Sigma}(y,v)\right).\nonumber
\end{eqnarray}

Assuming that both $U_m$, $B$ and $B'$ (resp. $C'$) are of order greater than or equal to $m+1$
(resp. $2m+1$), it is straightforward to notice that the right hand side of equation $(\ref{action-diff})$ is of order greater than or equal to $2m+1$ modulo a vector field $r'''_{\Sigma}$ vanishing on $\Sigma$.

Let $\alpha$ be a nonzero weight of $S$ of degree less than or equal to $2m$. Let us find a
solution $U_{m,\alpha}$ belonging to ${\cal O}_p({\cal U}_m)\otimes\vfo {n,\alpha} {m+1}$ of the {\bf cohomological equation} 
$$
[NF^m,U_{m,\alpha}](x,u)+D_uNF_m(x,u)D\pi(x)U_{m,\alpha}(x,u)=B_{\alpha}^{m+1,2m}(x,u).
$$
Here $B_{\alpha}^{m+1,2m}$ denotes the projection of $B^{m+1,2m}$ onto the $\alpha$-weight space.
Let us write the cohomological equation as follows:
\begin{eqnarray}
[NF^{m}, U_{m,\alpha}](x,u) & = & \sum_{j=1}^la_j^{m}(u)[S_j,U_{m,\alpha}]\nonumber\\
& =& {\left(\sum_{j=1}^la_j^{m-1}(u)\alpha(g_j)\right)}U_{m,\alpha}(x,u).\label{equ-modifie}
\end{eqnarray}
Let us set 
\begin{eqnarray*}
A_{m,\alpha}(u)& := & \sum_{j=1}^la_j^{m}(u)\alpha(g_j),\\
D_m(U_m) & := & D_uNF_m(x,u)D\pi(x)U_m(x,u)=\sum_{j=1}^lD_ua_j^{m}(x,u)D\pi(x)U_m(x,u)S_j,\\
{\cal U}_m''& := &{\cal U}_m'\setminus \cup_{\alpha\in {\cal
    W}_{n,*}^{m+1,2m}} \{u\in {\cal U}'_m\;|\;A_{m,\alpha}(u)=0\}.
\end{eqnarray*}
The operator $D_m$ is nilpotent with $D_m\circ D_m=0$. In fact, since the Lie derivative of (each component of) $\pi$ along the $S_i$'s vanishes, we have \\$D\pi(x)D_m(U_m)=0$. 


Let us set, on $D_n(0,R_m)\times {\cal U}''_m$
\begin{eqnarray}
U_{m,\alpha}(x,u) & := & \left(Id - \frac{D_m}{A_{m,\alpha}(u)}\right)\frac{B_{\alpha}^{m+1,2m}(x,u)}{A_{m,\alpha}(u)}.\label{equ-sol}
\end{eqnarray}
This vector field belongs to ${\cal O}_p({\cal U}''_m)\otimes\vfo {n,\alpha} {m+1}(D_n(0,R_m))$. First of all, since $B_{\alpha}^{m+1,2m}$ is of order greater than or equal to $m+1$, so is $D_m(B_{\alpha}^{m+1,2m})$. 
Since $A_{m,\alpha}(u)$ is a function of $u$, $U_{m,\alpha}$ is of order greater than or equal to $m+1$. 
%
%

According to $(\ref{equ-modifie})$ and the properties of the operator $D_m$, we have
\begin{eqnarray*}
[NF^m,U_{m,\alpha}]+D_uNF_m(x,u)D\pi(x)U_{m,\alpha}(x,u) & = & A_{m,\alpha}(x,u)U_{m,\alpha}+D_m(U_{m,\alpha})\\
& = & \left(\frac{B_{\alpha}^{m+1,2m}}{A_{m,\alpha}(u)} - \frac{D_m(B_{\alpha}^{m+1,2m})}{A^2_{m,\alpha}(u)}\right)\\
&&+ D_m\left(\frac{B_{\alpha}^{m+1,2m}}{A_{m,\alpha}(u)}
- \frac{D_m(B_{\alpha}^{m+1,2m})}{A_{m,\alpha}^2(u)}\right)\\
& = & B_{\alpha}^{m+1,2m}.
\end{eqnarray*}
Let us set $U_m:=\sum_{\alpha\neq 0}U_{m,\alpha}$, the sum being taken over the set of nonzero weights of $S$ into $\pvf n {m+1} {2m}$. Then, we have 
$$
[NF^m, U_m]+D_uNF_m(x,u)D\pi(x)U_{m}(x,u)=B_{+}^{m+1,2m},
$$
where $B_{+}^{m+1,2m}$ denotes $B^{m+1,2m}-B_{0}^{m+1,2m}$.
Let us set $$NF^{2m}(x,u):= NF^m(x,u)+B^{m+1,2m}_{0|\Sigma}(x,u).$$ Then, $(\tilde \Phi_m)_*\tilde X_m(y,v)=NF^{2m}(y,v)+\tilde C'(y,v)+\tilde r'_{\Sigma}(y,v)$ is normalized up to order $2m$. Here, $\tilde C'(y,v)$ stands for a fibered vector field over a vector field of order greater than or equal to $2m+1$. The vector field $\tilde r'_{\Sigma}$ is defined to be 
$$
(\tilde \Phi_m)_*\tilde r_{\Sigma}+ \tilde B^{m+1,2m}_{0}-\tilde B^{m+1,2m}_{0|\Sigma},
$$
and vanishes on ${\Sigma}$.

\subsection{Estimate for the solution of the cohomological equation}

We recall that $$NF^m(x,u)=\sum_{j=1}^la^{m}(u)S_j$$ where the $a^{m}_j$'s belongs to ${\cal O}_p({\cal U}_m)$ and where we have $m=2^k$.
By assumption, ${\cal K}_{k-1}$ is a nonvoid compact set of $\pi(D_n(0,\rho))$. 
Let us define the compact set on $\pi(D_n(0,\rho))$
$$
{\cal K}_{k}=\left\{b\in  {\cal K}_{k-1}\,|\,\forall \alpha\in {\cal W}_{n,*}^{2^k+1,2^{k+1}}\left|\alpha\left(\sum_{j=1}^l a^{m}_j(b)g_j\right)\right|\geq \gamma \omega_{k+1}\right\}.
$$
Let us {\bf assume} that ${\cal K}_{k}\neq \emptyset$. Let us set 
$$
\Lambda:=\max_{1\leq i\leq n,1\leq j\leq l}|\lambda_i(g_j)|. 
$$
%

\begin{figure}[hbtp]
  \begin{center}
    \leavevmode
    \input{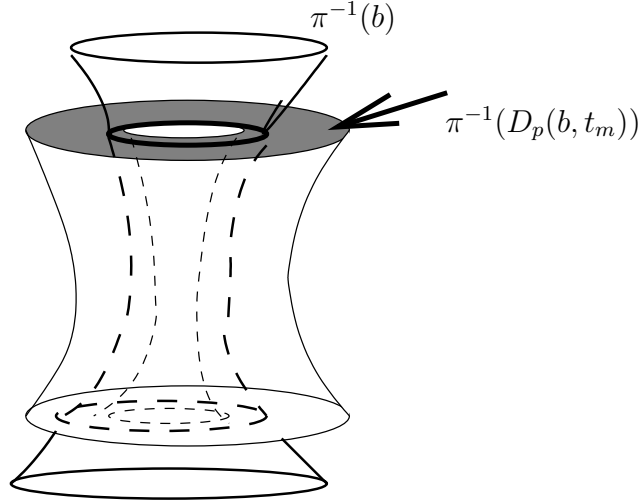} 
    \caption{The domain $\pi^{-1}(D_p(b,t_m))$ is a thickening neighborhood of the
    fiber $\pi^{-1}(b)$}
  \end{center}
\end{figure}

\begin{props}\label{cohom-estimate}
Let $r>1/2$, and $b\in{\cal K}_{k}$ which assumed to be nonvoid. Let us set 
$$
t_m := \gamma\frac{\omega_{k+1}}{2l\Lambda(2m+1)}.
$$
If $\|D(a^{m}(u))\|_{D_p(b,t_m)}\leq 1$, then there exists a positive number $c_1$ such that, for
any nonzero weight $\alpha$ of $S$ in $\pvf n {m+1} {2m}$, the solution of the cohomological equation $(\ref{equ-sol})$ satisfies
$$
|U_{m,\alpha}|_{D_p(b,t_m),r} \leq \frac{c_1}{\gamma^2\omega_{k+1}^2}|B_{\alpha}^{m+1,2m}|_{D_p(b,t_m),r}.
$$
\end{props}
\begin{proof}
Let $\alpha$ be a nonzero weight of $S$ into $\pvf n {m+1} {2m}$. Then there exists 
a multiindex $Q=(q_1,\ldots, q_n)\in {\bf N}^n$, $m+1\leq |Q|\leq 2m$ and an index $1\leq i\leq n$ such that $$\alpha(g)=\alpha_{Q,i}(g)=\sum_{j=1}^nq_j\lambda_j(g)-\lambda_i(g).$$
On the set $D_p(b,t_m)$, we have the following estimate:
\begin{eqnarray*}
\left|\sum_{j=1}^l(a_j^{m}(u)-a^{m}_j(b))\alpha(g_j)\right| & \leq & \gamma\frac{\omega_{k+1}}{4\Lambda(2m+1)}\max_j|\alpha(g_j)|.
\end{eqnarray*}
Since $|\alpha(g_j)|\leq (2m+1)\max_i|\lambda_i(g_j)|$, then we have, on $D_p(b,t_m)$,
\begin{eqnarray}
\left|\sum_{j=1}^l(a_j^{m}(u)-a^{m}_j(b))\alpha(g_j)\right| & \leq & \gamma\frac{\omega_{k+1}}{2}.\label{estim1}
\end{eqnarray}
Since $b\in {\cal K}_{k}$, then 
$$
\left|\sum_{j=1}^l a_j^{m}(b)\alpha(g_j)\right|=\left|\alpha\left(\sum_{j=1}^la^{m}_j(b)g_j\right)\right|\geq \gamma \omega_{k+1}.
$$
Therefore, on $D_p(b,t_m)$, we have the following estimate
\begin{eqnarray}
|A_{m,\alpha}(u)|& = & \left|\sum_{j=1}^la_j^{m}(u)\alpha(g_j)\right| \nonumber\\
& \geq & \left| \left|\sum_{j=1}^la^{m}_j(b)\alpha(g_j)\right| - \left|\sum_{j=1}^l(a_j^{m}(u)-a_j^{m}(b))\alpha(g_j)\right| \right|\nonumber\\
& \geq & \gamma\frac{\omega_{k+1}}{2}.\label{estim2}
\end{eqnarray}
Let us set $m_r:= \max_{1\leq j\leq l}|S_j|_{r}|D\pi|_r$. Then, we have
$$
|D_uNF_m(x,u)D\pi(x)U_{m,\alpha}(x,u)|_{D_p(b,t_m),r}\leq lpnm_r |U|_{D_p(b,t_m),r}\|D_u(a^{m})(u)\|_{D_p(b,t_m)}.
$$
By equation $(\ref{equ-sol})$, we obtain the following estimate:
$$
|U_{m,\alpha}|_{D_p(b,t_m),r}\leq \frac{4|B_{\alpha}^{m+1,2m}|_{D_p(b,t_m),r}}{\gamma^2\omega_{k+1}^2}\left( \gamma\frac{\omega_{k+1}}{2} 
+ pnlm_r\|D(a^{m}(u)\|_{D_p(b,t_m)}\right).
$$
Since $\omega_{k+1}\leq 1$ and $\gamma\leq\gamma'$ then, 
there exists a positive constant $c_1$ such that
\begin{eqnarray}
|U_{m,\alpha}|_{D_p(b,t_m),r} & \leq & \frac{c_1}{\gamma^2\omega_{k+1}^2}|B_{\alpha}^{m+1,2m}|_{D_p(b,t_m),r}.\label{estim-cohom}
\end{eqnarray}
This ends the proof of the proposition.
\qed\end{proof}
In what follows, we shall set
\begin{eqnarray}
\gamma_k^{-m}& := &\frac{c_1}{\gamma^2\omega_{k+1}^2}, \label{def-gamma}
\end{eqnarray}
and we may assume that $\gamma_k\leq 1$.


\section{The induction process}
%
We assume that $1/2< r\leq 1$ and let $m=2^k$ be an integer greater than or equal to some positive integer $N_0$ ($N_0$ is greater than $1$ and is to be set in the sequel). We define the positive numbers
$$
\gamma_k=\left(\frac{c_1}{\gamma^2\omega_{k+1}^2}\right)^{-1/m},\quad \theta_k := \gamma_km^{-2/m},\quad r_i:=\theta_k^ir,\quad i=1,\ldots 4,
$$

Let us set
\begin{eqnarray*}
{\cal NF}_{m,b}(r)& = & \left\{X\in {\cal O}_p(D_p(b,t_{m/2}))\otimes\hvf n {1} \;|\; |X|_{D_p(b,t_m),r}<1-\frac{1}{m^3},\right.\\
& & \left. |D_u(X)|_{D_p(b,t_m),r}<1/(2l|L^{-1}|)-\frac{1}{m^2}\right\},\\
{\cal B}_{m+1,b}(r) & = & \left\{X\in {\cal O}_p(D_p(b,t_m))\otimes\vfo n {m+1}\;|\; |X|_{D_p(b,t_m),r}<\frac{2^5n}{m^4}\right\}.\\
\end{eqnarray*}
We assume that ${2^5n}/{m^4}$ is less than $1$. This can always be achieved for $m$ is sufficiently large.

{\bf We shall assume that $|b|+t_{2m}\leq (r_3)^{|R_i|}$, for all $i=1,\ldots, p$.}
This will enable us to apply the remark \ref{importante_remarque}. 

The aim of this section is to prove the following result:
\begin{props}\label{reccurence}
With the above notation, let $b\in {\cal K}_{k-1}$. We assume $\tilde X_m$ is normalized up to order $m$ along $\Sigma$. Hence, we have 
$$pr_1(\tilde X_m)= NF^m(x,u)+R_{m+1}+r_{\Sigma}$$ with 
$$NF^m(x,u)=\sum_{j=1}^la_j^m(u)S_j,$$ and where $R_{m+1}$ is of order greater than or equal to $m+1$ and $r_{\Sigma}$ vanishes on $\Sigma$. Let $\tilde \Phi_m$ be the diffeomorphism
which normalizes $X_m$ up to order $2m$ as defined above. Let us write
$$pr_1((\tilde \Phi_m)_*\tilde X_m)= NF^{2m}+ R'_{2m+1} + r'_{\Sigma}.$$

Let us assume that $(NF^m, R_{m+1})\in {\cal NF}_{m,b}(r)\times {\cal B}_{m+1,b}(r)$. If $b\in {\cal K}_k$ and if $m$ is sufficiently large (say $\geq 2^{k_0}$ independent of $r$ and $b$), then 
\begin{enumerate}
\item $D_n(0,r_4)\times D_p(b,t_{2m})\subset \tilde \Phi_m(D_n(0,r)\times D_p(b,t_{m}))$,
\item $(NF^{2m}, R'_{2m+1})\in {\cal NF}_{2m,b}(r_5)\times {\cal B}_{2m+1,b}(r_5)$.
\end{enumerate}
\end{props}

\subsection{From estimates on $NF^m$ to estimates on $a^{m-1}$}

As we have seen above, our proposition rests on the assumption that the quantity $\|D(a^{m})(u)\|_{D_p(b,t_{m})}$ is less than $1$. Nevertheless, the only quantity 
which can be easily estimated (in particular, through the induction process) is $|D_u(NF^m)|_{D_p(b,t_{m}),r}$ as well as 
$|NF^m|_{D_p(b,t_{m}),r}$. The next lemma translates the estimates $D_u(NF^m)$ and $NF^m$ into an estimate for $\|D(a^{m})\|_{D_p(b,t_{m})}$.

By definition, we have, for any integer $1\leq i\leq l$,
\begin{eqnarray*}
S_j & = & \sum_{k=1}^n{\lambda_{j,k}x_k\frac{\partial}{\partial x_k}},\\
NF^{m}(x,u) & = & \sum_{j=1}^l{a_{j}^{m}(u)S_j}:=\sum_{k=1}^n{x_k g_{k}(u)\frac{\partial}{\partial x_k}},
\end{eqnarray*}
with
$$
g_{k}(u)= \left(\sum_{j=1}^l{\lambda_{j,k}a_{j}^{m}(u)}\right).
$$
This can be rewritten under the following form
$$
\begin{pmatrix} g_{1}(u)\\ \vdots\\ \vdots\\ g_{n}(u) \end{pmatrix}
=
\begin{pmatrix} \lambda_{1,1} & \ldots & \lambda_{l,1}\\ \vdots & & \vdots \\ \vdots & & \vdots \\
\lambda_{1,n} & \ldots & \lambda_{l,n} \end{pmatrix}
\begin{pmatrix} a_{1}^{m}(u)\\ \vdots\\ a_{l}^{m}(u) \end{pmatrix}.
$$
We recall that $l$ is less than or equal to $n$.
Since the $S_i$'s are linearly independent over ${\bf C}$, the matrix 
$(\lambda_{j,i})_{\substack{1\leq i\leq n\\ 1\leq j\leq l}}$ has rank $l$. 
Without any loss of generality, 
we can assume that the matrix \\$L:= (\lambda_{j,i})_{1\leq i,j\leq l}$ is invertible with $L^{-1}:= (\tilde\lambda_{i,j})_{1\leq i,j\leq l}$ as inverse.

\begin{lemms}\label{a-nf}
Let $1/2<r\leq 1$ and let $\eta_1$ be a positive number. Let us set $\eta:=\frac{1}{2l|L^{-1}|}$.
If $|NF^{m}|_{D_p(b,t_{m}),r}$ is less than $\eta\eta_1$ and $\left|D_u(NF^m)\right|_{D_p(b,t_{m}),r}$ is less than $\eta$, then we have
\begin{eqnarray*}
\max_j\|D_u a^{m}_j(u)\|_{D_p(b,t_{m})} & < & 1\\
\max_j\|a_{j}^{m}(u) \|_{D_p(b,t_{m})}& < & \eta_1.
\end{eqnarray*}
Moreover, we have 
$$
\max_j\|a_j^{2m}(u)- a_j^{m}(u)\|_{D_p(b,t_{m})}\leq 2l|L^{-1}||NF^{2m}-NF^{m}|_{D_p(b,t_{m}),r}.
$$
\end{lemms}
\begin{proof}
 
We can write, for any integer $1\leq j\leq l$,
$$
a_{j}^{m}(u) =  \sum_{k=1}^{l}{\tilde \lambda_{j,k}g_{k}(u)}.
$$
Since $r$ is greater than $1/2$, we have
$$
\|g_{k}(u)\|_{D_p(b,t_{m})} \leq  2r \|g_{k}(u)\|_{D_p(b,t_{m})}
 \leq  2\left|x_k\tilde g_{k}(u)\right|_{D_p(b,t_{m}),r} \leq  2|NF^{m}|_{D_p(b,t_{m}),r}.
$$
As a consequence, we obtain 
\begin{equation}\label{estim-A}
\max_j\|a_{j}^{m}(u)\|_{D_p(b,t_{m})}\leq 2l|L^{-1}||NF^{m}|_{D_p(b,t_{m}),r}.
\end{equation}
On the other hand, for any integer $1\leq k\leq n$, we have
\begin{eqnarray*}
\frac{\partial NF^m}{\partial u_k} &= &\sum_{j=1}^l\frac{\partial a_{j}^{m}(u)}{\partial u_k}S_j =: \sum_{i=1}^nx_iG_i^k(u)\frac{\partial}{\partial x_i}.
\end{eqnarray*}
As above,  for any integers $1\leq j\leq l$ and $1\leq k\leq n$, we have
$$
\frac{\partial a_{j}^{m}(u)}{\partial u_k}= \sum_{i=1}^{l}{\tilde \lambda_{j,i}G_{i}^k(u)}.
$$
Therefore, we obtain the following estimate: 
\begin{eqnarray*}
\left\|\frac{\partial a_{j}^{m}(u)}{\partial u_k}\right\|_{D_p(b,t_{m})} &\leq & l|L^{-1}|\max_i\|G_{i}^k(u)\|_{D_p(b,t_{m})}\\ 
& \leq & 2l|L^{-1}|\max_i|x_iG_{i}^k(u)|_{D_p(b,t_{m}),r}\\ 
& \leq & 2l|L^{-1}|\left|\frac{\partial NF^m}{\partial u_k}\right|_{D_p(b,t_{m}),r}.\\
\end{eqnarray*}

Let us set $\eta:=\frac{1}{2l|L^{-1}|}$. Let us assume that 
$|NF^{m}|_{D_p(b,t_{m}),r}$ is less than $\eta\eta_1$ and that $\left|\partial NF^m/\partial u_k\right|_{D_p(b,t_{m}),r}$ is less the $\eta$. Then $\left|\partial a_{j}^{m}(u)/\partial u_k\right|_{D_p(b,t_{m}),r}$ is less than $1$ and $\max_j|a_{j}^{m}(u)|_{D_p(b,t_{m}),r}$ is less than $\eta_1$. The last statement of the lemma is proved by estimates $(\ref{estim-A})$ applied to 
$a_j^{2m-1}-a_j^{m}$ instead of $a_j^{m}$.
\qed\end{proof}

\subsection{The image of $D_n(0,r)\times D_p(b,t_{m})$ by $\tilde \Phi_m$ and its inclusions}

The aim of this section is to prove the following proposition:
\begin{props}\label{inclusions}
Let $1/2<r\leq 1$ and let $b$ belongs to ${\cal K}_k$. We assume that 
$$
\max_j\|D(a_j^{m})\|_{D_p(b,t_{m})}< 1
\text{ and }\max_j\|a_j^{2m}-a_j^{m}\|_{D_p(b,t_{m})}< 2l|L^{-1}|.
$$
If $m$ is sufficiently large then, for any positive numbers $\nu,\nu'$ and $\nu''$ such that $\nu<2\nu'$ and $\nu'+1/6<\nu''\leq 1$ and for any integer $q=0,1$, we have
$$
D_n(0,r_{q+2})\times D_p(b,t_{2m})\subset \tilde\Phi_m\left(D_n(0,r_{q+1})\times
D_p(b,t_{m})\right)
$$
\end{props}

\begin{lemms}\label{epsilon-vois}
Let $q$ be a nonnegative integer.
Under the assumptions of proposition \ref{inclusions} and if $m$ is sufficiently large (say $m\geq m_2(q)$) then, 
for all $0\leq \nu', \nu''\leq 1$ such that $84\nu''-36\nu'> 9$, 
the $\epsilon$-neighborhood of $D_n(0,r_{q+1})\times D_p(b,\nu't_{2m})$ is included in $D_n(0,r_{q})\times D_p(b,\nu''t_{m})$
with 
$$
\epsilon=\frac{\gamma\omega_{k+1}}{24l\Lambda(2m+1)}.
$$
\end{lemms}
\begin{proof}
First of all, let us show that $r_p-r_{p+1}$ is greater than $\epsilon$. In fact, since $r$ is less than $1/2$, we have 
\begin{eqnarray*}
\theta_k^{p}r-\theta_k^{p+1}r & = & \theta_k^{p}r(1-\theta_k)\\
&=& \left(\gamma_km^{-2/m}\right)^{p}r\left(1-\gamma_km^{-2/m}\right)\\
& =  & \left(\frac{\gamma^2\omega_{k+1}^2}{c_1}\right)^{p/m}m^{-2p/m}r\left(1-\left(\frac{\gamma^2\omega_{k+1}^2}{c_1}\right)^{1/m}m^{-2/m}\right)\\
& > & 1/2\left(\frac{\gamma^2\omega_{k+1}^2}{c_1}\right)^{p/m}m^{-2p/m}\left(1-\left(\frac{\gamma^2\omega_{k+1}^2}{c_1}\right)^{1/m}m^{-2/m}\right).
\end{eqnarray*}
We want to show that, if $m$ is sufficiently large,  
\begin{eqnarray}\label{maj-epsilon}
\frac{2d}{24pl\Lambda(2m+1)}& < & m^{-2p/m}\left(\frac{d^2}{c_1}\right)^{p/m}\left(1-\left(\frac{d^2}{c_1}\right)^{1/m}m^{-2/m}\right)
\end{eqnarray}
wher we have set $d:= \gamma\omega_{k+1}$. By assumption, the serie $$-\sum_{k\in {\bf N}^*}\frac{\ln\omega_{k+1}}{2^k}$$ converges; thus its general term tends to zero as $k$ tends to infinity. By applying the logarithm, we conclude that 
$$
\lim_{k\rightarrow +\infty}\left(\frac{\gamma^2\omega_{k+1}^2}{c_1}\right)^{1/m}=1.
$$
We recall that $m=2^k$. 
Therefore, there exists an integer $m_1=2^{k_1}$ such that, for any $m=2^k\geq m_1$, we have 
$$
\left(\frac{d^2}{c_1}\right)^{1/m}>1/2.
$$
On the other hand, we have 
$$
\left(\frac{d^2}{c_1}\right)^{1/m}m^{-2/m}< 1.
$$
Thus, in order to prove inequality $(\ref{maj-epsilon})$, it is sufficient to prove, that if $m\geq m_1$ is sufficiently large, then
$$
\frac{2d}{24pl\Lambda(2m+1)}<\frac{m^{-2p/m}}{2^{p}}\left(1-\left(\frac{d^2}{c_1}\right)^{1/m}m^{-2/m}\right).
$$
Let $f$ be the function of $x>0$ defined to be
$$
f(x)= \frac{2}{24pl\Lambda(2x+1)}+\frac{1}{2^{p}}x^{-2p/x}\left(x^{-2/x}-1\right).
$$
Since $d^2/c_1\leq 1$ and $d\leq 1$, we have, for all $m\in {\bf N}^*$,
$$
\frac{2d}{24pl\Lambda(2m+1)}+\frac{m^{-2p/m}}{2^{p}}\left(\left(\frac{d^2}{c_1}\right)^{1/m}m^{-2/m}-1\right)\leq f(m).
$$
Since both $x^{-2/x}$ and $x^{-2p/x}$ tends to $1$ as $x$ tends to plus infinity, we conclude that $f(x)$ tends to zero as $x$ tends to plus infinity. The computation of the derivative of $f$ shows that:
\begin{eqnarray*}
x^2f'(x) & = & \frac{1}{2^{p}}\left(-2(p+1)(1-\ln x)x^{-2(p+1)/x}+2p(1-\ln x)x^{-2p/x}\right)\\
& & -\frac{2d}{12pl\Lambda(2+1/x)^2}\\
& = & (1-\ln x)x^{-2p/x}(2p-2(p+1)x^{-2/x})-\frac{2d}{12pl\Lambda(2+1/x)^2}.
\end{eqnarray*}
Therefore when $x$ tends to plus infinity, so does $x^2f'(x)$, so $f$ is increasing from a certain point on and vanishes at infinity. Thus, there is a nonnegative integer $m_0$ from which $f$ is negative; $m_0$ depends only on $l$, $\Lambda$ and $p$. 

As a conclusion, if $m$ is greater than or equal to $m_2=\sup(m_0,m_1)$, then 
$$
\frac{2d}{24pl\Lambda(2m+1)}+\left(\frac{d^2}{c_1}\right)^{(p+1)/m}\frac{m
^{-2(p+1)/m}}{2^{p}}<\left(\frac{d^2}{c_1}\right)^{p/m}\frac{m^{-2p/m}}{2^{p}}.
$$ 
This means than $r_{p}-r_{p+1}$ is greater than $\epsilon$. Since $$\nu'/(4m+1)+ 1/(12(2m+1))<\nu''/((2m+1))$$ and $\omega_{k+2}\leq \omega_{k+1}$, we have
$$
\nu't_{2m} +\epsilon \leq \left(\frac{\nu'}{4m+1}+ \frac{1}{12(2m+1)}\right)\frac{\gamma\omega_{k+1}}{2l\Lambda}<\nu'' t_m.
$$
This concludes the proof.
\qed\end{proof}

\begin{lemms}\label{imageincluse}
Let $q$ be a nonnegative integer. Under the assumptions of proposition \ref{inclusions} and
if $m$ sufficiently large (say $m\geq m_3$) then, for any $0\leq \nu',
\nu''\leq 1$ such that $84 \nu''-36\nu'> 9$,  we have $$\tilde\Phi_m^{-1}(D_n(0,r_{q+2})\times D_p(b,\nu't_{2m}))\subset D_n(0,r_{q+1})\times D_p(b,\nu''t_{m}).$$
\end{lemms}
\begin{proof}
First of all, for any $(u,x)\in D_n(0,r_{p+2})\times D_p(b,\nu't_{2m})$, we have 
\begin{eqnarray}
|U_m(x,u)| & \leq & |U_m|_{D_p(b,\nu't_{2m}), r_{p+2}} \leq  |U_m|_{D_p(b,t_{m}), r_{p+2}}\nonumber\\
& \leq & \left(\theta_k^{p+2}\right)^{m+1}|U_m|_{D_p(b,t_{m}), r}\nonumber\\
& \leq & \theta_k^{p(m+1)+2}\left(\gamma_k^2 m^{-4/m}\right)^{m}|U_m|_{{\cal V}_{m,b}^1, r}\label{u1}\\
& \leq & \theta_k^{p+2}\left(\gamma_k^2 m^{-4/m}\right)^{m}|U_m|_{D_p(b,t_{m}), r}\nonumber\\
& \leq & \theta_k^{p+2}\left(\gamma_k^2 m^{-4/m}\right)^{m}\gamma_k^{-m}\nonumber\\
& \leq & \theta_k^{p+2}\gamma_k^{m}m^{-4}\nonumber\\
& \leq & \theta_k^{p+2}\frac{\gamma^2\omega_{k+1}^2}{c_1m^{4}}.\label{u2}
\end{eqnarray}
The third (resp. fifth, sixth) inequality is due to the fact that $U_m$ is of
order greater than or equal to $m+1$ (resp. $\theta_k^{p(m+1)}$ less than or equal to $\theta_k^{p}$, $|B^{m+1,2m}|_{D_p(b,t_{m}), r}$ less than $1$).
Let us set 
\begin{equation}
M  :=  \sup_{x\in D_{n}(0,1)} n|D\pi(x)|.\label{M}
\end{equation}
If $m$ is sufficiently large (say $m$ greater than or equal to $m_2'$), we have 
$$
\max(1,M)\left(\gamma_k^2 m^{-4/m}\right)\frac{\gamma^2\omega_{k+1}^2}{c_1m^{4}}<\frac{\gamma\omega_{k+1}}{24l\Lambda(2m+1)}=\epsilon,
$$
since $\omega_{k+1}$ is less than or equal to $1$. On the other hand, we have
$$
|\pi(x+U_m(u,x))-\pi(x)|\leq M |U_m(u,x)|
$$
since the point $x+U_m(u,x)$ belongs to $D_{n}(0,1)$.

Let us set $m_3=\max(m_2,m_2')$. We can conclude the proof by applying the previous lemma to 
$$
\tilde \Phi_m(u,x)=(x+U_m(u,x), u+\pi(x+U_m(u,x))-\pi(x)).
$$
\qed\end{proof}


\subsection{Estimates from the induction process}

Let us assume that $\tilde X_m$ is normalized up to order $m=2^k$ along $\Sigma$ and let us set $pr_1(\tilde X_m)=NF^m+B+C+r_{\Sigma}$ where $B$ belongs to ${\cal O}_p(D_p(b,\nu't_{m}))\otimes\vfo n {m+1}$ and $C$ belongs to ${\cal O}_p(D_p(b,\nu't_{m}))\otimes\vfo n {2m+1}$. We assume that $NF^m$ belongs to ${\cal NF}_{m}(r)$ and that $B+C$ belongs to ${\cal B}_{m+1}(r)$.

Let $\tilde \Phi_m$ be the normalizing diffeomorphism of the previous section. Let us set $y=pr_1\tilde\Phi_m(x,u)$ as well as 
$$
\tilde\Phi_m^*X(y,v)= NF^m(y,v)+B'(y,v)+C'(y,v)+r'_{\Sigma}(y,v)=:Z(y,v)+r'_{\Sigma}(y,v)
$$ 
where $B'$ belongs to ${\cal O}_p(D_p(b,\nu''t_{m}))\otimes\pvf n {m+1} {2m}$ and $C'$ belongs
to ${\cal O}_p(D_p(b,\nu''t_{m}))\otimes\vfo n {2m+1}$.
In order to give an estimate for $C'$, we shall use the following equality:
\begin{eqnarray*}
C'(y,v) & = &  \left(NF^m(\tilde \Phi_m^{-1}(y,v))- NF^m(v,y)\right)+(B+C)\left(\tilde \Phi_m^{-1}(y,v)\right)-B'(y,v)\\
& & -\left(D_y(U_m)(y,v)+D_v(U_m)(y,v)D\pi(y)\right)Z(y,v).
\end{eqnarray*}
Therefore, we obtain the following estimate:
\begin{eqnarray*}
\|C'\|_{D_p(b,t_{2m}),r_3}& \leq & \|NF^m\circ\tilde\Phi_m^{-1}- NF^m\|_{D_p(b,t_{2m}),r_3}+\|(B+C)\circ\tilde\Phi_m^{-1}\|_{D_p(b,t_{2m}),r_3}\\
& & +\|D_y(U_m)(NF^m+B'+C')\|_{D_p(b,t_{2m}),r_3}+|B'|_{D_p(b,t_{2m}),r_3}\\
& & + \|D_v(U_m)(y,v)D\pi(y)(B'+C')\|_{D_p(b,t_{2m}),r_3}.
\end{eqnarray*}
According to lemma $\ref{imageincluse}$, we have
$$
\tilde\Phi_m^{-1}(D_n(0,r_3)\times D_p(b,t_{2m}))\subset D_n(0,r_2)\times D_p(b,\nu't_{m}).
$$ 
Therefore we have 
\begin{eqnarray*}
\|(B+C)\circ\tilde\Phi_m^{-1}\|_{D_p(b,t_{2m}),r_3}& \leq & \|(B+C)\|_{D_p(b,\nu't_{m}),r_2}\leq |(B+C)|_{D_p(b,\nu't_{m}),r_2}\\
& \leq & \left(\frac{r_2}{r}\right)^{m+1}|B+C|_{D_p(b,t_{m}),r}\\
& \leq & \left(\frac{\gamma^2\omega_{k+1}^2}{c_1}\right)^{2+2/m}m^{-4-4/m}\\
& \leq & \frac{1}{m^4}.
\end{eqnarray*}
Let us set $z:=U_m(y,v)=x-y$, $w:=\pi(y+U_m(y,v))-\pi(y)=u-v$ and $M'=(n+p)\max(1,M)$.
For the same reason as above, we have the following estimates: 
\begin{eqnarray}
\|NF^m\circ\Phi_m^{-1}- NF^m\|_{D_p(b,t_{2m}),r_3} & \leq & 
\left\|\int_0^1D(NF^m)\left(y+tz,v+tw\right)(z,w)dt\right\|_{D_p(b,t_{2m}),r_3}\nonumber\\
& \leq & |D(NF^m)|_{D_p(b,\nu't_{m}),r_2},M'|U_m|_{D_p(b,t_{2m}),r_3}\nonumber\\
& \leq & M'|D(NF^m)|_{D_p(b,\nu't_{m}),r_2}\frac{\gamma^2\omega_{k+1}^2\gamma_k^3 m^{-6/m}}{c_1m^{4}}\label{maj-nf}\\
& \leq & M'|D(NF^m)|_{D_p(b,\nu't_{m}),r_2}\frac{1}{m^4}\nonumber\\
& \leq & M'(2l|L^{-1}|\max_j|S_j|_1+1/(2l|L^{-1}|))\frac{1}{m^4}.\nonumber
\end{eqnarray}
The second inequality comes from proposition $\ref{epsilon-vois}$ while the last one comes from the fact that, by assumption, 
$$
D(NF^m)(y,v)=\sum_{j=1}^l(a_j(v)D_y(S_j))+D_v(NF^m)(y,v).
$$

In order to get a good estimate for $C'$, we need to give an estimate of
$\|D(U_m)\|_{D_p(b,\nu' t_{m}),r_2}$. In the next few lines, we shall relate 
$|D(U_m)|_{D_p(b,\nu' t_{m}),r_2}$ to $|U_m|_{D_p(b,t_{m}),r}$.

Let $f$ belongs to ${\cal O}_p(D_p(b,t))\otimes {\cal M}_n^{m+1}(D_n(0,r_p))$. Let us
assume that $|f|_{D_p(b,t),r_p}$ is finite. It can be written $$f(y,v)=\sum_{m+1\leq |Q|}f_Q(v)y^Q$$ where the $f_Q$'s belong to ${\cal O}_p(D_p(b,t))$. For $1\leq i\leq n$ and $1\leq j\leq p$, we have
$$
\frac{\partial f}{\partial y_i}= \sum_{m+1\leq |Q|}q_if_Q(v)\frac{y^Q}{y_i},\quad\quad\frac{\partial f}{\partial v_j}= \sum_{m+1\leq |Q|}\frac{\partial f_Q}{\partial v_j}y^Q.
$$
According to Cauchy estimates, we have
$$
\|f_Q\|_{D_p(b,t)}\leq  \frac{|f|_{D_p(b,t),r_p}}{r_p^{|Q|}}.
$$
Hence, 
\begin{eqnarray*}
\left|\frac{\partial f}{\partial y_i}\right|_{D_p(b,t),r_{p+2}}& \leq & \frac{|f|_{D_p(b,t),r_{p}}}{r_{p+2}}\sum_{m+1\leq |Q|}q_i\left(\frac{r_{p+2}}{r_p}\right)^Q\\
& \leq & \frac{|f|_{D_p(b,t),r_{p}}}{r_{p+2}}\sum_{m+1\leq |Q|}q_i\theta_k^{2|Q|}.
\end{eqnarray*}
The number of monomials of order $k$ in $n$ variables is equal to
$\frac{(k+1)\cdots(k+n)}{n!}$ and it is less than or equal to $(2k)^n/n!$ if
$k\geq n$. In this case, we have
$$
\sum_{|Q|\geq m+1}q_i\theta_k^{2|Q|}\leq \sum_{q\geq m+1}\frac{(2q)^{n+1}\theta_k^{2q}}{n!}.
$$
If $m$ is large enough, $(2q)^{n+1}\theta_k^q/n!$ is less than or equal to $1$ when $q\geq m+1$. In fact, we have 
\begin{eqnarray*}
(2q)^{n+1}\theta_k^q/n! & \leq & \frac{2^{n+1}}{n!}\left(-\frac{n+1}{\ln \theta_k}\right)\theta_k^{-(n+1)/\ln \theta_k}\\
&\leq &\frac{2^{n+1}}{n!}\frac{m(m+1)}{2\ln m+\ln (c/\gamma^2\omega_{k+1}^2)}e^{-(m+1)}.
\end{eqnarray*}
Hence, we have
$$
\sum_{|Q|\geq m+1}\theta_k^{2|Q|}\leq \sum_{q\geq m+1}\theta_k^{q}=\frac{\theta_k^{m+1}}{1-\theta_k}.
$$
As a consequence, we have 
\begin{equation}
\left|\frac{\partial f}{\partial y_i}\right|_{D_p(b,t),r_{p+2}}\leq |f|_{D_p(b,t),r_p}\frac{\theta_k^{m+1}}{r_{p+2}(1-\theta_k)}.\label{deriveex}
\end{equation}
Let $0<t'<t$ be such that $t-t'>\epsilon$, we have
$$
\left|\frac{\partial f}{\partial v_j}\right|_{D_p(b,t'),r_p}  = 
\sum_{m+1\leq |Q|}\left\|\frac{\partial f_Q}{\partial v_i}\right\|_{D_p(b,t')}r_p^{|Q|}.
$$
By Cauchy integral formula, we obtain 
$$
\left\|\frac{\partial f_Q}{\partial v_i}\right\|_{D_p(b,t')}\leq \frac{\|f\|_{D_p(b,t)}}{\epsilon},
$$
hence
\begin{equation}
\left|\frac{\partial f}{\partial v_i}\right|_{D_p(b,t'),r_p}\leq  \frac{|f|_{D_p(b,t),r_{p}}}{\epsilon}\label{deriveeu}.
\end{equation}
We recall that 
$$
\epsilon=\frac{\gamma\omega_{k+1}}{24l\Lambda(2m+1)}.
$$
Since $B+C$ is of order greater than or equal to $m+1$, we have
\begin{eqnarray}
|B'|_{D_p(b,t_{2m}),r_3} & \leq & |B+C|_{D_p(b,t_{2m}),r_3}\nonumber\\
& \leq & \theta_k^{3(m+1)}|B+C|_{D_p(b,t_{m}),r}\nonumber\\
& \leq & \frac{\gamma^6\omega_{k+1}^6}{c_1^3m^6}|B+C|_{D_p(b,t_{m}),r}\label{maj-b}\\
& \leq & \frac{1}{m^6}.\nonumber
\end{eqnarray}
From this estimate and using the remark \ref{importante_remarque}, we obtain: 
\begin{eqnarray*}
|NF^{2m}|_{D_p(b,t_{2m}),r_3}& \leq & |NF^{m}|_{D_p(b,t_{2m}),r_3}+|B'_{|\Sigma}|_{D_p(b,t_{2m}),r_3}\\
& \leq & |NF^{m}|_{D_p(b,t_{m}),r}+|B'|_{D_p(b,t_{2m}),r_3}\\
& \leq & 1-\frac{1}{(m)^3}+ \frac{1}{m^6}\\
&\leq & 1-\frac{1}{(2m)^3}\quad\text{if}\quad m>\sqrt[3]{8/7}.
\end{eqnarray*}
Let us set $Z:=NF^m+B'+C'$. On the other hand, we have
\begin{eqnarray*}
\left\|D_yU_mZ\right\|_{D_p(b,t_{2m}),r_3} & \leq & n\|D_yU_m\|_{D_p(b,t_{2m}),r_3}\left(\|Z\|_{D_p(b,t_{2m}),r_3}\right)\\
& \leq & \frac{n\left(\theta_k\right)^{m+1}}{r_3(1-\theta_k)}|U_m|_{D_p(b,t_{m}),r_1}\\
& & \times \left(1-\frac{1}{(2m)^3}+\|C'\|_{D_p(b,t_{2m}),r_3}\right)\\
& \leq & \frac{2n\left(\theta_k\right)^{2m-1}}{(1-\theta_k)}|U_m|_{D_p(b,t_{m}),r}\\
& & \times \left(1-\frac{1}{(2m)^3}+\|C'\|_{D_p(b,t_{2m}),r_3}\right)\\
& \leq & \frac{2n\left(\theta_k\right)^{m-1}}{m^2(1-\theta_k)}\left(1-\frac{1}{(2m)^3}+\|C'\|_{D_p(b,t_{2m}),r_3}\right).
\end{eqnarray*}
The last inequality comes from the fact that $|U_m|_{D_p(b,t_{m}),r}$ is less than or equal to $\gamma_k^{-m}$ as well as $\theta_k^m=\gamma_k^{m}m^{-2}$. 
Let us show that, if $m\geq m_4''$ is sufficiently large, then we have  
\begin{equation}
\left|1-\frac{\gamma^2\omega_{k+1}^2}{c_1m^{2}}\right| >  \theta_k. \label{theta}
\end{equation}
Since $\frac{\gamma^2\omega_{k+1}^2}{c_1}\leq 1$, we have
$$
\left|1-\frac{\gamma^2\omega_{k+1}^2}{c_1m^{2}}\right|\geq \left|1-m^{-2}\right|
$$
as well as
$$
m^{-2/m}\geq \left(\frac{c_1}{\gamma^2\omega_{k+1}^2}\right)^{-1/m}m^{-2/m}.
$$
Therefore, it is sufficient to show that: 
$$
\left|1-m^{-2}\right| > m^{-2/m}.
$$
For that purpose, let $f$ be the function of the real variable $x$ which is assumed to be greater than or equal to $2$ and defined to be $$f(x)=1-x^{-2}- x^{-2/x}.$$ We have $x^{2}f'(x)=2x^{-1}+2(1-\ln x )x^{-2/x}$. 
Since $x^{2}f'(x)$ tends to minus infinity as $x$ tends to plus infinity, the function $f$ is
decreasing from a certain point on. But since $f(x)$ tends to zero as $x$ tends to plus infinity, $f$ is positive from a certain point $m_3''$ on. Therefore, if $m\geq m_4''\geq m_3''$, then 
\begin{eqnarray*}
\frac{2n\left(\theta_k\right)^{m-1}}{m^2(1-\theta_k)}& \leq & \frac{2nc_1}{\gamma^2\omega^2_{k+1}}m^{-2+2/m}\left(\frac{\gamma^2\omega_{k+1}^2}{c_1}\right)^{1-1/m}\\
& \leq & \frac{2n}{m^2}.
\end{eqnarray*}
As a consequence, if $m$ is large enough, we have
\begin{equation}
\left\|D_yU_mZ\right\|_{D_p(b,t_{2m}),r_3} \leq \frac{2n}{m^2}\left(1-\frac{1}{(2m)^3}+\|C'\|_{D_p(b,t_{2m}),r_3}\right).
\end{equation}
Using $(\ref{deriveeu})$, we have
\begin{eqnarray*}
\left\|D_v(U_m)D\pi(x)(B'+C')\right\|_{D_p(b,t_{2m}),r_3} & \leq & \frac{pM}{\epsilon}|U_m|_{D_p(b,t_{m}),r_3}\|B'+C'\|_{D_p(b,t_{2m}),r_3}\\
&\leq &
\frac{24pMl\Lambda(2m+1)\theta_k^{3(m+1)}\gamma_k^{-m}}{\gamma\omega_{k+1}}\\
&&\times\left(1/m^6+\|C'\|_{D_p(b,t_{2m}),r_3}\right)\\
&\leq & \frac{72pMl\Lambda}{m^5c_1^6}\gamma^5\omega_{k+1}^5\left(1/m^6+\|C'\|_{D_p(b,t_{2m}),r_3}\right).
\end{eqnarray*}
Since $B'$ is of order greater than or equal to $m+1$, we have 
\begin{eqnarray*}
|D_v(B'_{|\Sigma})|_{D_p(b,t_{2m}),r_3} & \leq & \epsilon^{-1}|B'_{|\Sigma}|_{D_p(b,t_{m}),r_3}\\
& \leq & \epsilon^{-1}\theta_k^{3(m+1)}|B'|_{D_p(b,t_{m}),r}\\
& \leq & \frac{72l\Lambda\gamma^5\omega^5_{k+1}}{c_1^6m^5}.
\end{eqnarray*}
Therefore, we have
\begin{eqnarray*}
|D_v(NF^{2m})|_{D_p(b,t_{2m}),r_3} & \leq & |D_v(NF^m)|_{D_p(b,t_{2m}),r_3}+|D_v(B'_{|\Sigma})|_{D_p(b,t_{2m}),r_3}\\
& \leq & 1/(2l|L^{-1}|)-1/m^2+\frac{72l\Lambda\gamma^5\omega^5_{k+1}}{c_1^6m^5}\\
& < & 1/(2l|L^{-1}|)-1/(2m)^2\quad\text{if $m^3>\frac{144l\Lambda\gamma^5}{c_1^6}$}.
\end{eqnarray*}
We just have shown that $NF^{2m}\in {\cal NF}_{2m}(r_3)$.

At last, we have
\begin{eqnarray*}
\|C'\|_{D_p(b,t_{2m}),r_3}& \leq & \left(1+M'(2l|L^{-1}|\max_j|S_j|_1+1/(2l|L^{-1}|))\right)\frac{1}{m^4}\\
& & +\frac{2n}{m^2}\left(1-\frac{1}{(2m)^3}+\|C'\|_{D_p(b,t_{2m}),r_3}\right)\\
& & +\frac{72pMl\Lambda}{m^5c_1^6}\gamma^5\omega_{k+1}^5\left(1/m^6+\|C'\|_{D_p(b,t_{2m}),r_3}\right).
\end{eqnarray*}
If $m$ is large enough, we have
$$
\frac{2n}{m^2}+\frac{72pMl\Lambda}{m^5c_1^6}\gamma^5\omega_{k+1}^5\leq 1/2,
$$
as well as 
$$
\frac{A}{m^4}-\frac{2n}{8m^5} + \frac{72pMl\Lambda}{m^{11}c_1^6}\gamma^5\omega_{k+1}^5\leq \frac{2n}{m^2},
$$
where $A:=\left(1+M'(2l|L^{-1}|\max_j|S_j|_1+1/(2l|L^{-1}|))\right)$.
Hence, we have 
$$
\|C'\|_{D_p(b,t_{2m}),r_3}<\frac{2n}{m^2}.
$$
Using Cauchy estimates with $C'$, which is of order greater than or equal to $2m+1$, we have
\begin{eqnarray*}
|C'|_{D_p(b,t_{2m}),r_5}& \leq & \sum_{
|Q|\geq 2m+1}\|c_Q\|_{D_p(b,t_{2m}),r_5}r_5^Q\\
&\leq & \|C'\|_{D_p(b,t_{2m}),r_3}\sum_{|Q|\geq
 2m+1}\left(\frac{r_5}{r_3}\right)^Q\\
&\leq & \|C'\|_{D_p(b,t_{2m}),r_3}\sum_{|Q|\geq
 2m+1}\theta_k^{2|Q|}\\
&\leq & \frac{\|C'\|_{D_p(b,t_{2m}),r_3}\theta_k^{2m+1}}{1-\theta_k}\\
& \leq & \|C'\|_{D_p(b,t_{2m}),r_3}\frac{c_1m^2}{\gamma^2\omega_{k+1}^2}\theta_k^{2m}\\
& \leq & \|C'\|_{D_p(b,t_{2m}),r_3}\frac{1}{m^2}\leq \frac{2n}{m^4}.
\end{eqnarray*}
Therefore, we have shown that $C'\in {\cal B}_{2m+1}(r_5)$ provided that $m$ is large enough.


\section{Proof of the existence of an invariant analytic set}

Let $1/2<r \leq 1$ be a positive number and let $\{R_k\}_{k\geq 0}$ be the sequence of positive real numbers defined by induction as follows:
\begin{eqnarray*}
R_0 & = & r\\
R_{k+1} & = & \gamma_k^5m^{-10/m}R_k\quad\text{where}\quad m=2^k.
\end{eqnarray*}
\begin{lemms}\label{rayon-conv}
The sequence $\{R_k\}_{k\geq 0}$ converges to a positive number and there exists an integer $k_1$ such that, for all $k>k_1$, $R_k>R_{k_1}/2$.
\end{lemms}
\begin{proof}
We recall that 
$$
\gamma_k=\left(\frac{c_1(\eta_1)}{\gamma^2\omega_{k+1}^2}\right)^{-1/2^k}.
$$
Since we have  
$$
R_{k+1}=r\prod_{i=1}^k\gamma_i^5(2^i)^{-2^{-i}10}
$$
then, by applying the logarithm, we obtain
$$
\ln R_{k+1} = \ln r +10\sum_{i=1}^k{\frac{\ln \omega_{i+1}}{2^i}}
-5\ln (c_1/\gamma^2)\sum_{i=1}^k{\frac{1}{2^i}}-10\ln 2\sum_{i=1}^k{\frac{i}{2^i}}.
$$
The last two sums of the right hand side are convergent series, so is the first one since $\omega$ is a diophantine sequence. 
Therefore, there exists an integer $k_1$ such that 
$$
\prod_{i=k_1+1}^{+\infty}\gamma_i^5(2^i)^{-2^{-i}10}>1/2.
$$
Thus, if $k>k_1$ then, we have 
$$
R_k=R_{k_1}\prod_{i=k_1+1}^k\gamma_i^5(2^i)^{-2^{-i}10}>\frac{R_{m_1}}{2}.
$$
\qed\end{proof}

\subsection{The sequence of inverse diffeomorphisms converges to a holomorphic map $\Theta_{b}$ on $D_n(0,1/2)\times \{b\}$}

Let us assume that $X=X_0+R_{M_0+1}$ is a perturbation of order $M_0+1=2^K_0+1$ of $X_0$. Let $\tilde X= NF^{M_0}(x,u)+ \tilde R_{M_0+1}(x,u)$ be fibered along $\Sigma$ over $X$ such that $\tilde R_{M_0+1|\Sigma}=R_{M_0+1}$ and $NF^{M_0}(x,u)=\sum_{j=1}^la_j(u)S_j(x)=X_{0|\Sigma}$. 

{\bf We assume that $\tilde X$ is a good perturbation of $NF^{M_0}$ and that $ {\cal K}_{\infty}$ is not empty}.
Let $b\in {\cal K}_{\infty}$. {\bf Let us assume that}
$$
(NF^{M_0}, R_{M_0+1})\in {\cal NF}_{M_0,b}(1)\times {\cal B}_{m+1,b}(1).
$$
As above, we may define the sequence of positive real numbers $\{R_k\}_{k\geq K_0}$, with $R_{K_0}=1$. Thus, for any integer $k$ greater than $K_0$, we have $1/2<R_k\leq 1$. 
Let us prove by induction on $k\geq K_0$, that there exists a fibered diffeomorphism $\tilde\Psi_k$ of $({\bf C}^{n+p},(0,b))$ such that the vector field $$\tilde\Psi_k^*(NF^{M_0}+\tilde R_{M_0+1}):= NF^{2^{k+1}}+\tilde R_{2^{k+1}+1}\mod \Sigma$$ is normalized up to order $2^{k+1}$ along $\Sigma$, $(NF^{2^{k+1}},R_{2^{k+1}+1})$ belongs to the space ${\cal NF}_{2^{k+1},b}(R_{k+1})\times {\cal B}_{2^{k+1}+1,b}(R_{k+1})$ and 
$$
\|\text{Id} -\tilde\Psi_k^{-1}\|_{D_p(b,t_{2M_0}),R_{K_0+1}}\leq\max(1,M)\sum_{p=K_0}^k\frac{1}{2^{4p}}.
$$
Here we have set $\|.\|_{D_p(b,t_{2M_0}),R_{K_0+1}}:=\|.\|_{D_p(b,t_{2M_0})\times D_n(0,R_{K_0+1})}$.
\begin{itemize}
\item For $k=K_0$: according to proposition $(\ref{reccurence})$, there exists a diffeomorphism $\tilde\Phi_{M_0}$ such that 
$\tilde\Phi_{M_0}^*(NF^{M_0}+\tilde R_{M_0+1})=NF^{2M_0}+\tilde R_{2M_0+1}\mod \Sigma$ is normalized up to order $2M_0$ along $\Sigma$ and
$(NF^{2M_0},R_{2M_0+1})$ belongs to ${\cal NF}_{2M_0,b}(R_{K+1})\times {\cal B}_{2M_0+1,b}(R_{K_0+1})$. Moreover, we have 
\begin{eqnarray*} 
\left\|\text{Id}-\tilde\Phi_{M_0}^{-1}\right\|_{D_p(b,t_{2M_0}),R_{K_0+1}}& \leq  &  \left\|(U_{M_0},\pi(Id+U_{M_0})-\pi)\right\|_{D_p(b,t_{2M_0}),R_{K_0+1}}\\
&\leq & \max(1,M)\left\|U_{M_0}\right\|_{D_p(b,t_{2M_0}),R_{K_0+1}}\\
&\leq & \max(1,M)/m^4.
\end{eqnarray*}
The last inequality is due to estimate $(\ref{u2})$.

\item Let us assume that the result holds for all integers les than or equal to $k-1$: by assumptions, 
the vector field $$(\tilde\Psi_{k-1})_*(NF^{K}+\tilde R_{K+1}\mod \Sigma)= NF^{2^{k}}+\tilde R_{2^{k}+1}\mod \Sigma$$ is normalized up to order $2^{k}$ along $\Sigma$. Moreover, 
$(NF^{2^{k}},\tilde R_{2^{k}+1})$ belongs to ${\cal NF}_{2^{k},b}(R_{k})\times {\cal B}_{2^{k}+1,b}(R_{k})$.
According to lemma \ref{rayon-conv}, we have $1/2<R_k\leq 1$ and $b$ belongs to ${\cal K}_{k}$. Thus, we may apply proposition $(\ref{reccurence})$: there exists a diffeomorphism $\Phi_{2^k}$ such that the vector field
$$
\widetilde{(\Phi_{2^k}\circ\Psi_{k-1})}_*(NF^{M}+\tilde R_{M+1}\mod \Sigma)=NF^{2^{k+1}}+\tilde R_{2^{k+1}+1}\mod \Sigma
$$
is normalized up to order $2^{k+1}$ along $\Sigma$. Moreover, $(NF^{2^{k+1}},\tilde R_{2^{k+1}+1})$ belongs to ${\cal NF}_{2^{k+1},b}(R_{k+1})\times {\cal B}_{2^{k+1}+1,b}(R_{k+1})$.
Let us set $\tilde \Psi_{k}=\tilde \Phi_{2^k}\circ \tilde\Psi_{k-1}$. According to the first point of proposition $(\ref{reccurence})$ and estimate $(\ref{u2})$, 
we have 
$$\|\text{Id} - \tilde\Phi_{k}^{-1}\|_{D_p(b,t_{2^{k+1}})\times D_n(0,R_{k+1})}<\max(1,M)/2^{4{k}}.
$$
It follows that 
\begin{eqnarray*}
\|\text{Id} -\tilde\Psi_{k}^{-1}\|_{D_p(b,t_{2^{k+1}}),R_{k+1}}& \leq & \left\|(\text{Id}-\tilde\Psi_{k-1}^{-1})\circ \tilde\Phi_{2^k}^{-1}+(\text{Id} -\tilde\Phi_{2^k}^{-1})\right\|_{D_p(b,t_{2^{k+1}}),R_{k+1}}\\
& \leq & \left\|(\text{Id}-\tilde\Psi_{k-1}^{-1})\circ \tilde\Phi_{2^k}^{-1}\right\|_{D_p(b,t_{2^{k+1}}),R_{k+1}}\\
& & +\left\|\text{Id} -\tilde\Phi_{2^k}^{-1}\right\|_{D_p(b,t_{2^{k+1}}),R_{k+1}}\\
&\leq & \left\|(\text{Id}-\tilde\Psi_{k-1}^{-1})\right\|_{\tilde\Phi_{2^k}^{-1}(D_n(0,R_{k+1})\times D_p(b,t_{2^{k+1}}))}\\
&&+\left\|\text{Id} -\tilde\Phi_{2^k}^{-1}\right\|_{D_p(b,t_{2^{k+1}}),R_{k+1}}\\
&\leq &\max(1,M)\sum_{p=K}^{k-1}\frac{1}{2^{4p}}+\frac{\max(1,M)}{2^{4k}}\\
&\leq & \max(1,M)\sum_{p=1}^{+\infty}\frac{1}{2^{4p}}.\\
\end{eqnarray*} 
This ends the proof of the inductive step.
\end{itemize}
According to lemma \ref{rayon-conv}, if $R_K=1$, $D_n(0,1/2)$ is contained in $D_n(0,R_{k})$ for all integers $k$ greater than or equal to $K$.

Let us choose a positive number $\rho<1/2$ such that, if $m$ is large enough, 
$$
t_m+\rho^{|R_i|}\leq \left(\frac{1}{2}\right)^{|R_i|},\quad\quad i=1,\ldots,p.
$$
It follows that $D_n(0,\rho)\times\{b\}$ is contained in $D_n(0,R_{k+1})\times D_p(b,t_{2^{k+1}})$.
According to the previous estimate, the sequence 
$$
\left\{\tilde\Psi^{-1}_{k|D_n(0,\rho)\times\{b\}}\right\}_{k\geq K_0}
$$
of holomorphic functions is a uniformly bounded sequence (for the sup-norm) of
holomorphic maps on $D_n(0,\rho)\times\{b\}$. Therefore, according to Montel theorem, we can extract a subsequence which converges to a {\bf
  holomorphic map} $\Theta_b$ on $D_n(0,\rho)\times\{b\}$. The set
$\Theta_b(\pi^{-1}(b)\cap D_n(0,\rho)\times\{b\})$ is {\bf included in the
  manifold $\Sigma$}. In fact, by construction, $\tilde\Phi_m$ leaves $\Sigma$ invariant (globally), so does $\tilde\Psi_k$ as well as $\tilde\Psi_k^{-1}$. 
%

On the other hand, we have
$$ 
J^{2^{k+1}}\left(\tilde\Psi_{k+1}^{-1}(x,u)-\left(\begin{array}{c}x\\u\end{array}\right)\right)=\sum_{\substack{Q\in {\bf N}^n\\K_0\leq |Q|\leq 2^{k+1}}}\psi_Q(u)x^Q.
$$ 
where the $\psi_Q$'s are vector valued functions which don't depend on $k$. This follows from the definition by induction of $\tilde\Psi_{k+1}^{-1}$ and the fact that $U_{2^{k+1}}(x,u)$ is of order greater than or equal to $2^{k+1}+1$ (in $x$).
Anyway, we can write
$$
\tilde\Psi_{k+1}^{-1}(x,u)=\sum_{\substack{Q\in {\bf N}^n\\K_0\leq |Q|}}\psi_{Q,k}(u)x^Q.
$$
Let $x$ be a point in $D_n(0,\rho)$. Let us apply Cauchy inequalities to the $\psi_Q$'s:
$$
|\psi_{Q,k}(b)|\leq \frac{\left\|\tilde\Psi^{-1}_{k+1}-\text{Id }\right\|_{D_p(b,t_{2^{k+2}}),R_{k+2}}}{R_{k+2}^{|Q|}}.
$$
Therefore, we have the estimate:
\begin{eqnarray*}
\sum_{\substack{Q\in {\bf N}^n\\2^{k+1}+1\leq |Q|}}|\Psi_{Q,k}(b)x^Q| &\leq & \sum_{\substack{Q\in {\bf N}^n\\2^{k+1}+1\leq |Q|}}\frac{\left\|\tilde\Psi_{k+1}^{-1}-\text{Id }\right\|_{D_p(b,t_{2^{k+2}}),R_{k+2}}}{R_{k+2}^{|Q|}}|x|^{|Q|}.
\end{eqnarray*}
According to lemma \ref{rayon-conv}, we have $R_{k+2}>1/2$ and $|x|<1/2$. Thus, the series 
$$
\sum_{\substack{Q\in {\bf N}^n\\K_0\leq |Q|\leq 2^{k+1}}}\left(\frac{|x|}{R_{k+2}}\right)^{|Q|}
$$
converges.

Therefore, the series $J^{2^{k+1}}(\tilde\Psi_{k+1}^{-1}(x,b)-(x,b))$ converges on $D_n(0,\rho)$ to \\$\Theta_b(x)-(x,b)$. Hence, 
$\Theta_b$ is a biholomorphism in a neighborhood of the origin in ${\bf C}^n$. We shall denote by $\Psi_b$ its inverse. It follows that 
$$
{\cal V}_b:=\Theta_b(\pi^{-1}(b)\cap D_n(0,\rho))
$$
is a analytic subset of a neighborhood of the origin.


\subsection{The biholomorphism $\Psi_{b}$ conjugates the restriction of $X$ on 
${\cal V}_{b}$ to the restriction of a linear diagonal vector field on $\pi^{-1}(b)\cap D_n(0,\rho)$}

The sequence of linear vector fields
$$
\left\{NF^{2^k}_{|D_n(0,1/2)\times\{b\}}\right\}_{k\geq K}
$$
is uniformly bounded on $D_n(0,1/2)\times\{b\}$. This is due to the second
point of proposition \ref{reccurence} and to the fact that 
$$
\left|NF^{2^k}_{|D_n(0,1/2)\times\{b\}}\right|\leq \left|NF^{2^k}\right|_{D_p(b,t_{2^k}),R_k}.
$$
Therefore, there is a subsequence which converges to a linear vector field
$$
NF_b(x)=\sum_{i=}^n a_i(b)S_i(x).
$$
 
Let us show that the restriction of $\Psi_{b}$ to $V_{b}=\Theta_{b}(\pi^{-1}(b)\cap D_0(\rho)\times\{b\})$ conjugates the restriction of $\tilde X$ to $V_{b}$ to the restriction of $NF_b(y)$ to the set $\pi^{-1}(b)\cap D_0(\rho)\times\{b\}$. 

For any point $y$ in $\pi^{-1}(b)\cap D_n(0,\rho)$, let us set $\tilde
\Psi_k^{-1}(y,b)=(x,u)$. By construction, the point $(x,u)$ still belongs to
$\Sigma$. Moreover, we have
$$
(\tilde \Psi_k^{-1})_*(NF^{2^k}+\tilde R_{2^k+1}+\tilde r_{\Sigma}^k)(x,u)=\tilde X(x,u).
$$
Hence, we have 
$$
\tilde X(x,u)= D\tilde \Psi_k^{-1}(y,b)NF^{2^k}(y,b)+D\tilde \Psi_k^{-1}(y,b)\tilde R_{2^k+1}(y,b).
$$
Since $\{\tilde \Psi_k^{-1}(y,b)\}$ converges to $\Theta_{b}(y,b)$,
$\{\tilde NF^{2^k}(y,b)\}$ converges to $NF_b(y)$ and $\{\tilde R_{2^k+1}\}$
converges to $0$, the right hand side converges to $D\Theta_{b}(y,b)NF_b(y)$. Thus, we have, for any point $y$ in $\pi^{-1}(b)\cap D_0(\rho)$, 
$$
(\Theta_{b})_*NF_b(\Theta_{b}(y,b))=D\Theta_{b}(y,b)NF_b(y,b)= \tilde X(\Theta_{b}(y,b)).
$$ 

%
%


\section{Diophantine approximations on complex manifold}

The aim of this section is to prove theorem \ref{theo3}. It gives a sufficient condition that will ensure that the compact set ${\cal K}_{\infty}$ is not empty.
One way to achieve this, is to prove that it will have a positive measure. This section is an adaptation of the fourth part of R\"ussmann work \cite{russmann-weak}.

\subsection{The measure of diophantine points of the image of a small perturbation of nondegenerate map is positive}

The main goal of this section is to give an upper bound for the measure of 
the set of points whose values by a perturbation of a nondegenerate map are ``small". 

First of all, let us set some notation.
Let ${\cal B}$ be an open set in ${\bf R}^n$, let $p$ be a nonnegative integer, and let $f\in C^{p}({\cal B},{\bf R})$ be a $p$-times continuously differentiable function on ${\cal B}$. If $x=(x_1,\ldots, x_n)\in {\bf C}^n$, we shall set 
$$
|x|_2:=\sqrt{|x_1|^2+\cdots + |x_n|^2}.
$$
Let $a$ be a point in ${\bf R}^n$, let $k$ be a nonnegative integer less than or equal to $p$ and let $y$ belongs to ${\cal B}$. We shall set 
$$
D^kf(y)(a^k):= D^kf(y)(a,\ldots, a),\quad\|D^kf(y)\|:=\sup_{a\in {\bf R}^n, |a|_2=1}|D^pf(y)(a,\ldots, a)|
$$
as well as 
$$
\|f(y)\|_{p}:=\max_{0\leq k\leq p}\|D^kf(y)\|, \quad\|f\|_{{\cal B},p}=\max_{y\in {\cal B}}\|D^kf(y)\|.
$$
It is known that $|D^kf(y)(a_1,\ldots, a_k)|\leq \frac{k^k}{k!}\|D^kf(y)\|$ (see \cite{russmann-weak}). First of all, let us recall one of the results of R\"ussmann.
\begin{theos}[theorem 17.1\cite{russmann-weak}]\label{mes-epsilon}
Let ${\cal K}$ be a compact set ${\bf R}^n$ with diameter $d=\sup_{x,y\in {\cal K}}|x-y|_2$. Let $\vartheta$ be a positive number and let ${\cal B}$ be the $\vartheta$-neighborhood of ${\cal K}$ in ${\bf R}^n$. 
Let $g\in C^{\mu_0+1}({\cal B},{\bf R})$ be a function such that 
$$
\min_{y\in {\cal K}}\max_{0\leq \nu\leq \mu_0}\|D^{\nu}g(y)\|\geq \beta
$$
for some $\mu_0\in {\bf N}$ and for some positive number $\beta$. Then, for any function $\tilde g\in C^{\mu_0}({\cal B},{\bf R})$ satisfying $\|g-\tilde g\|_{{\cal B},\mu_0}\leq \beta/2$, we have the estimate
$$
\text{mes}_{n}\{y\in {\cal K}\,|\,|\tilde g(y)|\leq \epsilon\}\leq Bd^{n-1}\left(\frac{1}{\sqrt{n}}+2d+\frac{d}{\vartheta}\right)\left(\frac{\epsilon}{\beta}\right)^{\frac{1}{\mu_0}}\beta^{-1}\|g\|_{{\cal B},\mu_0+1},
$$
for all $0<\epsilon\leq \frac{\beta}{2\mu_0+2}$ and with $B= 3(2\pi e)^{n/2}(\mu_0+1)^{\mu_0+2}[(\mu_0+1)!]^{-1}$.
\end{theos}

\begin{lemms}\label{mubeta}
Let ${\cal U}$ be a connected open set on ${\bf C}^n$. Let $f:{\cal U}\rightarrow {\bf C}^l$ be a nondegenerate holomorphic map. Then, for any nonvoid compact set ${\cal K}\subset {\cal U}$, there is a positive integer $\mu_0$ and a positive number $\beta$ such that $$\left\||(c,f(y))|^2\right\|_{\mu_0}\geq \beta$$ for all $y$ in the compact set ${\cal K}$ and all $c$ in the unit sphere $$\Bbb S^l:=\{c\in {\bf C}^l\,|\,|c|_2=1\}.$$
\end{lemms}
\begin{proof}
We have to show that there exists $\mu_0$ and $\beta$ such that 
$$\max_{0\leq k\leq \mu_0}\left\|D^k|(c,f)|^2(y)\right\|\geq \beta$$ for all $y$ in the compact set ${\cal K}$ and all $c$ in the unit sphere $$\Bbb S^l:=\{c\in {\bf C}^l\,|\,|c|_2=1\}.$$ Let us assume that such $\mu_0$ and $\beta$ don't exist. Then, for any positive integer $\nu$, there would be $c_{\nu}\in \Bbb S^l$ and $y_{\nu}\in {\cal K}$ such that 
$$
\max_{0\leq k\leq \nu}\left\|D^k|(c_{\nu},f)|^2(y_{\nu})\right\|<\frac{1}{\nu};
$$
that is 
$$
\forall \nu\geq k+1, \forall a\in \Bbb S^n, \quad\left|D^k\left(|(c_{\nu},f)|^2\right)(y_{\nu})(a^k)\right|\leq \left\|D^k|(c_{\nu},f)|^2(y_{\nu})\right\|<\frac{1}{\nu}.
$$
Moreover, we have
\begin{eqnarray*}
D^k\left(|(c_{\nu},f)|^2\right)(y_{\nu})(a^k)& = &\sum_{i=0}^kC^k_iD^i\left((c_{\nu},f)\right)(y_{\nu})(a^i)D^{k-i}\left((\bar c_{\nu},\bar f)\right)(y_{\nu})(a^{k-i})\\
&= & \sum_{i=0}^kC^k_i\left((c_{\nu},D^if(y_{\nu})(a^i))\right)\left((\bar c_{\nu},D^{k-i}\bar f(y_{\nu})(a^{k-i}))\right).
\end{eqnarray*}
By compactness, we can extract a subsequence of $\{(c_{\nu},y_{\nu})\}_{\nu\geq 1}$ which converges to $(c,y)\in \Bbb S^l\times {\cal K}$. Thus, for any $a$ belonging to $\Bbb S^n$, 
the associated subsequence $\{D^k\left(|(c_{\nu},f)|^2\right)(y_{\nu})(a^k)\}$ converges to $D^k\left(|(c,f)|^2\right)(y)(a^k)$. According to the previous estimates, $D^k\left(|(c,f)|^2\right)(y)(a^k)$ vanishes for all $a\in {\bf C}^n$. Therefore, for any nonnegative integer $k$, $D^k\left(|(c,f)|^2\right)(y)=0$; that is, 
the real analytic function $|(c,f)|^2$ vanishes at $y\in {\cal K}$ as well 
as all its derivatives. Since ${\cal U}$ is connected, $|(c,f)|^2$ vanishes identically on ${\cal U}$; that is $(c,f)\equiv 0$. This contradicts the nondegeneracy of $f$.
\qed\end{proof}

Let ${\cal U}$ be a connected open set on ${\bf C}^n$, $f:{\cal U}\rightarrow {\bf C}^l$ a nondegenerate holomorphic map and ${\cal K}\subset {\cal U}$ nonvoid compact set. Let us set 
$$
\beta(\mu,f,{\cal K})=\min_{\substack{y\in {\cal K},\\ c\in \Bbb S^l}}\max_{0\leq k\leq \mu}\left|D^k|(c,f)|^2(y)\right|.
$$
Clearly, $\{\beta(\mu,f,{\cal K})\}_{\mu\geq 0}$ is a nondecreasing sequence
of nonnegative numbers. According to the previous lemma, there exists a positive integer $N_0$ such that $\beta(N_0,f,{\cal K})$ is positive. The smallest of these integers $N_0$ will be called the {\bf index of nondegeneracy} (of $f$ with respect to ${\cal K}$) and will be denoted by $\mu_0=\mu_0(f,{\cal K})$. The positive number 
$\beta(f,{\cal K}):= \beta(\mu_0(f,{\cal K}),f,{\cal K})$ will be called the {\bf amount of nondegeneracy}.
\begin{defis}
Let $\omega=\{\omega_k\}_{k\geq 1}$ be a diophantine sequence. We shall say that 
$S$ is {\bf strictly diophantine with respect to $(\omega,\mu_0)$} if 
$$
\lim_{k\rightarrow +\infty}\left(2^k+n+1\right)^{(n+1)}\left(\frac{\omega_k}{\omega_k(S)}\right)^{2/\mu_0}=0.
$$
In this case, we shall set 
$$
M_{\omega,\omega(S),2/\mu_0}:=\sup_{k\geq 1}\left(2^k+n+1\right)^{(n+1)}\left(\frac{\omega_k}{\omega_k(S)}\right)^{2/\mu_0} <+\infty,
$$
as well as $$M_{\omega,\omega(S)}:= \sup_{k\geq 1}\frac{\omega_k}{\omega_k(S)}.$$

\end{defis}
The next result shows the following: if $\tilde f$ is a small perturbation of a
nondegenerate holomorphic map $f$ from ${\bf C}^n$ to ${\bf C}^l$, then the
set of diophantine points (with respect to $S$) on the image of $\tilde f$ is
big in the sense that it has a positive measure. This kind of problems (for
real maps) is now classical and often called ``diophantine approximation on
manifolds" when the image of $\tilde f$ is (at least locally) a manifold. We
refer to \cite{bernik-book} for an up-to-date treatment of this topic (see
also \cite{margulis-dioph} for such a result in the real case).

We recall that $S: {\frak g}\rightarrow \hvf n 1$ is a Lie morphism from the
commutative Lie algebra ${\frak g}$. Let $\{g_1,\ldots, g_l\}$ be a basis of
${\frak g}$. We recall that ${\cal W}_{n,*}^{k,m}$ denotes the set of nonzero weights of $S$ into $\pvf n k m$.
\begin{props}\label{prop-mesure}
Let ${\cal U}$ be a connected open neighborhood of $0$ in ${\bf C}^n$, let $f:{\cal U}\rightarrow {\bf C}^l$ be a nondegenerate holomorphic map and let ${\cal K}\subset {\cal U}$ nonvoid compact set. Let $\mu_0$ be the index on nondegeneracy of $f$ with respect to ${\cal K}$. Assume that $S$ is strictly diophantine relatively to $(\omega=\{\omega_i\}_{i\geq 1},\mu_0)$. For any $\vartheta\in ]0, dist({\cal K}, {\bf C}^n\setminus {\cal U})[$, let ${\cal K}_{\vartheta}\subset {\cal U}$ denote the $\vartheta$-neighborhood of ${\cal K}$. Then, for any map $\tilde f\in C^{\mu_0}({\cal K}_{\vartheta}, {\bf C}^l)$ such that 
$$
\|f-\tilde f\|_{{\cal K}_{\vartheta}, \mu_0}\leq \frac{-\|f\|_{{\cal K}_{\vartheta},\mu_0}+\sqrt{\|f\|_{{\cal K}_{\vartheta},\mu_0}^2+\frac{\beta}{2^{\mu_0-1}}}}{2},
$$ 
the measure of the set 
$$
{\cal H}(\tilde f):=\left\{b\in {\cal K}\,|\,\forall i\in {\bf N}^*, \quad\forall \alpha\in {\cal W}_{n,*}^{2^{i-1}+1,2^{i}},\quad\left|\alpha\left(\sum_{j=1}^l{\widetilde 
{f_j}(b)g_j}\right)\right|\geq \gamma \omega_i\right\}
$$
satisfies $\text{mes}_{2n}{\cal H}(\tilde f)\geq mes_{2n}{\cal K}-\epsilon^*$ as soon as 
$0<\epsilon^*<mes_{2n}{\cal K}$ and $0<\gamma  \leq \gamma^*$ with
\begin{eqnarray*}
\gamma^*&=&\min\left[\left(\frac{\epsilon^*(n-1)!}{M\left((2^2+n)^na_2-(n+1)a_1+\frac{n}{4}M_{\omega,\omega(S),\mu_0}\right)}\right)^{\mu_0/2},
\frac{1}{M_{\omega,\omega(S) }}\sqrt{\frac{\beta}{2\mu_0+2}} \right]
\end{eqnarray*}
where
$M=Bd^{2n-1}(\frac{1}{\sqrt{2n}}+2d+\frac{d}{\vartheta})\beta^{-\frac{1}{\mu_0}-1}\sup_{c\in
  \Bbb S^l}\||(f,c)|^2\|_{{\cal K}_{\vartheta}, \mu_0+1}$, \\ $B=3(2\pi
e)^{n}(\mu_0+1)^{\mu_0+2}\left((\mu_0+1)!\right)^{-1}$ and
$a_i:=(\frac{\omega_i}{\omega_i(S)})^{2/\mu_0}$ for any positive integer $i$.
\end{props}
\begin{proof}

We have
$$
{\cal H}(\tilde f)=\left\{b\in {\cal K}\,|\,\forall i\geq 1,\forall \alpha\in {\cal W}^{2^{i-1}+1,2^i}_{n,*},
\left|\sum_{j=1}^l\tilde f_j(b)\alpha(g_j)\right|\geq \gamma\omega_i\right\}.
$$
In order to estimate $\text{mes}_{2n}{\cal H}(\tilde f)$ from below, it is sufficient to 
estimate \\$\text{mes}_{2n}{\cal K}-\text{mes}_{2n}{\cal H}(\tilde f)$ from above. In fact, we have 
$$
\text{mes}_{2n}{\cal H}(\tilde f)=\text{mes}_{2n}{\cal K}-\left(\text{mes}_{2n}{\cal K}-\text{mes}_{2n}{\cal H}(\tilde f)\right).
$$
Let $\alpha$ belong to ${\cal W}^{2^{i-1}+1,2^i}_{n,*}$, we set 
$$
(\tilde f(b),\alpha):= \sum_{j=1}^l\tilde f_j(b)\alpha(g_j)
$$
as well as
$$
|\alpha|_2:=\sqrt{|\alpha(g_1)|^2+\cdots+|\alpha(g_l)|^2}.
$$
Therefore, we have
\begin{eqnarray*}
{\cal K}\setminus {\cal H}(\tilde f) &= & \left\{b\in {\cal K}\,|\,\exists i\geq 1,\exists \alpha\in {\cal W}^{2^{i-1}+1,2^i}_{n,*}\text{ such that }\left|\sum_{j=1}^l\tilde f_j(b)\alpha(g_j)\right|< \gamma\omega_i\right\}\\
&= & \left\{b\in {\cal K}\,|\,\exists i\geq 1,\exists \alpha\in {\cal W}^{2^{i-1}+1,2^i}_{n,*}\text{ such that }\left|(\tilde f(b),\alpha)\right|< \gamma\omega_i\right\}\\
&=&\left\{b\in {\cal K}\,|\,\exists i\geq 1,\exists \alpha\in {\cal W}^{2^{i-1}+1,2^i}_{n,*}\text{ such that }\left|\left(\tilde f(b),\frac{\alpha}{|\alpha|_2}\right)\right|^2< \left(\frac{\gamma\omega_i}{|\alpha|_2}\right)^2\right\}.\\
\end{eqnarray*}
Let us set $c=\frac{1}{|\alpha|_2}(\alpha(g_1),\ldots,\alpha(g_l))$. It belongs to the unit sphere $\Bbb S^l$. Let us set $g(y)=|(f(y),c)|^2$ as well as $\tilde g(y)=|(\tilde f(y),c)|^2$. We have 
$$
\tilde g(y)-g(y)=|(\tilde f(y)-f(y),c)|^2+2Re\left((c,f(y))(c,\tilde f(y)-f(y))\right).
$$
By differentiation, we obtain for any nonnegative integer $\nu$ and for any $a$ in the unit sphere $\Bbb S^n$,
\begin{eqnarray*}
D^{\nu}(\tilde g-g)(y)(a^{\nu}) & = & \sum_{k=0}^{\nu}C^{\nu}_kD^{k}(\tilde f-f,c)(y)(a^k)D^{\nu-k}\overline{(\tilde f-f,c)}(y)(a^{\nu-k})\\
&&+2Re\left(\sum_{k=0}^{\nu}C^{\nu}_kD^{k}(f,c)(y)(a^k)D^{\nu-k}(\tilde f-f,c)(y)(a^{\nu-k})\right)\\
&=& \sum_{k=0}^{\nu}C^{\nu}_k\left(D^{k}(\tilde f-f)(y)(a^k),c\right)\left(D^{\nu-k}\overline{(\tilde f-f)}(y)(a^{\nu-k}),\bar c\right)\\
&&+2Re\left(\sum_{k=0}^{\nu}C^{\nu}_k\left(D^{k}f(y)(a^k),c\right)\left(D^{\nu-k}(\tilde f-f)(y)(a^{\nu-k}),c\right)\right).
\end{eqnarray*}
By Schwarz inequality, we obtain the following estimate for $\nu\leq \mu_0$:
\begin{eqnarray*}
|D^{\nu}(\tilde g-g)(y)(a^{\nu})|& \leq &\left(\sum_{k=0}^{\nu}C^{\nu}_k\right)\left(\max_{0\leq k\leq \nu}|D^{k}(\tilde f-f)(y)(a^k)|^2_2\right.\\
&&\left.+2\max_{0\leq k\leq \nu}|D^{k}(\tilde f-f)(y)(a^k)|_2\max_{0\leq k\leq \nu}|D^{k}(f)(y)(a^k)|_2\right)\\
& \leq & 2^{\nu}\left(\|\tilde f-f\|^2_{{\cal K}_{\vartheta},\nu}+\|\tilde f-f\|_{{\cal K}_{\vartheta},\nu}\|f\|_{{\cal K}_{\vartheta},\nu}\right)\\
& \leq & 2^{\mu_0}\left(\|\tilde f-f\|^2_{{\cal K}_{\vartheta},\mu_0}+\|\tilde f-f\|_{{\cal K}_{\vartheta},\mu_0}\|f\|_{{\cal K}_{\vartheta},\mu_0}\right).
\end{eqnarray*}
Therefore, we have 
$$
\|\tilde g-g\|_{{\cal K}_{\vartheta},\mu_0}\leq  2^{\mu_0}\left(\|\tilde f-f\|^2_{{\cal K}_{\vartheta},\mu_0}+\|\tilde f-f\|_{{\cal K}_{\vartheta},\mu_0}\|f\|_{{\cal K}_{\vartheta},\mu_0}\right).
$$
Let $P(X)=X^2+\|f\|_{{\cal K}_{\vartheta},\mu_0}X-\frac{\beta}{2^{\mu_0+1}}$
be a polynomial in the real indeterminate $X$. Its discriminant $\Delta=\|f\|_{{\cal K}_{\vartheta},\mu_0}^2+\frac{\beta}{2^{\mu_0-1}}$ is positive. Thus, $P$ has two real roots $r_{\pm}:=(-\|f\|_{{\cal K}_{\vartheta},\mu_0}\pm\sqrt{\Delta})/2$. 
The root $r_+$ is positive, whereas $r_-$ is negative. Clearly, $P$ is negative 
in $]r_-,r_+[$. As a consequence, 
$$
\|\tilde g-g\|_{{\cal K}_{\vartheta},\mu_0}\leq   2^{\mu_0}\left(\|\tilde f-f\|^2_{{\cal K}_{\vartheta},\mu_0}+\|\tilde f-f\|_{{\cal K}_{\vartheta},\mu_0}\|f\|_{{\cal K}_{\vartheta},\mu_0}\right)\leq \frac{\beta}{2},
$$
as soon as 
$$
\|\tilde f-f\|_{{\cal K}_{\vartheta},\mu_0} \leq \frac{-\|f\|_{{\cal K}_{\vartheta},\mu_0}+\sqrt{\|f\|_{{\cal K}_{\vartheta},\mu_0}^2+\frac{\beta}{2^{\mu_0-1}}}}{2}.
$$
In this case, we can apply R\"ussmann theorem \ref{mes-epsilon} to $\tilde g$: for any $\epsilon$ in $]0,\frac{\beta}{2\mu_0+2}]$, we have 
$$
\text{mes}_{2n}\{y\in {\cal K}\,|\,|(\tilde f(y),c)|^2\leq \epsilon\}\leq M\epsilon^{\frac{1}{\mu_0}}
$$
where 
\begin{eqnarray*}
M &= & Bd^{2n-1}\left(\frac{1}{\sqrt{2n}}+2d+\frac{d}{\vartheta}\right)\beta^{-1-\frac{1}{\mu_0}}\sup_{c\in \Bbb S^l}\||(f,c)|^2\|_{{\cal K}_{\vartheta},\mu_0+1}\\
B& =& 3(2\pi e)^{n}(\mu_0+1)^{\mu_0+2}[(\mu_0+1)!]^{-1}.
\end{eqnarray*}
We recall that 
\begin{eqnarray*}
{\cal K}\setminus {\cal H}(\tilde f) & = & \left\{b\in {\cal K}\,|\,\exists i\geq 1,\exists \alpha\in {\cal W}^{2^{i-1}+1,2^i}_{n,*}\text{ such that }\left|\left(\tilde f(b),\frac{\alpha}{|\alpha|_2}\right)\right|^2< \left(\frac{\gamma\omega_i}{|\alpha|_2}\right)^2\right\}\\
& = & \bigcup_{i\geq 1}\bigcup_{\alpha\in {\cal W}^{2^{i-1}+1,2^i}_{n,*}}\left\{b\in {\cal K}\,|\,\left|\left(\tilde f(b),\frac{\alpha}{|\alpha|_2}\right)\right|^2< \left(\frac{\gamma\omega_i}{|\alpha|_2}\right)^2\right\}.
\end{eqnarray*}
Let $\alpha$ be a weight in ${\cal W}^{2^{i-1}+1,2^i}_{n,*}$. Since, by definition, $\max_j|\alpha(g_j)|$ is greater than or equal to $\omega_i(S)$, we have
$$
|\alpha|_2=\sqrt{|\alpha(g_1)|^2+\cdots +|\alpha(g_l)|^2}\geq \omega_i(S).
$$
It follows that 
$$
\frac{\gamma\omega_i}{|\alpha|_2}\leq \frac{\gamma\omega_i}{\omega_i(S)}\leq \gamma M_{\omega,\omega(S) }.
$$
If 
$$
\gamma \leq \frac{1}{M_{\omega,\omega(S) }}\sqrt{\frac{\beta}{2\mu_0+2}},
$$
then 
$$
\left(\frac{\gamma\omega_i}{|\alpha|_2}\right)^2\leq \frac{\beta}{2\mu_0+2}.
$$
In this case, we obtain
$$
\text{mes}_{2n}\left\{b\in {\cal K}\,|\,\left|\left(\tilde f(b),\frac{\alpha}{|\alpha|_2}\right)\right|^2< \left(\frac{\gamma\omega_i}{|\alpha|_2}\right)^2\right\}\leq M\left(\frac{\gamma\omega_i}{|\alpha|_2}\right)^{\frac{2}{\mu_0}};
$$
so that 
$$
\text{mes}_{2n}{\cal K}\setminus {\cal H}(\tilde f)  \leq  M\sum_{i\geq 1}\sum_{\alpha\in {\cal W}^{2^{i-1}+1,2^i}_{n,*}}\left(\frac{\gamma\omega_i}{|\alpha|_2}\right)^{\frac{2}{\mu_0}}.
$$
Since ${\cal W}^{2^{i-1}+1,2^i}_{n,*}$ is isomorphic to a subset of $$\left\{(Q,j)\in {\bf N}^n\times \{1,\ldots,n\}\,|\,2^{i-1}+1\leq |Q|\leq 2^i\right\},$$
we have 
$$
\text{mes}_{2n}{\cal K}\setminus {\cal H}(\tilde f)  \leq  M\sum_{i\geq 1}\left(\sum_{\substack{Q\in{\bf N}^n\\ 2^{i-1}+1\leq |Q|\leq 2^i}}1\right)\left(\sum_{j=1}^n1\right)\left(\frac{\gamma\omega_i}{\omega_i(S)}\right)^{\frac{2}{\mu_0}}.
$$
Let us set 
$$
Z_i:= \sum_{\substack{Q\in{\bf N}^n\\ 0\leq |Q|\leq 2^i}}1.
$$ 
It is well known that $$Z_i=C^{n+2^i}_n=\frac{(n+2^i)(n-1+2^i)\cdots (2^i+1)}{n!}.$$ We have 
$$
\sum_{\substack{Q\in{\bf N}^n\\ 2^{i-1}+1\leq |Q|\leq 2^i}}1=Z_i-Z_{i-1}.
$$
The previous estimate can be written as 
$$
\text{mes}_{2n}{\cal K}\setminus {\cal H}(\tilde f)  \leq \gamma^{\frac{2}{\mu_0}}n M\sum_{i\geq 1}(Z_i-Z_{i-1})\left(\frac{\omega_i}{\omega_i(S)}\right)^{\frac{2}{\mu_0}}.
$$
Since $S$ is strictly diophantine relatively to $(\omega,\mu_0)$, we have
$$
\lim_{i\rightarrow +\infty}Z_{i}\left(\frac{\omega_i}{\omega_i(S)}\right)^{\frac{2}{\mu_0}}\leq \frac{1}{n!}\lim_{i\rightarrow +\infty}(n+2^i)^n\left(\frac{\omega_i}{\omega_i(S)}\right)^{\frac{2}{\mu_0}}=0.
$$
Let us set 
$$
a_i:=\left(\frac{\omega_i}{\omega_i(S)}\right)^{\frac{2}{\mu_0}}.
$$ 
Thus, we have
$$
\text{mes}_{2n}{\cal K}\setminus {\cal H}(\tilde f)  \leq \gamma^{\frac{2}{\mu_0}}n M\left[\sum_{i\geq 1}Z_i\left(a_i-a_{i+1}\right)-(n+1)a_1\right].
$$
Let $\psi:]1,+\infty[\rightarrow {\bf R}_+$ be the function defined to be 
$$
\psi(x):=\sum_{i\geq 1}\frac{\omega_i}{\omega_i(S)}\chi_{\scriptscriptstyle]2^{i-1},2^i]},
$$
where $\chi_{\scriptscriptstyle]2^{i-1},2^i]}$ denotes the characteristic function of $]2^{i-1},2^i]$. Since 
$$
\left(\frac{\omega_i}{\omega_i(S)}\right)^{\frac{2}{\mu_0}}-\left(\frac{\omega_{i+1
}}{\omega_{i+1}(S)}\right)^{\frac{2}{\mu_0}} = \int_{2^i}^{2^{i+1}}-d(\psi(t)^{\frac{2}{\mu_0}}),
$$
the previous estimate becomes
$$
\text{mes}_{2n}{\cal K}\setminus {\cal H}(\tilde f)  \leq \gamma^{\frac{2}{\mu_0}} \frac{M}{(n-1)!}\left[-\int_2^{+\infty}(2^t+n)^nd(\psi(t)^{\frac{2}{\mu_0}})-(n+1)a_1\right].
$$
Since 
$$
\lim_{t\rightarrow +\infty}(2^t+n)^n\psi(t)^{\frac{2}{\mu_0}}= \lim_{k\rightarrow +\infty}(2^k+n)^n\left(\frac{\omega_k}{\omega_k(S)}\right)^{\frac{2}{\mu_0}}=0, 
$$
After integrating by part, we obtain
\begin{eqnarray*}
\text{mes}_{2n}{\cal K}\setminus {\cal H}(\tilde f) &\leq &\gamma^{\frac{2}{\mu_0}} \frac{M}{(n-1)!}\left[(2^2+n)^na_2-(n+1)a_1\right.\\
&&\left.+n\ln 2\int_2^{+\infty}2^t(2^t+n)^{n-1}\psi(t)^{\frac{2}{\mu_0}}dt \right].
\end{eqnarray*}
By assumptions, for any positive integer $k$, $(2^k+n+1)^{n+1}a_k$ is less than or equal to $M_{\omega,\omega(S), 2/\mu_0}$. Thus, 
$$
\int_2^{+\infty}2^t(2^t+n)^{n-1}\psi(t)^{\frac{2}{\mu_0}}dt\leq M_{\omega,\omega(S),2/\mu_0}\int_2^{+\infty}\frac{dt}{2^t}=\frac{M_{\omega,\omega(S), 2/\mu_0}}{4\ln 2}.
$$
At the end, we obtain the estimate
$$
\text{mes}_{2n}{\cal K}\setminus {\cal H}(\tilde f)\leq \gamma^{\frac{2}{\mu_0}}\frac{M}{(n-1)!}\left[(2^2+n)^na_2-(n+1)a_1+\frac{n}{4}M_{\omega,\omega(S), 2/\mu_0}\right].
$$
Let $\epsilon^*$ be a postive number less than $\text{mes}_{2n}{\cal K}$. Let $\gamma$ denotes a positive number such that
\begin{eqnarray*}
\gamma & < &\max\left[ \left(\frac{\epsilon^*(n-1)!}{M\left[(2^2+n)^na_2-(n+1)a_1+\frac{n}{4}M_{\omega,\omega(S),2/\mu_0}\right]}\right)^{\frac{\mu_0}{2}}\right.,\\
& & \left. \frac{1}{M_{\omega,\omega(S) }}\sqrt{\frac{\beta}{2\mu_0+2}}\right].
\end{eqnarray*}
Then, the $2n$-measure $\text{mes}_{2n}{\cal K}\setminus {\cal H}(\tilde f)$ is less than or equal to $\epsilon^*$. Therefore, the $2n$-measure $\text{mes}_{2n}{\cal H}(\tilde f)$ is greater than or equal to $\text{mes}_{2n}{\cal K}-\epsilon^*$ and we are done.
\qed\end{proof}

\subsection{Application and proof of theorem \ref{theo3}}

The aim of this section is to give a sufficient condition which ensures that 
the sequence of compact sets $\{{\cal K}_k\}_{k\geq 1}$ will ``converge" to a
nonvoid compact set ${\cal K}_{\infty}$. More
precisely, we shall show that if $S$ is strictly diophantine with respect to $(\omega, \mu_0)$ and if the compact set ${\cal K}$ has a positive $2p$-measure, then ${\cal K}_{\infty}$ is  also of positive $2p$-measure.

Let us recall some facts: the vector field $\tilde X$ belongs to $\vfo n 1(D_n(0,1))$ and is assumed to be a good deformation of the nondegenerate vector field $X_0$. Let $\gamma$ be a positive number. Let ${\cal K}$ be a compact set of $\pi(D_0(\rho))$ of positive $2r$-measure. Let us consider the decreasing sequence $\{{\cal K}_k(NF,\omega,\gamma)\}_{k\geq k_0}$ of compact sets of $\pi(D_0(\rho))$ defined to be:
\begin{eqnarray*}
{\cal K}_{k_0} & = & {\cal K}\\
{\cal K}_k & = & \left\{b\in  {\cal K}_{k-1}\,|\,\forall \alpha\in {\cal W}_{n,*}^{2^{k}+1,2^{k+1}}\left|\alpha\left(\sum_{j=1}^l{a_j^{2^k}(b)g_j}\right)\right|\geq \gamma \omega_{k+1}\right\}.
\end{eqnarray*}
Here $$NF^m(x,u)=\sum_{j=1}^l a_j(u)^mS_j(x)$$ denotes the Lindstedt-Poincar\'e normal form of $\tilde X$ of order $m$.
We want to consider the map $a^{2^k}(b)$ as a perturbation of the nondegenerate map $a^{2^{k_0}}(b)$. First of all, we shall extend this function to a fixed neighborhood ${\cal B}$ of ${\cal K}$ (independent of $k$). Using a result of R\"ussmann, we can bound the norm of the extension $\tilde a^{2^k}(b)$ on ${\cal B}$ by a the norm of $a^{2^k}(b)$ on ${\cal K}+t_{2^{k+1}}$, the $t_{2^{k+1}}$-neighborhood of ${\cal K}$. Let $\mu_0$ be the amount of nondegeneracy of $a^{2^{k_0}}$. We shall estimate the $C^{\mu_0}$-norm of the difference on some well chosen neighborhood of ${\cal K}_{k-1}$. We shall show that this estimate is small enough so that we can apply proposition \ref{prop-mesure}. It will follow that the set ${\cal H}(\tilde a^{2^k})$ has positive measure and is contained in ${\cal K}_k$.

\begin{theos}[theorem 19.7 \cite{russmann-weak}]\label{russmann2}
Let ${\cal K}$ be a nonvoid set in ${\bf C}^n$, let $t$ be  positive number and let ${\cal K}+t$ denote the $t$-neighborhood of ${\cal K}$. Let $f:{\cal K}+t\rightarrow {\bf C}^p$ be a holomorphic and bounded function. Then, there is a $C^{\infty}$ function $\tilde f : {\bf C}^n\rightarrow {\bf C}^p$ such that $\tilde f(x)=f(x)$ for all $x\in {\cal K}$, and the estimates 
$$
\sup_{x\in {\bf C}^n}\sup_{a\in {\bf C}^n, |a|_2\leq 1}\left|D^{\nu}\tilde f(x)(a^{\nu})\right|\leq C(n,\nu)t^{-\nu}\sup_{x\in {\cal K}+t}|f(x)|,\quad \nu\in {\bf N}
$$
hold with constants $C(n,\nu)$ not depending on $f$ and satisfying the inequalities
$$
1=C(n,0)\leq C(n,1)\leq C(n,2)\leq \cdots.
$$
\end{theos}

From now on, we shall define $\tilde a^{2^k}$ to be the aforementioned extension of ${a^{2^k}}$ to ${\bf C}^n$.

\begin{lemms}
Let $\chi$ be a function in $C^{\mu_0}({\cal K}_{\vartheta}, {\bf C}^l)$ such that, for any integer $k_0\leq k\leq \nu$,
$$
\left\|\chi-\tilde a^{2^k}\right\|_{{\cal K}_{\vartheta}}\leq \frac{\gamma\omega_{k+1}}{l\Lambda(2^{k+1}+1)}.
$$
Assume that none of the ${\cal K}_k$'s is nonvoid when the interger $k$ ranges from $k_0$ to $\nu$. Then, ${\cal K}_{\nu+1}$ contains $\cap_{k=k_0+1}^{\nu}{\cal H}_{k}(\chi)$ where 
$$
{\cal H}_{k}(\chi):=\left\{b\in  {\cal K}\,|\,\forall \alpha\in {\cal 
W}_{n,*}^{2^{k}+1,2^{k+1}},\quad\left|\alpha\left(\sum_{j=1}^l{\chi_j(b)g_j}\right)\right|\geq 2\gamma \omega_{k+1}\right\}.
$$
\end{lemms}
\begin{proof}
Let $b$ be a point in ${\cal H}_{k}(\chi)$, for some $k_0\leq k\leq \nu$. Let $\alpha$ be a weight in ${\cal W}_{n,*}^{2^{k}+1,2^{k+1}}$. We have
\begin{eqnarray*}
\left|\alpha\left(\sum_{j=1}^l\tilde a^{2^k}_j(b)g_j\right)\right| & \geq & \left|\left|\alpha\left(\sum_{j=1}^l\chi_j(b)g_j\right)\right|-\left|\alpha\left(\sum_{j=1}^l(\tilde a^{2^k}_j-\chi_j)(b)g_j\right)\right|\right|\\
& \geq & \left| 2\gamma \omega_{k+1} - \left|\alpha\left(\sum_{j=1}^l(\tilde a^{2^k}_j-\chi_j)(b)g_j\right)\right|\right|.
\end{eqnarray*}
Moreover, for all $1\leq j\leq l$, $|\alpha(g_j)|\leq \Lambda(2^{k+1}+1)$. So, we obtain  
$$
\left|\alpha\left(\sum_{j=1}^l{\left(\tilde a^{2^k}_j-\chi_j\right)(b)g_j}\right)\right|\leq l\Lambda(2^{k+1}+1)\left\|\chi-\tilde a^{2^k}\right\|_{{\cal K}}.
$$
Since $\left\|\chi-\tilde a^{2^k}\right\|_{{\cal K}_{\vartheta}}$ is less than or equal to $\frac{\gamma\omega_{k+1}}{l\Lambda(2^{k+1}+1)}$, we obtain 
$$
\left|\alpha\left(\sum_{j=1}^l{\tilde a^{2^k}_j(b)g_j}\right)\right| \geq \gamma\omega_{k+1}.
$$
Thus, we have proved the inclusion
$$
{\cal H}_{k}(\chi)\subset \left\{b\in  {\cal K}\,|\,\forall \alpha\in {\cal W}_{n,*}^{2^k+1,2^{k+1}},\quad\left|\alpha\left(\sum_{j=1}^l\tilde a_j^{2^k}(b)g_j\right)\right|\geq \gamma \omega_{k+1}\right\}.
$$
If $k$ is greater than $k_0$, we have
\begin{eqnarray*}
{\cal K}_{k} & = & {\cal K}_{k-1}\cap \left\{b\in  {\cal K}\,|\,\forall \alpha\in {\cal W}_{n,*}^{2^k+1,2^{k+1}},\;\left|\alpha\left(\sum_{j=1}^l\tilde a_j^{2^k}(b)g_j\right)\right|\geq \gamma \omega_{k+1}\right\}\\
&\supset & {\cal K}_{k-1}\cap {\cal H}_{k}(\chi).
\end{eqnarray*}
By induction on the integer $k$ which is greater than or equal to $k_0$, we obtain 
$$
{\cal K}_{k}\supset {\cal K}_{k_0}\cap \bigcap_{j=k_0+1}^k {\cal H}_{j}(\chi)=\bigcap_{j=k_0+1}^k {\cal H}_{j}(\chi)\supset \bigcap_{j=1}^k {\cal H}_{j}(\chi).
$$
\qed\end{proof}

We want to apply the previous lemma to $\tilde a^{2^{\nu}}$. Let us define ${\cal B}_{k}:={\cal K}_k+t_{2^{k}}$ to be the $t_{2^{k}}$-neighborhood of ${\cal K}_k$. The inclusion ${\cal B}_{\nu}\subset{\cal B}_{k}$ holds whenever $k$ is less than or equal to $\nu$ and greater than or equal to $k_0$.
For any integer $k_0\leq k\leq \nu$, we have
\begin{eqnarray*}
\left\|\tilde a^{2^{\nu}}-\tilde a^{2^k}\right\|_{{\cal K}_{\vartheta}}& \leq & \sum_{j=k+1}^{\nu}\left\|\tilde a^{2^{j}}-\tilde a^{2^{j-1}}\right\|_{{\cal K}_{\vartheta}}\\
& \leq & \sum_{j=k+1}^{\nu}\left\|a^{2^{j}}-a^{2^{j-1}}\right\|_{{\cal B}_j}\\
& \leq & \sum_{j=k+1}^{\nu}2l|L^{-1}|\left|NF^{2^{j}}-NF^{2^{j-1}}\right|_{{\cal B}_j, R_j}\\
& \leq & 2l|L^{-1}|\sum_{j=k+1}^{\nu}\frac{\gamma^6\omega_{j}^6}{c_1^32^{6(j-1)}}\\
&\leq &  \frac{2l|L^{-1}|\gamma^6\omega_{k+1}^6}{c_1^3}\frac{2^{-6k}-2^{-6\nu}}{1-2^{-6}}.\\
\end{eqnarray*}
The third inequality is due to lemma \ref{a-nf} while the fourth one is due to inequality $(\ref{maj-b})$.
If $k$ is large enough, then 
$$
\frac{2l|L^{-1}|l\Lambda\gamma^5\omega_{k+1}^5}{c_1^3}\frac{2^{-6k}-2^{-6\nu}}{1-2^{-6}}(2^{k+1}+1)\leq 1.
$$ 
As a consequence, if $k_0$ is large enough, then for all integer $\nu$ greater than or equal to $k_0$,
$$
{\cal K}_{\nu+1}\supset \bigcap_{k=k_0+1}^\nu{\cal H}_{k}(\tilde a^{2^{\nu}})\supset \bigcap_{k=k_0+1}^{+\infty}{\cal H}_{k}(\tilde a^{2^{\nu}}).
$$

Let us show that the set $$\bigcap_{k=k_0+1}^{+\infty}{\cal H}_{k}(\tilde a^{2^{\nu}})$$
is of positive $2p$-measure. In order to do this, we shall apply proposition \ref{prop-mesure} to $\tilde a^{2^{\nu}}$. We have to obtain good estimates for the derivatives of the approximate function. In order to compensate the power of $t$ which arises in the inequality (cf. theorem \ref{russmann2}), we shall ``decrease" the radius on which we have obtained the estimates which led to the proof of the existence of invariant analytic subsets.

Let $\mu$ be an integer greater than or equal to $3$. Let us define the sequence of positive numbers $$1/2<\tilde R_0=r\leq 1, \quad\tilde R_{j+1}=\theta_j^{\mu}\tilde R_j.$$
As in lemma \ref{rayon-conv}, we show that if $k$ is greater than or equal to some positive integer $k_2$, then $\tilde R_k> \tilde R_{k_2}/2$. Let us set $R_{k_2}:=1$ as above. The basic sets on which the estimates are done are now the $D_n(0,\tilde R_{2^k})\times D_p(b,t_m)$'s. The convergence will follow from the analysis done above. Using the notation of the induction process section, we have the estimate (following (\ref{maj-b}))
\begin{eqnarray}
|B'|_{D_p(b,t_{2m}),\tilde R_{k+1}} & \leq & |B+C|_{D_p(b,t_{2m}),\tilde R_{k+1}}\nonumber\\
& \leq & \theta_k^{\mu(m+1)}|B+C|_{D_p(b,t_{m}),r}\nonumber\\
& \leq & \frac{\gamma^{2\mu}\omega_{k+1}^{2\mu}}{c_1^{\mu}m^{2\mu}}|B+C|_{D_p(b,t_{m}),r}\nonumber\\
& \leq & \frac{2^5n\gamma^{2\mu}\omega_{k+1}^{2\mu}}{c_1^{\mu}m^{2\mu+4}}.\label{maj-b2}
\end{eqnarray}
The second inequality is due to the fact that $B+C$ is of order greater than or equal to $m+1$.
Therefore, we obtain, for any nonnegative integer $k$ less than or equal to $\mu_0$,
\begin{eqnarray*}
\left\|D^k(\tilde a^{2^{\nu}}-\tilde a^{2^{k_0}})\right\|_{{\cal K}_{\vartheta}}& \leq & \sum_{j=k_0+1}^{\nu}\left\|D^k(\tilde a^{2^{j}}-\tilde a^{2^{j-1}})\right\|_{{\cal K}_{\vartheta}}\\
& \leq & \sum_{j=k_0+1}^{\nu}\left\|D^k(\tilde a^{2^{j}}-\tilde a^{2^{j-1}})\right\|_{{\cal K}_{\vartheta}}\\
& \leq & \sum_{j=k_0+1}^{\nu}C(p,k)t_{2^j}^{-k}\left\|a^{2^{j}}- a^{2^{j-1}}\right\|_{{\cal B}_j}\\
& \leq & \sum_{j=k_0+1}^{\nu}2l|L^{-1}|C(p,k)t_j^{-k}\left|NF^{2^{j}}-NF^{2^{j-1}}\right|_{{\cal B}_j, \tilde R_j}\\
& \leq & C(p,k)2l|L^{-1}|\frac{2^kl^k\Lambda^k 2^5n\gamma^{2\mu}}{c_1^{\mu}\gamma^k }\sum_{j=k_0+1}^{\nu}\frac{(2^{j+1}+1)^k\omega_{j+1}^{2\mu-k}}{2^{j2\mu}}.\\
\end{eqnarray*}
The last inequality is due to inequality $(\ref{maj-b2})$. Moreover, the last sum is not only convergent but also small if $\mu$ is well chosen (with respect to $\mu_0$). As a consequence, we obtain 
$$
\|a^{2^{k_0}}-\tilde a^{2^{\nu}}\|_{{\cal K}_{\vartheta}, \mu_0}\leq \frac{-\|a^{2^{k_0}}\|_{{\cal K}_{\vartheta},\mu_0}+\sqrt{\|a^{2^{k_0}}\|_{{\cal K}_{\vartheta},\mu_0}^2+\frac{\beta}{2^{\mu_0-1}}}}{2}.
$$ 
Hence, according to proposition \ref{prop-mesure}, the measure of the set 
$$
{\cal H}(\tilde a^{2^{\nu}}):=\left\{b\in {\cal K}\,|\,\forall i\in {\bf N}^*, \quad\forall \alpha\in {\cal W}_{n,*}^{2^{i-1}+1,2^{i}}\left|\alpha\left(\sum_{j=1}^l{\tilde a^{2^{\nu}}_j(b)g_j}\right)\right|\geq \gamma \omega_i\right\}
$$
satisfies $\text{mes}_{2p}{\cal H}(\tilde a^{2^{\nu}})\geq mes_{2n}{\cal
  K}-\epsilon^*$.
Therefore, ${\cal K}_{\nu+1}$ is nonvoid, its $2p$-measure satisfies 
$$
\text{mes}_{2p}{\cal H}(\tilde a^{2^{\nu}})\geq mes_{2p}{\cal K}-\epsilon^*,$$
and we are done.

\section{Where do the tori of the classical KAM theorem come from~?}

Let us give a taste of how we can recover the ``classical" KAM theory with
genuine real tori. Let us consider a real analytic hamiltonian $H$ in a neighborhood of the origin in ${\bf R}^{2n}$. We assume that 
$$
H(x,y)=\sum_{l=1}^{N_0}\sum_{i=1}^n\mu_{i,l}(x_i^2+y_i^2)^l+h_{M_0+1}(x,y),
$$
where $h_{M_0+1}(x,y)$ is a real analytic function of order greater than or
equal to $M_0+1\geq 2N_0+1$. Here, the $\mu_{i,l}$'s are real numbers and $\omega=\sum_{i=1}^ndx_i\wedge dy_i$ denotes the canonical symplectic form of ${\bf R}^{2n}$.
Let us write the hamiltonian the the complex coordinates $z_j=x_j+iy_j$, $j=1,\ldots, n$. We have
$$
H(x,y)=\tilde H(z,\bar z)=\sum_{l=1}^{N_0}\sum_{i=1}^n\mu_{i,l}(z_i\bar z_i)^l+\tilde h_{M_0+1}(z,\bar z),
$$
where $\tilde h_{M_0+1}(z,\bar z)=h_{M_0+1}(x,y)$.

Let us complexify the hamiltonian. We obtain a holomorphic hamiltonian $G$ in a
neighborhood of the origin in ${\bf C}^{2n}$ with $(z,w)$ as complex symplectic coordinates:
$$
G(z,w)=\sum_{l=1}^{N_0}\sum_{i=1}^n\mu_{i,l}(z_iw_i)^l+\tilde h_{M_0+1}(z,w).
$$
We recover $\tilde H$ (or $H$) by restricting $G$ to the set
$$
\bigcap_{i=1}^n\{w_i=\bar z_i\}.
$$
It is assumed that $G$ is a perturbation of the nondegenerate integrable hamiltonian 
$$
H_0=\sum_{l=1}^{N_0}\sum_{i=1}^n\mu_{i,l}(z_iw_i)^l.
$$
Let $\lie g$ be an $n$-dimensional commutative lie algebra. Let $S$ be the injective semi-simple linear morphism defined to be: $$S(g_i)=z_i\frac{\partial}{\partial z_i}-w_i\frac{\partial}{\partial w_i},\quad i=1,\ldots,n.$$ 
Its nonzero weights are integers so that $S$ is diophantine. Its ring of invariant is $\widehat{\cal O}_{2n}^{S}=\Bbb C[[u_1,\ldots,u_n]]$ with $u_i=z_iw_i$ and its centralizer $\left(\fvfo {2n} 1\right)^{S}$ is the $\Bbb C[[u_1,\ldots,u_n]]$-module generated by $z_i\frac{\partial}{\partial z_i}$ and $w_i\frac{\partial}{\partial z_i}$ with $1\leq i\leq n$. We refer to our previous \cite{stolo-ihes}[chapter 10] for more details.
Since the vector field $X_G$ associated to $G$ is symplectic, its Lindstedt-Poincar\'e normal form of order any order $m$ is of the form:
$$
\sum_{j=1}^na_j^m(u_1,\ldots, u_n)S(g_j).
$$
Therefore, we can apply our result: $X_G$ has invariant analytic subsets which are biholomorphic to
the intersection of $\cap_{i=1}^n\{w_iz_i=c_i\}$ with a fixed polydisc, for some well chosen constants $c_i$. 
The normalization process is certainly compatible with the restriction to the
set $\cap_{i=1}^n\{w_i=\bar z_i\}$; that is, it commutes with the complex
conjugacy. This can be done as in the case of Poincar\'e-Dulac normal form
(see \cite{bruno-hamilton} for instance). Therefore, if one chooses a set of
real constants, the hamiltonian vector field will have invariant (real) analytic subsets analytically isomorphic to the intersection of a fixed polydisc with
$$
\left(\bigcap_{i=1}^n\{w_iz_i=c_i\}\right)\bigcap \left(\bigcap_{i=1}^n\{w_i=\bar z_i\}\right)=\bigcap_{i=1}^n\{z_i\bar z_i=x_i^2+y^2_i=c_i\}.
$$
for some real constants $c_i$. These are the genuine real tori.

\subsection{The volume preserving case}
Let us consider a holomorphic volume preserving vector field $X$ which is a deformation of a nondegenerate volume preserving polynomial vector field $X_0$ in a neighborhood of the origin of ${\bf C^n}$. Let $\lie g$ be a $(n-1)$-dimensional commutative Lie algebra with a basis $G=\{g_1,\ldots,g_{n-1}\}$. Let $S$ be the linear semi-simple and injective morphism defined to be
$$
S(g_i)=x_i\frac{\partial}{\partial x_i}-x_{i+1}\frac{\partial}{\partial x_{i+1}},\quad i=1,\ldots,n-1.
$$
The values of the nonzero weights of $S$ are integers; thus, $S$ is diophantine. Moreover, if we set $u=x_1\cdots x_n$, then the ring of invariant of $S$ is defined to be 
$\widehat{\cal O}_n^{S}=\Bbb C[[u]]$ whereas its centralizer $\left(\fvfo n 1\right)^{S}$ is the $\Bbb C[[u]]$-module generated by the $x_i\frac{\partial}{\partial x_i}$'s, $1\leq i\leq n$. 
Since the vector field $X$ is volume preserving, its Lindstedt-Poincar\'e normal form of order any order $m$ is of the form:
$$
\sum_{j=1}^{n-1}a_j^m(u)S(g_j).
$$
Therefore, we can apply our result: $X$ has invariant analytic subsets which are biholomorphic to to the intersection of a fixed polydisc with
$$\left\{x_1\cdots x_n=c_i\right\}$$ for some well chosen constants $c_i$. 

\bibliographystyle{alpha}

\bibliography{normal,math,asympt,analyse,stolo,lie,kam}

\end{document}